\documentclass[a4paper,12pt]{article}
\usepackage[utf8]{inputenc}
\usepackage[T2A]{fontenc}
\usepackage[english]{babel}
\usepackage{amsthm,amsmath,amssymb,latexsym}
\newcommand{\ld}{\begin{smallmatrix}.\;\;\\\;\;.\end{smallmatrix}}

\begin{document}

{\begin{center}
\textbf{A.A. Tuganbaev}\\
\textbf{Centrally Essential Rings and Semirings}\\
\end{center}}

This work is a review of results about centrally essential rings and semirings. A ring (resp., semiring) is said to be centrally essential if it is either commutative or satisfy the property that for any non-central element $a$, there exist non-zero central elements $x$ and $y$ with $ax=y$. The class of centrally essential rings is very large; many corresponding examples are given in the work.\\ The author is very grateful for the great help in editing the manuscript to Adel Abyzov, Oleg Lyubimtsev and Danil Tapkin.

\tableofcontents

\addcontentsline{toc}{section}{Introduction}

\section*{Introduction}

In Introduction and Sections 1--5, the word a \textsf{ring} means an \textsf{associative ring}. By default, it is assumed that a ring has $1\ne 0$; the case of not necessarily unital rings is specially indicated.

In Section 6, the word a \textsf{ring} means a \textsf{not necessarily associative ring}.
 
A not necessarily unital ring $A$ is said to be a \textsf{ring with essential center} or a \textsf{centrally essential}\label{ceness} ring, or a \textsf{CE ring}\label{cering} if either $A$ is commutative or for any non-central element $a\in A$, there are non-zero central elements $x$ and $y$ with $ax=y$. 

It is clear that any commutative ring is centrally essential. An unital ring $A$ with center $C=Z(A)$ is centrally essential if and only if the $C$-module $A$ is an essential extension of $C_C$.

From the definition of a centrally essential ring $A$, it might seem that such a ring is possibly commutative. Indeed, $A$ satisfies many propeties of commutative rings. For example,
\begin{itemize}
\item
all idempotents of $A$ are central, see 1.1.4 below;
\item 
if $A$ is semiprime, then $A$ is commutative, see 1.2.2;
\item 
If $A$ is a centrally essential local ring, then the ring $A/J(A)$ is a field and, in particular, is commutative; see 1.3.2.
\item 
If $A$ is a right or left semi-Artinian centrally essential ring, then the factor ring $A/J(A)$ is commutative; see 1.4.5.
\end{itemize}

However, a centrally essential ring $A$ may be very far from a commutative ring. For example,
\begin{itemize}
\item
the factor ring $A/J(A)$ of $A$ with respect to the prime radical can be not centrally essential and, in particular, the semiprime ring $A/J(A)$ can be non-commutative; see 3.5.5;
\item 
for any ideal of $A$ generated by central idempotents, the corresponding factor ring is not necessarily centrally essential; see 2.2.5;
\item 
factor rings of the ring $A$ are not necessarily centrally essential; see two previous items;
\item 
there are finite non-commutative centrally essential unital the group algebras; see Example 1 below;
\item 
there are finite non-commutative centrally essential Grassmann algebras; see Example 2 below;
\item 
there are torsion-free Abelian groups $G$ of finite rank such that their endomor\-phism rings are non-commutative centrally essential rings; see 3.7.13(c).
\end{itemize}

\textbf{Example 1.} Let $F$ be the field of order $2$ and let $G=Q_8$ be the quaternion the group of order 8, i.e., $G$ is the group with two generators $a,b$ and three defining relations $a^4=1$, $a^2=b^2$ and $aba^{-1}=b^{-1}$; see \cite[Section 4.4]{Hal59}. Then the group algebra $FG$ is a non-commutative finite local centrally essential ring consisting of $256$ elements; this follows from Proposition 3.2.4 below.

We give some necessary notions. For a ring $A$, we denote by $Z(A)$ (or $C(A)$), $J(A)$, $N(A)$ and $K(A)$ the center, the Jacobson radical, the prime radical and the K\"othe radical (i.e., the sum of all nil-ideals which is the largest nil-ideal), respectively. We also set $[a,b]=ab-ba$ for any two elements $a,b$ of the ring $A$. For a group or a semigroup $X$, we denote by $Z(X)$ or $C(X)$ the center of $X$.

\textbf{Example 2.} We give one more example of a non-commutative finite centrally essential ring. Let $F$ be a field consisting of three elements, $V$ be a linear $F$-space with basis $e_1,e_2,e_3$, and let $\Lambda(V)$ be the Grassmann algebra \footnote{See 2.2.1.} of the space $V$. Since $e_1\wedge e_1=e_2\wedge e_2=e_3\wedge e_3=0$ and any product of generators is equal to the $\pm$product of generators with ascending subscripts, we have that $\Lambda(V)$ is a finite $F$-algebra of dimension 8 with basis
$$
\{1,e_1,e_2,e_3,e_1\wedge e_2,e_1\wedge e_3,e_2\wedge e_3,e_1\wedge e_2\wedge e_3\},
$$
$$
|\Lambda(V)|=3^8,\;e_k\wedge e_i\wedge e_j=-e_i\wedge e_k\wedge e_j=e_i\wedge e_j\wedge e_k.
$$
Therefore, if
$$
x=\alpha_0\cdot 1+\alpha_1^1e_1+\alpha_1^2e_2+\alpha_1^3e_3+\alpha_2^1e_1\wedge e_2+\alpha_2^2e_1\wedge e_3+$$
$$
+\alpha_2^3e_2\wedge e_3+\alpha_3e_1\wedge e_2\wedge e_3,
$$ 
then 
$$
\begin{array}{lll}[e_1,x]&=&2\alpha_1^2e_1\wedge e_2+2\alpha_1^3e_1\wedge e_3,\\
{}[e_2,x]&=-&2\alpha_1^1e_1\wedge e_2+2\alpha_1^3e_2\wedge e_3,\\
{}[e_3,x]&=-&2\alpha_1^1e_1\wedge e_3-2\alpha_1^2e_2\wedge e_3.\end{array}
$$
Thus, $x\in Z(\Lambda(V))$ if and only if $\alpha_1^1=\alpha_1^2=\alpha_1^3=0.$ In other words, the center of the algebra $\Lambda(V)$ is of dimension $5$. On the other hand, if $\alpha_1^1\neq 0$, then
$$
x\wedge(e_2\wedge e_3)=\alpha_0e_2\wedge e_3+\alpha_1^1e_1\wedge e_2\wedge e_3\in Z(\Lambda(V))\setminus \{0\}.
$$ 
In addition, $e_2\wedge e_3\in Z(\Lambda(V))$. A similar argument applies if $\alpha_1^2\neq 0$ or $\alpha_1^3\neq 0$. Consequently, $\Lambda(V)$ is a finite centrally essential non-commutative ring.

\textbf{Example 3.} This example is due to the reviewer of the paper \cite{MT19b} who has kindly suggested a series of examples of non-commutative centrally essential rings arising from a construction described in \cite{Jel16}. The following statement $(*)$ and its proof were also proposed by the reviewer:

$(*)$ If $B$ is an ideal of the ring $A$ such that $B\subseteq Z(A)$ and $A/B$ is a field, then $A$ is a centrally essential ring.

$\lhd$ We assume that $A$ is not commutative and $a$ is a non-central element of $A$. If $aB\neq 0$ then it is clear that $Z(A)\cap aZ(A)\neq 0$. We assume the contrary, i.e., $aB=0$. Since $a\notin B$ and $A/B$ is a field, the element $a$ is invertible modulo $B$, i.e., $sa=1-x$ for some $s\in A$ and $x\in B$. For any $y\in B$, we have $0=say=y-xy$; this implies that $xB=B=xA$, $x$ is a central idempotent, and $A$ has a Peirce decomposition $A=Ax\oplus A(1-x)$ where the both summands $Ax=B$ and $A(1-x)\cong A/B$ are commutative. Therefore, $A$ is commutative. This contradicts to the choice of $A$, and $(*)$ is true. 

It remains to consider the simplest case of the construction given in \cite[Proposition 7]{Jel16} (we reserve notation of this paper). Let $F=\mathbb{Q}(x,y)$ be the field of rational functions. We consider two partial derivations $d_1=\dfrac{\partial}{\partial x}$ and $d_2=\dfrac{\partial}{\partial x}$. Then the ring $A=T(F,F)$ consisting of matrices
$$
\left\lbrace\left.\begin{pmatrix}
f&d_1(f)&g\\
0&f&d_2(f)\\
0&0&f
\end{pmatrix}\;\;\right|\;\; f,g\in F \right\rbrace
$$
and its ideal
$$
B=\hat F=\left\lbrace\left.\begin{pmatrix}
0&0&g\\
0&0&0\\
0&0&0
\end{pmatrix}\;\;\right|\;\; g\in F \right\rbrace
$$
satisfy conditions of $(*)$.~$\rhd$

We give some definitions.
For a module $M$, the \textsf{socle}\label{soc} $\text{Soc}\,M$ is the sum of all simple submodules of $M$; if $M$ does not contain simple submodules, then $\text{Soc}\,M = 0$ by definition. A module $M$ is said to be \textsf{finite-dimensional}\label{findimmod} (in the sense of Goldie) if $M$ does not contain submodule which is an infinite direct sum of non-zero submodules. A module $M$ is said to be \textsf{Noetherian}\label{netermod} (resp., \textsf{Artinian})\label{artinmod} if $M$ does not contain an infinite properly ascending (resp., properly descending) chain submodules. Direct summands of free modules are called \textsf{projective}\label{promod} modules. A module $M$ is said to be \textsf{hereditary}\label{hermod} if all submodules of the module $M$ are projective. 
A module $M$ is said to be \textsf{distributive}\label{dismod} (resp., \textsf{uniserial}),\label{unisermod}
 if the submodule lattice of the module $M$ is distributive (resp., it is a chain). 
We recall that a module $X$ is called an \textsf{essential extension}\label{essext} submodule $Y$ of the module $X$ if $Y\cap Z\ne 0$ for any non-zero submodule $Z$ in $X$. In this case, $Y$ is called an \textsf{essential submodule}\label{esssub} of the module $X$. A submodule $Y$ of the module $X$ is said to be \textsf{closed}\label{closed} (in $X$) if $Y=Y'$ for any submodule $Y'$ of the module $X$ which is an essential extension of the module $Y$.

A ring $A$ is called a \textsf{domain}\label{dom} if $A$ does not have non-zero zero-divisors. A commutative domain $A$ is called a \textsf{Dedekind} domain if $A$ is a commutative hereditary Noetherian domain. If $A$ is a ring, then a proper ideal $B$ of the ring $A$ is said to be \textsf{completely prime} if the factor ring $A/B$ is a domain. A ring $A$ is said to be \textsf{right invariant}\label{invRing} (resp., \textsf{left invariant}) if all right (resp., left) ideals of the ring $A$ are ideals. A ring $R$ is said to be \textsf{semiprime}\label{semipr} (resp., \textsf{prime}),\label{priring} if $R$ does not have nilpotent non-zero ideals (resp., the product of any two non-zero ideals of the ring $R$ is not equal to zero).
A ring $R$ is said to be \textsf{arithmetical}\label{ariring} if the lattice of its two-sided ideals is distributive, i.e., $X\cap (Y+Z)=X\cap Y+X\cap Z$ for any three ideals $X,Y,Z$ of the ring $R$. 
It is clear that a commutative ring is right (resp., left) distributive if and only if the ring is arithmetical.
For a ring $R$, an element $r$ is called a \textsf{left non-zero-divisor}\label{nzd} or a \textsf{right regular}\label{regele} element if the relation $rx = 0$ implies the relation $x = 0$ for any $x\in R$.
We note that one-sided zero-divisors are two-sided zero-divisors in a centrally essential ring; see 1.1.2(a).
A ring $R$ has the \textsf{right (resp., left) classical ring of fractions}\label{clfra} $Q_{\text{cl}}(R_r)$ (resp., $Q_{\text{cl}}(R_l)$) if and only if for any two elements $a, b\in R$ such that $b$ is a non-zero-divisor, there exist elements $c, d\in R$ such that $d$ is a non-zero-divisor and $bc = ad$ (resp., $cb = da$). If the rings $Q_{\text{cl}}(R_r)$ and $Q_{\text{cl}}(R_l)$ exist, then they are isomorphic to each other over $R$. In this case, one says that there exists the two-sided ring of fractions $Q_{\text{cl}}(R)$.

For a ring $R$ and a subset $S$ in $R$, we denote by $\ell_R(S)$ the \textsf{left annihilator}\label{ann} $\{r\in R\,|\, rS=0\}$ of the set $S$. The right annihilator $\text{r}_R(S)$ is defined similarly. For a right (resp., left) $R$-module $M$, its fully invariant submodule consisting of all elements whose annihilators are essential right (resp., left) ideals in $R$, is called the \textsf{singular submodule} for $M$; it is denoted by $\text{Sing} M$.\label{sinsub} For $M=R_R$ (resp., $M={}_RR$), the ideal $\text{Sing} M$ is called the \textsf{right} (resp., \textsf{left}) \textsf{singular ideal}\label{sinide} of the ring $R$.

Nesessary information from ring theory is contained in \cite{Row88}, \cite{Bou89}, \cite{Lam09}, \cite{Tug98}, \cite{Her05}, \cite{Lam01}. See \cite{Fuc15} and \cite{KMT03} for necessary information on Abelian groups.

\section[Semiprime, Local, Perfect\\ and Semi-Artinian Rings]{Semiprime, Local, Perfect and Semi-Artinian Rings}\label{section1}

In Section 1, the word a \textsf{ring} means an \textsf{associative ring}. By default, it is assumed that the ring has a non-zero identity element; the case of not necessarily unital rings is specified separately.

\subsection{General Properties}\label{subsection1.1}

\textbf{1.1.1. Remark.}\label{zd-ideal}\\
If $A$ is a ring such that the set $B$ of all left zero-divisors is an ideal, then $B$ is a completely prime ideal.

$\lhd$ Let $a,b\in A$ and $ab\in B$. Then there exists an element $x\in A\setminus\{0\}$ such that $abx=0$. If $bx=0$, then $b\in B$. Otherwise, it follows from relation $a(bx)=0$ that $a\in B$.~$\rhd$

\textbf{1.1.2. Non-zero-divisors in centrally essential rings.} \\
Let $A$ be a centrally essential ring.

\textbf{a.} Every left (resp., right) non-zero-divisor $a$ of the ring $A$ is a right (resp., left) non-zero-divisor of the ring $A$.

\textbf{b.} The ring $A$ is left uniform\footnote{A module $M$ is said to be \textsf{uniform}\label{unifo} if any two its non-zero submodule have non-zero intersection.} if and only if $A$ is right uniform.

\textbf{c.} If the ring $A$ is right uniform and $B=\text{Sing}\,A_A$, then $B$ is the set of all (left or right) zero-divisors of the ring $A$ and $B$ is a completely prime ideal of the ring $A$.

\textbf{d.} If the ring $A$ has a proper ideal $B$ containing all left zero-divisors of the ring $A$, then the factor ring $A/B$ is commutative.

\textbf{e.} If an ideal $B$ of the ring $A$ contains all central zero-divisors of the ring $A$, then $\ell.Ann_A(B)\subseteq{}Z(A)$.

$\lhd$ \textbf{a.} We consider only the case, where $a$ is a left non-zero-divisor. We can assume that $a$ is a central element of the ring $A$. We assume the contrary. Then $ba=0$ for some non-zero element $b$ of the ring $A$. Since $b\ne 0$, there exist non-zero central elements $x,y$ of the ring $A$ such that $bx=y\ne 0$. Then $ya=bxa=bax=0$. This is a contradiction.

\textbf{b.} We assume that the ring $A$ is right uniform and $a_1,a_2$ are non-zero elements of the ring $A$. There exist non-zero central elements
$x_1,x_2,y_1,y_2$ of the ring $A$ such that $a_1 x_1=y_1$ and $a_2 x_2=y_2$. Then 
$$
Aa_1\cap Aa_2\supseteq Ax_1a_1\cap Ax_2a_2=a_1x_1A\cap a_2x_2A=y_1A\cap y_2A\neq 0.
$$
 
\textbf{c.} By the definition of the right singular ideal, all its elements are left zero-divisors. Conversely, let $a$ be a left or right zero-divisor of the ring $A$. Then $\text{r}(a)\neq 0$ by the first assertion of the lemma. In a right uniform ring, this means that $\text{r}(a)$ is an essential right ideal, i.e., $a\in B$. Now we use Remark 1.1.1.

\textbf{d.} Let $a,b\in A\setminus B$. There exist non-zero central elements $x,y\in A$ such that $b x=y$. Then $[a,b]x=[a,bx]=0$, i.e., $[a,b]$ is a left zero-divisor. Therefore, $[a,b]\in B$.

\textbf{e.} Let $r\in \ell.Ann_A(B)$. There exist non-zero central elements $x,y$ of the ring $A$ such that $rx=y$. It is clear that $x\not\in B$, whence we have $x$ is not a zero-divisor. Therefore, for every element $a\in A$, it follows from relations $0=[a,y]=[a,rx]=[a,r]x$ that $[a,r]=0$.~$\rhd$

\textbf{1.1.3. Closed right ideals.}\\ 
Let $A$ be a centrally essential ring and let $B$ be its right ideal.\\
\textbf{a.} If the right ideal $B$ is not essential (this is the case, if $B$ is a proper closed right ideal), then there exists a non-zero central element $y$ of the ring $A$ such that $B\cap yA=0$ and, consequently, $yB=By=0$. In particular, all elements of the right ideal $B$ are zero-divisors.

\textbf{b.} There exists a centrally essential finite-dimensional algebra over a field which has a closed right ideal which is not ideal.

$\lhd$ \textbf{a.} Since $B$ is not essential, $B\cap dA=0$ for some non-zero $d\in A$. Since $A$ is centrally essential, $dx=y$ for some non-zero central elements $x,y\in A$. Then $B\cap yA=0$.

\textbf{b.} See Example 3.6.8.~$\rhd$

\textbf{1.1.4. Central idempotents.}\\ 
If a not necessarily unital ring $A$ is centrally essential, then every idempotent $e\in A$ is contained in the center $Z(A)$.

$\lhd$ We can assume that $A$ is not commutative. Let $a\in A$. We have to prove that $ae=ea=eae$. First, we prove the relation $e(a-ae)=0$. We assume the contrary, $e(a-ae)\ne 0$. Since the ring $A$ is centrally essential, there exist $x,y\in Z(A)$ such that 
$$
xe(a-ae)=y=ey=ye\ne 0.
$$ 
Then
$$
0\ne y=ye=xe(a-ae)e=x(eae-eae)=0.
$$
This is a contradiction. Therefore, $e(a-ae)=0$. Similarly, we have $(a-ea)e=0$. Therefore, the idempotent $e$ is central.~$\rhd$

\textbf{1.1.5; \cite{LT22b}.} Let a ring $A$ be a not necessarily unital, centrally essential ring, $e=e^2\in A$, $a,x_1,\ldots, x_n,y_1,\ldots,y_n\in A$, and let
$$
\begin{cases}
x_1y_1 +\ldots + x_ny_n = e\\
x_1aey_1 +\ldots + x_naey_n = 0\\
\end{cases}.
$$
Then $ae = 0$.

$\lhd$ We assume that $ae\neq 0$. If the element $ae$ is central, then
$$
ae = ae^2 = ae(x_1y_1 +\ldots + x_ny_n) = x_1aey_1 +\ldots + x_naey_n = 0;
$$
this is a contradiction.

Now we assume that the element $ae$ is not central. Since the ring $A$ is centrally essential and $ae\ne 0$, there exist non-zero central elements $x,y\in A$ such that $xae = y$. We note that $y = ye$. Therefore,
$$
0\ne y = ye = y(x_1y_1 +\ldots + x_ny_n) = x(x_1yey_1 +\ldots + x_naey_n) = 0;
$$
this is a contradiction. Therefore, $ae = 0$.~$\rhd$

\textbf{1.1.6. Maximal right ideals.}\\ If $A$ is a centrally essential ring with $1\ne 0$ and $M$ is a maximal right ideal of the ring $A$, then either $M$ is an ideal or there exists a non-zero central element $x\in \left(\cap_{n\ge 1}M^n\right)$.

$\lhd$ We assume the contrary. Then there exist non-zero elements $m\in M$ and $a\in A$ such that $am\notin M$. Since $M$ is a maximal right ideal, there exist elements $b\in A$ and $m'\in M$ such that $1 = amb + m'$. Since the ring $A$ is centrally essential, there exist non-zero central elements $x,y\in A$ such that $ax = y$. Then
$$
x = (amb + m')x = (ax)mb + m'x = mby + m'x\in M
$$
and $(ax)mb\in M^2$ and $m'x\in M^2$. Therefore, $x = (ax)mb + m'x\in M^2$ and $(ax)mb, m'x\in M^3$. Then $x\in M^3$. By repeating a similar argument, we obtain that $0\ne x\in \left(\cap_{n\ge 1}M^n\right)$.~$\rhd$

\textbf{1.1.7. Proposition.}\\ Let $R$ be a ring and let $A$ be a subring in $R$ such that there exists a basis of the module $R_A$ contained in $Z(R)$. If the ring $A$ is centrally essential, then the ring $R$ is centrally essential, as well.

$\lhd$ Let $B$ be a basis of the module $R_A$ and $B\subseteq Z(R)$. Every element $r\in R\setminus \{0\}$ has the unique decomposition of the form 
$$
r=\sum_{i=1}^n b_is_i,\mbox{ where }b_1\ldots b_n\in B\mbox{ and }s_1\ldots s_n\in R\setminus\{0\}.\eqno(*)
$$
We define a function $k\colon R\rightarrow \mathbb{Z}$ by equating $k(r)$ to the number of coefficients $s_i$ in the above decomposition $(*)$ that are contained in $Z(A)$ for $r\neq 0$ and $k(0)=0$. It is clear that $r\in Z(A)$ if and only if $k(r)=0$. Now let $x\in R\setminus \{0\}$. In the set $xZ(A)\setminus \{0\}$, we take an element $r$ such that the integer $k(r)$ is minimal. We prove that $k(r)=0$. We assume the contrary. Then we can assume that $s_1\in Z(A)$ in $(*)$. Since the ring $A$ is centrally essential, there exist non-zero central elements $x,y\in Z(A)$ such that $xs_1=y$. Then $xr\in xZ(A)$, $xr\neq 0$ and $k(xr)<k(r)$; this contradicts to the choice of the element $r$. We obtain that $0\neq r\in xZ(A)\cap Z(R)\subseteq xZ(R)\cap Z(R)$.~$\rhd$

\textbf{1.1.8. Proposition.}\\
Let $F$ be a field and let $R$ be a centrally essential $F$-algebra. Then for any commutative $F$-algebra $A$, the algebra $A\otimes_F R$ is centrally essential.

$\lhd$ If $B$ is an $F$-basis of the commutative of the algebra $A$, then $\{b\otimes 1|b\in B\}$ is a basis of the free module $(A\otimes R)_R$ which satisfies conditions of Proposition 1.1.7.~$\rhd$

\textbf{1.1.9. Remark.}\\ 
If $A$ is a centrally essential ring with center $C=Z(A)$, then every its right ideal $B$ is an essential extension of the ideal $M=\oplus _{i\in I}c_iA$, $c_i\in C$ generated by central elements.

$\lhd$ Let $\mathcal M$ be a non-empty set of all ideals of the ring $A$ which are contained in $B$ and are direct sums of principal ideals generated by a central element. We define a partial order on $\mathcal M$ such that $M_1\le M_2$ $\Leftrightarrow$ $M_2 =M_1\oplus X$, $X\in \mathcal M$. By the Zorn lemma, $\mathcal M$ contains a maximal element $M$. We assume that $B_A$ is not an essential extension of $M_A$. Then there exists a non-zero element $b\in B$ such that $M\cap bA=0$. Since $A$ is a centrally essential ring, $bc=d$ for some non-zero central of elements $c,d\in A$. Then $M\cap dA=0$ and $M\oplus dA$ is an element of the set $\mathcal M$ which exceeds maximal element $M$. This is a contradiction.~$\rhd$

\textbf{1.1.10. Remark.}\\ 
It is directly verified that any filtered product of centrally essential rings is a centrally essential ring. In particular, ultra-degrees of a centrally essential ring are centrally essential rings.

\textbf{1.1.11. Open question.}\\ 
Is it true that any tensor product of centrally essential algebras is centrally essential?

\subsection{Semiprime and Nonsingular Rings}\label{subsection1.2}

We recall that a ring $A$ is said to be \textsf{semiprime} if $A$ does not have non-zero nilpotent ideals.
A ring $A$ with non-zero $1$ is said to be \textsf{right nonsingular}\label{nonsin} if the right annihilator $r_A(a)$ of any non-zero element $a\in A$ is not essential.

\textbf{1.2.1. Lemma.}\\
Let $A$ be a centrally essential ring with center $C=Z(A)$ and let $a$ be a non-zero element of the ring $A$. If $a^n=0$ ($n\in \mathbb{N}$), then there exists a non-zero central element $y$ of the ring $A$ such that $y\in (aC)\cap (Ca)$, $(AyA)^n=0$ and $(yC)^n=0$. Consequently, if at least one of the rings $A$ and $C$ is semiprime, then $A$ does not have non-zero nilpotent elements.

$\lhd$ Since $a\ne 0$ and the ring $A$ is centrally essential, $ax=xa=y$ for some non-zero central elements $x$ and $y$ of the ring $A$. Then
$$
(AyA)^n=y^nA^n=(ax)^nA^n=a^nx^nA^n=0.\quad \rhd
$$

\textbf{1.2.2. Theorem; see \cite[Theorem 1.3(a)]{MT19b}.}\\
Let $A$ be a centrally essential ring. If at least one of the rings $A$ and $Z(A)$ is semiprime, then $A$ is a commutative ring without non-zero nilpotent elements.

$\lhd$ By Lemma 1.2.1, the ring $A$ does not have non-zero nilpotent elements. We assume that the ring $A$ is not commutative. Then $ab-ba\neq 0$ for some $a,b\in A$. Let $C=Z(A)$ be the center of the ring $A$ and let $E = \{c\in C\,|\, ac \in C\}$. We have that $E$ is an ideal of the ring $C$. We take any element $d\in C$ with $dE = 0$. If $xd\neq 0$, then $xdz \in C\setminus \{0\}$ for some $z\in C$. Therefore, $dz \in E$, whence we have $d(dz) = 0$ and $(dz)^2 =0$. Therefore, $dz = 0$ and $xdz = 0$; this is a contradiction. Therefore, $xd = 0$, whence we have $d\in E$. Therefore, $d^2 = 0$ and $d =0$. Then we obtain that $\text{Ann}_C(E) = 0$. For any $i\in E$, we have $xi = ix\in C$, whence we have 
$$
[x, y]i= (xy - yx)i = x(yi) - y(xi) = xiy-xiy =0
$$
and $[x,y]E = 0$. However $c_1[x, y] = c_2$ for some $c_1,c_2\in C\setminus \{0\}$, whence we have $c_2E = 0$ and therefore, $\text{Ann}_C(E)\neq 0$; this is a contradiction. Therefore, the ring $A$ is commutative.~$\rhd$

\textbf{1.2.3. Remark.}\\
In connection to Theorem 1.2.2, we note that a ring $A$ with semiprime center $Z(A)$ is not necessarily commutative. The corresponding example is the ring $A$ of all $2\times 2$ matrices over $\mathbb{R}$; the center of the ring $A$ consists of scalar matrices.

\textbf{1.2.4. Corollary.}\\
If $A$ is a centrally essential, right nonsingular ring, then $A$ is commutative and does not have non-zero nilpotent elements.

$\lhd$ By Theorem 1.2.2, it is sufficient to prove that $A$ is a ring without non-zero nilpotent elements. We assume the contrary. There exists a non-zero element $a$ of the ring $A$ with $a^2=0$. Since $A$ is centrally essential, there exist non-zero central elements $x,y\in A$ with $ax=y$. It follows from of the Zorn lemma that there exists a right ideal $B$ of the ring $A$ such that $B\cap yA=0$ and right ideal $B\oplus yA$ is an essential. Since $yB=By\subseteq B\cap yA=0$ and $y^2=a^2x^2=0$, we have $y(B\oplus yA)=0$. Since right ideal $B\oplus yA$ is an essential, $y=0$; this is a contradiction.~$\rhd$

In connection to Theorem 1.2.2, we prove the following proposition.

\textbf{1.2.5. Proposition; \cite{LT22b}.}\\
If $A$ is a not necessarily unital, centrally essential ring and its the center is a semiprime ring, then the ring $A$ is commutative.

$\lhd$ We assume that the ring $A$ is not commutative, i.e. there exist elements $a,b\in A$ such that $ab-ba\neq 0$. Since the ring $A$ is centrally essential, there exist non-zero central elements $x$ and $y$ such that $(ab-ba)x=y$. We note that $ay\neq 0$; otherwise, 
$$
y^2 = (ab-ba)xy = ((ay)b-b(ay))x = 0;
$$
this is impossible, since $y\ne 0$. 

If $ay\notin Z(A)$, then there exist non-zero central elements $z,t\in A$ such that $ayz=t$. 
We consider the set $W = \{w\in Z(A)\,|\, aw\in Z(A)\}$. It is clear that $yz\in W$. Now we assume that $yW = 0$. Then $y(yz) = 0$, $(yz)^2 = 0$ and $yz = 0$; this is a contradiction. Therefore, $yw\neq 0$ for some $w\in W$. However, 
$$
yw = (ab - ba)yw = ((wa)b - b(wa))x = 0,
$$
and this is a contradiction, as well.

Therefore, we have $0\ne ay\in Z(A)$

We assume that $at\in Z(A)$. Then $ayb\neq 0$; otherwise, 
$$
y^2 =(ayb-bay)x=-bayx=aybx= 0.
$$
In addition, $(ab)y = (ba)y$. Therefore, $(ab - ba)y = 0$. However, $y^2 = (ab-ba)xy = 0$; this is a contradiction. Therefore, the ring $R$ is commutative.~$\rhd$

\textbf{1.2.6. Remark.} If $A$ is a ring and the factor ring $A/J(A)$ is centrally essential, then all maximal right ideals of the ring $A$ are ideals.

$\lhd$ Since $A/J(A)$ is a centrally essential semiprime ring, it follows from Theorem 1.2.2 that the ring $A/J(A)$ is commutative. In particular, all maximal right ideals of the ring $A/J(A)$ are ideals. Then all maximal right ideals of the ring $A$ are ideals.~$\rhd$

\subsection{Local and semiperfect Rings}\label{subsection1.3}

Let $A$ be a ring with Jacobson radical $J(A)$. The ring $A$ is said to be \textsf{local}\label{locring} if the factor ring $A/J(A)$ is a division ring.
The ring $A$ is said to be \textsf{semiperfect}\label{semiperfring} if the factor ring $A/J(A)$ is isomorphic to a finite direct product of matrix rings over division rings and every idempotent the factor of the ring $A/J(A)$ is the image of the idempotent $e\in A$ under the natural epimorphism $A\to A/J(A)$.

\textbf{1.3.1. Remark.}\\
It is clear that any finite direct product of local rings is a semiperfect ring. In addition, all idempotents of any centrally essential of the ring are central by 1.1.4. Therefore, centrally essential semiperfect rings coincide with finite direct products of centrally essential of local rings, and their the study is reduced to the study of centrally essential of local rings.

\textbf{1.3.2. Theorem.}\\
Let $A$ be a centrally essential local ring with Jacobson radical $J(A)$. Then the the factor ring $A/J(A)$ is a field (in particular, it is commutative) and $M\cap Z(A)\ne 0$ for every minimal right ideal $M$.

$\lhd$ Let $a,b\in A$ and $ab-ba\notin J(A)$. An element $ab-ba$ is invertible, since $A$ is local. Since $a\ne 0$ and $A$ is centrally essential, $ax=y$ for some non-zero $x,y\in Z(A)$. Then 
$$
x=x(ab-ba)(ab-ba)^{-1}=(yb-by)(ab-ba)^{-1}=0;
$$
this is a contradiction. Therefore, $ab-ba\in J(A)$ and the ring $A/J(A)$ is commutative.

Now we assume that $M\cap Z(A) = 0$ for some minimal right ideal $M$ of the ring $A$. Let $m$ be a non-zero element of $M$. By assumption, there exist non-zero central elements $x$ and $y$ of the ring $A$ such that $mx = y$. Since $x\notin J(A)$ (otherwise, $mx = 0$), the element $x$ is invertible in $A$ and $m = x^{-1}y\in Z(A)$; this is a contradiction.~$\rhd$

\textbf{1.3.3. Theorem.}\\ 
Let $A$ be a centrally essential semiperfect ring with center $C=Z(A)$.
Then $A/J(A)$ is a finite direct product of fields. In particular, the ring $A/J(A)$ is commutative.
In addition, $A$ is a finite direct product of centrally essential of local rings and $\text{Soc}(A_C)\subseteq C$.

$\lhd$ By the definition of a semiperfect ring, the ring $A/J(A)$ is the direct sum of simple Artinian rings and each of them is isomorphic to a matrix ring over a division ring. Let $\bar{e}_1,\ldots,\bar{e}_n$ be a complete system of indecomposable orthogonal idempotents of the ring $\bar A=A/J(A)$. Then there exists a complete system of indecomposable orthogonal idempotents $e_1,\ldots,e_n$ in $A$ such that $e_i+J(A)=\bar{e}_i$, $i=1,\ldots,n$. By 1.1.4, all idempotents $e_1,\ldots,e_n$ are central. Therefore, $A=\oplus_{i=1}^n A_ie_i$ be a decomposition of the ring $A$ into direct sum of local centrally essential rings. Consequently, all rings $A_i/J(A_i)$ are commutative by Theorem 1.3.2. It is directly verified that all rings $A_i=Ae_i$ are centrally essential; therefore, division ring $A_i/J(A_i)$ is commutative. Then the ring $R/J(R)=\oplus_{i=1}^n R_i/J(R_i)$ is commutative, as well.

It follows from the above that, without loss of generality, we can assume that the ring $A$ is local.
We note that $J(C)=C\cap J(A)$ and $C$ is a local ring.

Now let $s$ be a non-zero element of $\text{Soc}\,(R_C)$. There exist non-zero central elements $x,y$ such that $sx=y$. It is clear that $x\notin J(R)$, since $J(C)\text{Soc}\,A_C=0$. Consequently, $x$ is an invertible element and $s=x^{-1}y\in C$.~$\rhd$

\textbf{1.3.4. Remark.}\\ 
It follows from the above that if $A$ is a centrally essential semiperfect ring, then $\text{Soc}\,{}_AA=\text{Soc}\,A_A$.

\subsection{Perfect and Semi-Artinian Rings}\label{subsection1.4}

Let $A$ be a ring with Jacobson radical $J(A)$.

The ring $A$ is said to be \textsf{left perfect}\label{perf} if $A$ is semiperfect and the radical $J(R)$ is \textsf{left $T$-nilpotent},\label{tnilp} i.e. for any sequence $x_1,x_2,\ldots$ of elements in $J(A)$, there exists a subscript $n$ such that $x_1x_2\ldots x_n=0$. Right perfect rings are similarly defined.

The ring $A$ is said to be \textsf{semilocal}\label{semiloc} if the factor ring $A/J(A)$ is isomorphic to a finite direct product of matrix rings over division rings.

A module $M$ is said to be \textsf{semi-Artinian}\label{semiart} if either $M=0$ or every non-zero factor module of the module $M$ is an essential extension of a semisimple module.

\textbf{1.4.1. Theorem.}\\
Let $A$ be a right or left perfect ring with center $C=Z(A)$.

\textbf{a.} The ring $A$ is centrally essential if and only if $\text{Soc}\,A_C\subseteq C$ and all idempotents of the ring $A$ are central.

\textbf{b.} Assume that all idempotents of the ring $A$ are central, the factor ring $A/J(A)$ is commutative, $\text{Soc}\,A_C = \text{Soc}\,A_A$, and $M\cap C\neq 0$ for every minimal right ideal $M$. Then the ring $A$ is centrally essential.

$\lhd$ \textbf{a.} If $A$ is centrally essential, then $\text{Soc}\,A_C\subseteq C$ and all idempotents of the ring $A$ are central by 1.1.4 and Theorem 1.3.3.

Conversely, let $\text{Soc}\,A_C\subseteq C$ and let all idempotents of the ring $A$ be central. Since all idempotents are central, we can assume that $A$ is a local ring. Then $J(C)=C\cap J(A)$ and $C/J(C)$ is a field.

Let $x$ be a non-zero element of the ring $A$. If $J(C)x=0$, then $x\in\text{Soc}\,A_C$; therefore, $x\in C$. Otherwise, there exists an element $c_1\in J(C)$ such that $c_1x\neq 0$. If $J(C)c_1x=0$, then $c_1x\in\text{Soc}\,A_C$ and $c_1x\in C$; otherwise, we take an element $c_2\in J(C)$ such that $c_2c_1x\neq 0$, and so on. Since the radical $J(A)$ of the right perfect or left ring $A$ is a $T$-nilpotent right or left and elements $c_i$ are central, this process will stopped at some finite step.

\textbf{b.} By \textbf{a}, it is sufficient to prove the relation $\text{Soc}\,A_C\subseteq C$ which is equivalent to the property that $M\subseteq C$ for any minimal right ideal $M$. By assumption, $M\cap C\neq 0$ and the ring $A/J(A)$ is commutative by assumption; therefore, we have $ab - ba\in J(A)$, for all $a,b\in A$. For every $m\in M\cap C$ we have $m(ab - ba) = 0$. On the other hand, since $m\in C$, we have 
$$
(ma)b = mba = b(ma), \quad ma\in C.
$$
In addition, $ma\in M$. Consequently, $M\cap C$ is a non-zero right ideal of the ring $A$. Since $M$ is a minimal right ideal, $M\cap C = M$ and $M\subset C$. Therefore, $\text{Soc}\,A_C = \text{Soc}\,A_A\subseteq C$.~$\rhd$

\textbf{1.4.2. Remark.}\\
In Theorem 1.4.1, we cannot omit the condition that $R$ is right or left perfect, since every non-commutative local domain (for example, the formal power series ring in one variable over the Hamiltonian quaternion division ring) satisfies all remaining conditions of this theorem \label{thisorem} but this ring is not centrally essential.

\textbf{1.4.3. Lemma.}\\
Let $A$ be a semiprime ring and let $S=\text{Soc}\,A_A$ be the right socle. If $S$ is an essential right ideal of the ring $A$ and $st=ts$ for all $s,t\in S$, then the ring $A$ is commutative.

$\lhd$ We prove the following property $(*)$ of our ring $A$:

If $e=e^2\in A$ and $eA$ is a minimal right ideal, 
then the idempotent $e$ is central. Indeed, let $a\in A$. Since $S$ is an ideal, $ae\in S$. By assumption, 
$$
eae=e\cdot ae=ae\cdot e=ae.
$$
Similarly, $ea=eae=ae$ and the idempotent $e$ is central.

We prove that $ab-ba=0$ for any elements $a,b\in A$. We assume that $ab-ba\ne 0$. It is well known that every minimal right ideal of a semiprime ring is generated by an idempotent. Since $S$ is an essential right ideal and is generated, by the above, by idempotents, $S\cap (ab-ba)A\ne 0$ and $e(ab-ba)\ne 0$ for some idempotent $e\in A$. Then $eab\ne eba$ and $(ea)(eb)=(eb)(ea)$ by assumption. 

By properties $(*)$, the idempotent $e$ is central. Then 
$$
eab=eeab=eaeb=ebea=eba;
$$
this is a contradiction. Therefore, $A$ is commutative.~$\rhd$

\textbf{1.4.4. Lemma.}\\
Let $A$ be a centrally essential ring and let $P$ be a semiprime nil-ideal of the ring $A$ such that the right socle $S/R$ of the ring $A/P$ is an essential right ideal of the ring $A/P$. Then the ring $A/P$ is commutative.

$\lhd$ We use the following well known facts.

\textbf{a.} In any semiprime ring $R$, the set of all minimal right ideals coincides with the set of all minimal left ideals and this set coincides with the set of all right ideals $eR$ such that $e=e^2$ and $eRe$ is a division ring; in addition, $\text{\text{Soc}\,}R_R=\text{\text{Soc}\,}_RR$.

\textbf{b.} If $R$ is a ring and $P$ is a nil-ideal of the ring $R$, then every idempotent $\overline e$ of the ring $R/P$ is of the form $e+P$, where $e=e^2\in R$.

Let $h\colon A\to A/P$ be the natural epimorphism. For every subset $X$ in $A$, we write $\overline X$ instead of $h(X)$. By \textbf{a}, there exists an ideal $S$ of the ring $A$ such that $P\subset S$ and $\overline S=\text{\text{Soc}\,}_{\overline A}\overline A=\text{\text{Soc}\,}\overline A_{\overline A}$.

First, we show that the ideal $\overline S$ is commutative. By \textbf{a} and \textbf{b}, any minimal left ideal $V$ of the ring $\overline A$ is generated by some primitive idempotent $\overline e$ which is of the form $\overline e=e+P$ for some primitive idempotent $e$ of the ring $A$. By 1.1.4, the idempotent $e$ is central. Therefore, $V$ is an ideal of $\overline A$, $eA$ and $(1-e)A$ are ideals in $A$, and $A=eA\oplus (1-e)A$. Therefore, the ring $eA$ is centrally essential. In addition, $\overline e\overline A\overline e=\overline e\overline A=V$ and $V=(eA+P)/P\cong eA/(P\cap eA)$. Therefore, $J(eA)\subseteq P\cap eA$. But $P$ is a nil-ideal, whence we have $P\cap eA\subseteq J(eA)$ and $P\cap eA=J(eA)$. By Theorem 1.3.2, the ring $V$ is commutative. Therefore, $\text{\text{Soc}\,}(\overline A)$ is a commutative ring, as a direct sum of commutative rings. In addition, $\text{\text{Soc}\,}(\overline A)$ is an essential right ideal of the semiprime ring $\overline A$. Then $\overline A$ is commutative by Lemma 1.4.3.~$\rhd$
 
\textbf{1.4.5. Theorem.}\\ 
If $A$ is a centrally essential, left or right semi-Artinian ring, then $A/J(A)$ is a commutative (von Neumann)\label{vnrring} regular ring.

$\lhd$ Let $A$ be a centrally essential semi-Artinian right or left ring and $\overline A=A/J(A)$. Since $A$ is a right or left semi-Artinian ring, $J(A)$ is a nil-ideal by \cite[Proposition 3.2]{NasP68}. By Lemma 1.4.4, the ring $A/J(A)$ is commutative. Every commutative semi-Artinian semiprimitive ring is von Neumann regular by \cite[Theorem 3.1]{NasP68}.~$\rhd$

\section{Graded Rings and Grassmann Algebras}\label{section2}

Subsections $2.1$ and $2.2$ are based on \cite{MT19a}.

\subsection{Graded Rings}\label{subsection2.1}

\textbf{2.1.1. Graded rings and homogeneous elements.}\\
Let $(S,+)$ be a semigroup. A ring $A$ is said to be \textsf{$S$-graded}\label{sgrad} if $A$ is a direct sum of additive subgroups $A_s$, $s\in S$, and $A_sA_t\subseteq A_{s+t}$ for any elements $s,t\in S$. 

For any $s\in S$, elements of the subgroup $A_s$ are called \textsf{homogeneous}\label{homogen} elements of degree $s$.

If $S=\mathbb{N}\cup\{0\}$, then $S$-graded rings are called \textsf{graded}\label{grad} rings. It is directly verified that the identity element of a graded ring is contained in the subgroup $A_0$. On an arbitrary graded ring $A=\oplus_{n\in\mathbb{N}\cup \{0\}}A_n$, we can define $\mathbb{Z}_2$-graduation:
$$
A=A_{(0)}\oplus A_{(1)}, \text{ where } A_{(i)}=\bigoplus_{k\in
\mathbb{N}\cup \{0\}}A_{2k+i},\; i\in \{0,1\}.
$$

\textbf{2.1.2. Generalized anti-commutative and\\ homogeneously faithful rings.}\\
One says that a graded ring $A=\oplus_{n\in\mathbb{N}\cup \{0\}}A_n$ is \textsf{generalized anti-commutative}\label{genanti} if for any integers $m,n\in\mathbb{N}\cup \{0\}$ and arbitrary elements $x\in A_m$ and $y\in A_n$, the relation $yx=(-1)^{mn}xy$ holds.

If the graded ring $A=\oplus_{n\in\mathbb{N}\cup \{0\}}A_n$ satisfies the condition
$$
\forall m,n\in \mathbb{N}\cup \{0\},\, A_{m+n}\neq 0 \Rightarrow A_m\neq 0 \&\ \forall x\in{A_m}\setminus\{0\},\, xA_n\neq 0,\eqno{(*)}
$$ 
then one says that $R$ is a \textsf{homogeneously faithful}\label{homfai} ring.

\textbf{2.1.3. The center of a graded ring.}\\ 
In any graded ring $A=\oplus_{n\in\mathbb{N}\cup \{0\}}A_n$, the relation $Z(A)=\oplus_{n\in\mathbb{N}\cup\{0\}}(A_n\cap Z(A))$ holds.

(\textbf{Remark.} If $S$ is a commutative cancellative semigroup, then the proof below remains true for every $S$-graded ring.)

$\lhd$ The inclusion $\oplus_{n\in\mathbb{N}\cup \{0\}}(A_n\cap Z(A))\subseteq Z(A)$ is obvious. 

Let $x=x_0+x_1+\ldots x_n\in Z(A)$, where $x_i\in A_i$, $i=0,1,\ldots,n$. If $y\in A_m$ for some $m\in\mathbb{N}\cup \{0\}$, then $0=[x,y]=[x_0,y]+\ldots+[x_n,y]$ and summands of the last sum are contained in distinct direct summands $A_{m},A_{m+1},\ldots,A_{m+n}$. Therefore, $[x_i,y]=0$ for any homogeneous element $y$ and of all $i=0,1,\ldots,n$. Then $x_i\in Z(A)$, since any element of the ring is a sum of homogeneous elements.~$\rhd$

\textbf{2.1.4.}
Let $A=\oplus_{n\in\mathbb{N}\cup \{0\}}A_n$ be a graded generalized anti-commutative homogeneously faithful ring which does not have additively 2-torsion elements. If there exists an odd positive integer $n$ such that $A_n\neq 0$ and $A_{n+1}=0$, then $Z(A)=A_{(0)}+A_n$. Otherwise, $Z(A)=A_{(0)}$.

$\lhd$ It follows from the generalized anti-commutativity relation that $A_{(0)}\subseteq Z(A)$. The following property follows from $(*)$:\\
if such an integer $n$ exists, then $A_m=0$ for $m>n$ and $A_m\neq 0$ for $0\leq m \leq n$; in addition, if $x\in A_n$ and $y=y_0+z\in A$, where $y_0\in A_0$ and $z\in\oplus_{m>0}A_m$, then $[x,y]=[x,y_0]=0$, i.e., $A_n\subseteq Z(A)$. Conversely, let $x\in Z(A)$. By 2.1.3, we can assume that $x$ is a homogeneous element of odd of degree $i$. Let $x \neq 0$ and $A_{i+1}\neq 0$. Then it follows from $(*)$ that there exists an element $y\in A_1$ such that $xy\neq 0$. We obtain that $0=[x,y]=2xy$; this is a contradiction. Therefore, either $x=0$ or $x\neq 0$ but $A_{i+1}=0$, i.e., $i=n$.~$\rhd$

\textbf{2.1.5. Theorem.}\\ 
Let $A=\oplus_{n\in\mathbb{N}\cup \{0\}}A_n$ be a graded generalized anti-commutative homogeneously faithful ring without additively $2$-torsion elements. The ring $A$ is centrally essential if and only if either $A=A_0$ or there exists an odd positive integer $n$ such that $A_n\neq 0$ and $A_{n+1}=0$.

$\lhd$ Let $A$ be a centrally essential ring, $C=Z(A)$, and let $A\neq A_0$. By $(*)$, we have $A_1\neq 0$. We take an element $x\in A_1\setminus\{0\}$ and assume that such an integer $n$ does not exist. By 2.1.3, we have $C=A_{(0)}$ and $xC\subseteq A_{(1)}$, whence we have $xC\cap C\subseteq A_{(1)}\cap A_{(0)}=0$; this is a contradiction.

Conversely, if $A=A_0$, then $C=A$, since the ring $A_0$ is commutative. We assume that there exists an odd positive integer $n$ such that $A_n\neq 0$ and $A_{n+1}=0$. Let $0\neq x\in A\setminus C$. We have $x=x_0+\ldots+x_n$, where $x_i\in A_i$, and we take the least odd positive integer $m$ such that $x_m\neq 0$. It is clear that $1\leq m\leq n$. We set $k=n-m$ and take an element $y\in A_k$ such that $x_my\neq 0$. It is clear that $y\in C$. In addition, $xy$ is a sum of homogeneous elements of even degree and the element $x_my$ of odd degree $n$. Therefore, $xy\in C$ by 2.1.3 and $xy\neq 0$.~$\rhd$

\subsection{Grassmann Algebras over Fields}\label{subsection2.2}

\textbf{2.2.1.}\label{Grassmann} Let $F$ be a field of characteristic $0$ or $p>2$, $V=F^n$ be a vector space over $F$ of dimension $n>0$, and let $\Lambda(V)$ be the Grassmann algebra of the space $V$ \cite[\S III.5]{Bou89} which is defined as a unital $F$-algebra with respect to multiplication operation $\wedge$ with generators $e_1,\ldots,e_n$ and defining relations $e_i\wedge e_j+e_j\wedge e_i=0$ for all $i,j\in\{1,\ldots,n\}$.

The algebra $\Lambda(V)$ has a natural graduation:
$$
\Lambda(V)=\bigoplus_{p\in\mathbb{N}\cup \{0\}}\Lambda^p(V),
$$
where $\Lambda^p(V)$, $1\leq p\leq n$, is a vector space with basis
$$
\{e_{i_1}\wedge \ldots \wedge e_{i_p}:1\leq i_1<\ldots<i_p\leq
n\},
$$
$\Lambda^0(V)=F$ and $\Lambda^p(V)=0$ for $p>n$.

It is well known that Grassmann algebras are generalized anti-commutative.

\textbf{2.2.2.} The graded algebra $R=\Lambda(V)$ is a homogeneously faithful ring.

$\lhd$ Let $p,q\in\{0,\ldots,n\}$ and $p+q\leq n$. If $pq=0$, then the condition $(*)$ holds. Now let $0<p<n$ and $0\neq x\in R_p$. We take a basis element $e_{i_1}\wedge \ldots \wedge e_{i_p}$ which has a non-zero coefficient in the representation of $x$. Since $p+q\leq n$, there exist subscripts $j_1,\ldots,j_q\in\{1,\ldots,n\}$ such that $1\leq j_1< \ldots <j_q\leq n$ and $\{i_1,\ldots,i_p\}\cap\{j_1,\ldots,j_q\}=\varnothing$.
We set $y=e_{j_1}\wedge\ldots\wedge e_{j_q}$ and note that the basis element $\pm e_{i_1}\wedge \ldots \wedge e_{i_p}\wedge
e_{j_1}\wedge\ldots\wedge e_{j_q}$ of the space
$\Lambda^{p+q}(V)$ has the non-zero coefficient in the representation of the element $xy$, since products of remaining basis elements of the space $\Lambda^{p}(V)$ by the element $y$ are equal to either $0$ or $\pm$other basis elements of the space $\Lambda^{p+q}(V)$.~$\rhd$

\textbf{2.2.3. Theorem.}\\
Let $F$ be a field of characteristic $0$ or $p>2$ and let $V$ be a finite-dimensional vector $F$-space. The Grassmann algebra $\Lambda(V)$ of the space $V$ is a centrally essential ring if and only if $V$ is of odd dimension. 

Theorem 2.2.3 follows from 2.2.2 and 2.1.4.

\textbf{2.2.4.} 
If $F$ is a finite field of odd characteristic and $\dim V$ is an odd positive integer exceeding $1$, then $\Lambda(V)$ is a centrally essential non-commutative finite ring. Thus, if $F$ is the finite field of order $3$ and $\Lambda(V)$ is an $8$-dimensional $F$-algebra with basis
$$
\{1,e_1,e_2,e_3,e_1\wedge e_2,e_1\wedge e_3,e_2\wedge e_3,e_1\wedge e_2\wedge e_3\},
$$
then $\Lambda(V)$ is a centrally essential non-commutative finite ring of order $3^8$.

\textbf{2.2.5. Example.} 
If $R$ is a centrally essential ring and $B$ is a proper ideal of the ring $R$ generated by some infinite set of central idempotents and the factor ring $R/B$ does not have non-trivial idempotents, then the ring $R/B$ is not necessarily centrally essential.

$\lhd$ Let $F$ be a field of order three, $A=\Lambda(F^3)$ be the Grassmann algebra of the three-dimensional vector $F$-space $F^3$, and let $S=\Lambda(F^2)$ be the Grassmann algebra of two-dimensional vector $F$-of the space $F^2$ considered as a subalgebra of the algebra $A$. We consider the direct product $P=A^{\mathbb N}=\{(a_1,a_2,\ldots)\,|\,a_i\in A\}$ of countable set of copies of the ring
$A$ and its subring $R$ consisting of all eventually constant sequences $(a_1,a_2,\ldots)\in P$ which stabilize at finite step on elements of $S$ depending on the sequence.

Let $e_i$ be a central idempotent which has the identity element of the field $F$ on $i$-th position and zeros on remaining positions. We denote by $B$ the ideal of the ring $R$ generated by all
idempotents $\{e_i\}$. It follows from Theorem 2.2.3 that $R$ is a centrally essential ring and the factor ring $R/B$ is isomorphic to the ring $S$ which is not centrally essential and does not have non-trivial idempotents.~$\rhd$

\subsection{Grassmann Algebras over Rings}\label{subsection2.3}

This subsection is based on \cite{MT19b}.

\textbf{2.3.1.} 
Let $A$ be a not necessarily commutative ring with center $C=Z(A)$ and let $A^n$ be a finitely generated free module of rank $n$. We define the algebra $\Lambda(A^n)$ of the module $A^n$. Namely,~ $\Lambda(A^n)=A\otimes_C\Lambda(C^n)$, where $\Lambda(C^n)$ is the Grassmann algebra of the free module $C^n$ over the commutative ring $C$; see \cite[\S III.5]{Bou89}.

Let $\{e_1,\ldots,e_n\}$ be a basis of the module $C^n$. For all $x\in \Lambda(C^n)$, we identify $1\otimes x$ with $x$ and obtain that the set 
$$
B_n=\{e_{i_1}\wedge\ldots\wedge e_{i_s}|0\leq s\leq n,\;\; 1\leq i_1<\ldots<i_s\leq n\}
$$
is a basis of the $A$-module $\Lambda(A^n)$ (we assume that the product is equal to $1$ for $s=0$). It is clear that the ring $R=\Lambda(A^n)$ has a natural graduation $R=\oplus_{s\geq 0} R_s$, where $R_0=A$, $R_s=\oplus_{1\leq i_1<\ldots<i_s\leq n}Ae_{i_1}\wedge\ldots\wedge e_{i_s}$ for
$1\leq s\leq n$, and $R_s=0$ for $s>n$.

\textbf{2.3.2. Theorem.} 
For positive integer $n$ and a ring $A$ with center $C=Z(A)$, the ring $\Lambda(A^n)$ is centrally essential if and only if $A$ is centrally essential and at least one of the following conditions holds.
\begin{itemize}
\item[\textbf{a)}]
$n$ is an odd integer.
\item[\textbf{b)}]
The ideal $\text{Ann}_A(2)$ is an essential submodule of the module $A_C$.
\end{itemize}

$\lhd$ We set $R=\Lambda(A^n)$. Let $R$ be centrally essential.

Let $a\in A\setminus \{0\}$ and $a'=ae_1\wedge\ldots\wedge e_n$. Then $0\neq a'\in R$. Therefore, there exists an
element $c\in Z(R)$ such that $0\neq ca'\in Z(R)$. We have $c=c_0+c'$, where $c_0\in A$, and $c'\in\oplus_{s>0}R_s$.
It is directly verified that $c_0\in Z(A)=Z(R)\cap R_0$. In addition, it is clear that $ca'=c_0a'=c_0ae_1\wedge
\ldots \wedge e_n$, whence we have $c_0a\neq 0$. For every $b\in A$, we have
$$
0=[b,c_0ae_1\wedge\ldots\wedge e_n]=[b,c_0a]e_1\wedge\ldots
\wedge e_n,
$$
whence we have $c_0a\in Z(A)$, i.e., $A$ is a centrally essential ring.

We assume that the ideal $\text{Ann}_A(2)$ is not an essential submodule of the module $A_C$ and $n$ is an even integer.

We take an element $a\in A$ such that $a\neq 0$ and $Ca\cap\text{Ann}_A(2)=0$. We consider the element $x=ae_2\wedge\ldots\wedge e_n$.

Let $c\in Z(R)$ and $0\neq cx\in Z(R)$. We have $c=c_0+c_1e_1+c'$, where $c'$ is a linear combination of elements of the basis $B_n$ which are equal to 1 and $e_1$. It is clear that $c_0,c_1\in C$ and $cx=c_0ae_2\wedge\ldots\wedge
e_n+c_1ae_1\wedge\ldots\wedge e_n$, where the both summands are contained in the center of the ring $R$.

We prove that $c_0a=0$. Indeed, 
$$
0=[e_1, c_0ae_2\wedge\ldots\wedge e_n]=c_0ae_1\wedge\ldots\wedge e_n-c_0ae_2\wedge\ldots\wedge e_n\wedge e_1=
$$
$$
=c_0a(1-(-1)^{n-1})e_1\wedge\ldots\wedge e_n=2c_0ae_1\wedge\ldots\wedge e_n,
$$
whence we have $c_0a\in \text{Ann}_A(2)\cap Ca=0$ by the choice of $a$. Then $c_1a\neq 0$ and $c_1e_1\in Z(R)$. However $c_1\in C$, $0=[c_1e_1,e_2]=2c_1e_1\wedge e_2$. Then $c_1a\in Ca\cap\text{Ann}_A(2)=0$. This is a contradiction.

Now we assume that $A$ is centrally essential and at least one of the above conditions \textbf{a)} or \textbf{b)} holds.

Let \textbf{a)} hold. We set $N=\text{Ann}_C(2)=C\cap\text{Ann}_R(2)$. We note that $N$ is an essential submodule in $A_C$. We consider an arbitrary non-zero element $x\in R$. We have
$$
x=\sum_{s=0}^n\sum_{1\leq i_1<\ldots<i_s\leq n}a_{i_1,\ldots,i_s}e_{i_1}\wedge\ldots\wedge e_{i_s},
$$
where coefficients $a_{i_1,\ldots,i_s}$ are contained in $A$. We can multiply $x$ by elements $C\subseteq Z(R)$ and obtain a situation, where all coefficients in representation of $x$ are contained in $N$. Indeed, if some coefficient $a_{i_1,\ldots,i_s}$ is not contained in $N$, then there exists an element $c\in C$ such that $0\neq ca_{i_1,\ldots,i_s}\in N$, i.e. under multiplication by $c$, the number of coefficients, contained in $N$, decreases. It remains to note that
$x\in Z(R)$ if all coefficients of their representation of $x$ are contained in $N$. Indeed, $[x,a]=0$ for any $a\in A$, since $N\subseteq Z(A)$ and
$$
[x,e_i]=\sum_{s=0}^n\sum_{i_1<\ldots<i_s}a_{i_1,\ldots,i_s}[e_{i_1}\wedge\ldots\wedge e_{i_s},e_i].
$$
We note that if the number $s$ is even or $i\in\{i_1,\ldots,i_s\}$, then 
$[e_{i_1}\wedge\ldots\wedge e_{i_s},e_i]=0$. Otherwise,
$$
[e_{i_1}\wedge\ldots\wedge e_{i_s},e_i]=\alpha e_{i_1}\wedge\ldots\wedge e_{i_s}\wedge
e_i,
$$
where $\alpha\in\{0,2\}$, i.e. we have $[a_{i_1,\ldots,i_s}e_{i_1}\wedge\ldots\wedge e_{i_s},e_i]=0$.
Since elements of the ring $A$ and $e_1,\ldots,e_n$ generate the ring $R$, we have $x\in Z(R)$, which is required.

Now we assume that the condition \textbf{b)} holds. We consider an arbitrary non-zero element $x\in R$. By repeating argument from the previous case, we can use multiplication by elements of $C$ to obtain such a situation that all coefficients $x$ with respect to the basis $B_n$ are contained in $C$.

We take the least odd $k$ such that the element $e_{i_1}\wedge\ldots\wedge e_{i_k}$ of the basis $B_n$ is contained in the representation of $x$ with non-zero coefficient $a$ (if this is impossible, then $x\in Z(R)$). Let 
$$
m=n-k,\; \{j_1,\ldots,j_q\}=\{1,\ldots,n\}\setminus \{i_1,\ldots,i_k\}.
$$ 
It is clear that integer $m$ is even, whence we have $c=e_{j_1}\wedge\ldots\wedge e_{j_m}\in Z(R)$. Then it is directly verified that $cx=\pm ae_1\wedge\ldots\wedge e_n+x'$, where $x'$ is a linear combination of elements of the basis $B_n$ with even degree $s$ and coefficients in $C$. Therefore, we repeat the argument from the previous case and obtain that $x'\in Z(R)$. Finally, it is directly verified that $ae_1\wedge\ldots\wedge e_n\in Z(R)$. The assertion is proved.~$\rhd$

\textbf{2.3.3. Lemma.}
If $A$ is a ring of finite characteristic $s$ and $C=Z(A)$, then the following conditions are equivalent.

\textbf{a)} The ideal $\text{Ann}_A(2)$ is an essential submodule of the module $A_C$.

\textbf{b)} $s=2^m$ for some $m\in\mathbb{N}$.

$\lhd$ a)\,$\Rightarrow$\,b). We assume the contrary. Then there exists an odd prime integer $p$ dividing $s$. The non-zero ideal $\text{Ann}_A(p)$ of the ring $A$ has the zero intersection with the ideal $\text{Ann}_A(2)$. Therefore, the ideal $\text{Ann}_A(2)$ is not an essential submodule of the module $A_C$. This is a contradiction.

b)\,$\Rightarrow$\,a). Since $s=2^m$, 
we have that for every $a\in A\setminus\{0\}$, the relation $\text{ord}\,{a}=2^k$ holds for some $k\in \mathbb{N}$. Then $0\neq 2^{k-1}a\in Ca\cap\text{Ann}_A(2)$. Therefore, the ideal $\text{Ann}_A(2)$ is an essential submodule of the module $A_C$.~$\rhd$

If $A$ is a ring of finite of characteristic or $A$ does not have zero-divisors, then the formulation of Theorem 2.3.2 can be simplified; see Theorem 2.3.4.

\textbf{2.3.4. Theorem.} 
Let $A$ be a ring with center $C=Z(A)$ and let $n$ be a positive integer.

The ring $\Lambda(A^n)$ is centrally essential if and only if $A$ is centrally essential and at least one of the following conditions holds.

\textbf{a.} If $A$ is a ring of finite characteristic $s$ (this is the case if the ring $A$ is finite), then the ring $\Lambda(A^n)$ is centrally essential if and only if the ring $A$ is centrally essential and at least one of the following conditions holds.
\begin{itemize}
\item
$n$ is an odd integer.
\item
$s=2^m$ for some $m\in\mathbb{N}$.
\end{itemize}

\textbf{b.} If $A$ is a domain, then the ring $\Lambda(A^n)$ is centrally essential if and only if the ring $A$ is centrally essential and at least one of the following conditions holds.
\begin{itemize}
\item
$n$ is an odd integer.
\item
$A$ is a ring of characteristic $2$
\end{itemize}

$\lhd$ We set $R=\Lambda(A^n)$.

\textbf{1.} The assertion follows from Theorem 2.3.2 and Lemma 2.3.3.

\textbf{2.} If $A$ is a ring of characteristic $2$ or $n$ is an odd integer, then $R$ is a centrally essential ring by \textbf{a}.

Now we assume that $A$ is a domain and the ring $R$ is centrally essential. By Theorem 2.3.2, the ring $A$ is centrally essential and it is sufficient to consider the case, where $n$ is an even integer and ideal $\text{Ann}_A(2)$ is an essential submodule of the module $A_C$. Since $A$ is a domain, $\text{Ann}_A(2)=A$. Therefore, $A$ is a ring of characteristic $2$.~$\rhd$

\section{Constructions of Rings}\label{section3}

\subsection{Polynomials, Series and Fractions}\label{subsection3.1}

Subsection 3.1 is based on \cite{MT20a}.

For arbitrary finite subset $S$ of the monoid $G$ and any ring $A$, we denote by $\Sigma_S$\label{sigmaS} the element $\sum_{x\in S}x$ of the monoid ring $AG$. For any element $r=\sum_{g\in G}a_g\cdot g\in AG$, we say that the set $\{g\in G|a_g\neq 0\}$ is the \textsf{support}\label{supp} of the element $r$;
we denote this set by $\text{supp}\,(r)$.\label{supp}

\textbf{3.1.1. Monoid rings.}\\
If $A$ is a centrally essential ring and $G$ is a commutative monoid, then the monoid ring $R=AG$ is centrally essential.

$\lhd$ For any non-zero element $r=\sum_{g\in G}r_g\cdot g\in R$, let 
$$
k(r)=|\{g\in G|r_g\in Z(A)\}|.
$$ 
It is clear that $k(r)\leq|\text{supp}\,(r)|<\infty$. With the use of the induction on $k$, we prove that for $k(r)=k$, there exist non-zero central elements $x$ and $y$ such that $rx=y$. 

If $k=0$ then $r_g\in Z(A)$ for all $g\in G$ and, therefore, $r\in Z(R)$. 

Otherwise, if $k>0$ and $k(r)=k$, then we can take an element $h\in G$ with $r_h\in Z(A)$. Since the ring $A$ is centrally essential, there exist non-zero central elements $x$ and $y$ with $xr_h=y$. It is clear that $0\neq xr=\sum_{g\in G}xr_g\cdot g$ and $k(xr)<k(r)$. By the induction hypothesis, there exist non-zero central elements $u$ and $v$ of the ring $R$ such that $uxr=v$. Since $ux\in Z(R)$, the proof is completed.~$\rhd$

The following assertion is a corollary of 3.1.1.

\textbf{3.1.2. Polynomial rings.}\\
For any centrally essential ring $A$, the polynomial ring $A[x]$ and the polynomial Laurent ring $A[x,x^{-1}]$ are centrally essential.

\textbf{3.1.3. Central rings of fractions.}\\
Let $A$ be a centrally essential ring, $Q$ be the ring of fractions of the ring $A$ with respect to some central multiplicative system $S$ consisting of non-zero-divisors, and let $0\ne s^{-1}a=as^{-1}\in Q$, $s\in S$. By assumption, there exist non-zero central elements $x,y\in A$ with $ax=y$. Then $0\ne as^{-1}x=s^{-1}y$ is a central element of the ring $Q$ and the ring $Q$ is centrally essential. If $Q$ is a centrally essential ring, then it is similarly proved that $A$ is a centrally essential ring. 

\textbf{3.1.4. Remark.}\\
Since the formal Laurent series ring $A((x))$ is the ring of fractions of the ring formal power series $A[[x]]$ with respect to the central multiplicative system $\{x_k\}_{k=0}^{\infty}$, it follows from 3.1.3 that the ring $A((x))$ is centrally essential if and only if the ring $A[[x]]$ is centrally essential. 

\textbf{3.1.5. Proposition.}\label{series}
If $R$ is a finite-dimensional centrally essential algebra, then the formal power series ring $R[[x]]$ is centrally essential.

$\lhd$
It is sufficient to prove that the ring $R[[x]]$ is isomorphic to $F[[x]]\otimes R$. First, we prove injectivity of a natural homomorphism $\varphi\colon F[[x]]\otimes R\rightarrow R[[x]]$ defined by the relation $\varphi(f(x)\otimes r)=f(x)r$ for every $f(x)\in F[[x]]$ and $r\in R$. Indeed, every element of the algebra $F[[x]]\otimes R$ can be represented in the form
$r=\sum_{i=1}^n f_i(x)\otimes r_i$, where $f_1(x),\ldots,f_n(x)\in F[[x]]$, $r_1,\ldots,r_n$ are linearly independent elements of the algebra $R$ (for example, $\{f_1(x),\ldots,f_n(x)\}$ can be subset of some fixed finite or infinite of the basis for $R$).
By assuming $f_i(x)=\sum_{j=0}^\infty x^j\alpha_{ij}$ for some 
$\alpha_{ij}\in F$, we have that $\varphi(r)=\sum_{i=1}^n(\sum_{j=0}^\infty x^j\alpha_{ij})r_i=
\sum_{j=0}^\infty x^j(\sum_{i=1}^n\alpha_{ij}r_i)$. Therefore, if $\varphi(r)=0$, then for all $j\geq 0$, we have
$\sum_{i=1}^n\alpha_{ij}r_i=0$, whence we have $\alpha_{ij}=0$ and $f_i(x)=0$ for all $i=1,\ldots,n$; consequently $r=0$.

If $R$ is a finite-dimensional algebra with basis $r_1,\ldots,r_n$, then for any series $f(x)=\sum_{j=0}^\infty x^jt_j$ with coefficients $t_j\in R$, we have $t_j=\sum_{i=1}^n\alpha_{ij}r_i$, whence we have
$$
f(x)= \sum_{j=0}^\infty x^j(\sum_{i=1}^n\alpha_{ij}r_i)=\sum_{i=1}^n(\sum_{j=0}^\infty x^j\alpha_{ij})r_i\in \varphi(F[[x]]\otimes R).\;\rhd
$$

\textbf{3.1.6. Theorem.}\\ If $R$ is a finite-dimensional centrally essential algebra, then the following conditions are equivalent.

\textbf{1)} The ring $R$ is centrally essential.

\textbf{2)} The power series ring $R[[x]]$ is centrally essential.

\textbf{3)} The ring Laurent series $R((x))$ is centrally essential.

$\lhd$ The implication 1)\,$\Rightarrow$\,2) follows from Proposition 3.1.5. 

The implication 2)\,$\Rightarrow$\,1) is directly verified. 

The equivalence 2)\,$\Leftrightarrow$\,3) follows from 3.1.4.~$\rhd$

\subsection{Group Rings}\label{subsection3.2}

Subsection 3.2 is based on \cite{MT20a}.

Let $R$ be a ring and let $G$ be a group. We set $(x,y)=x^{-1}y^{-1}xy$ for any elements $x,y$ of the group $G$; additive commutators and multiplicative commutators are denoted differently, since elements of the group are considered as elements of the group ring, as well. For any element $g$ of the group $G$, we denote by $g^G$\label{gG} class of conjugated elements which contains $g$. For the group $G$, the \textsf{upper central series}\label{upcen} of the group $G$ is a chain of subgroups $\{1\}=Z_0(G)\subseteq Z_1(G)\subseteq\ldots$, where $Z_i(G)/Z_{i-1}(G)$ is the center of the group $G/Z_{i-1}(G)$, $i\geq 1$. We denote by $\text{NC}(G)$ the \textsf{nilpotence class}\label{nilpCla} of the group $G$, i.e. the least positive integer $n$ with $Z_n(G)=G$ (if it there exists a).

A group $G$ is called an \textsf{$FC$-group}\label{FCgr} if all classes of conjugated elements in $G$ are finite.

\textbf{3.2.1. Proposition.}\label{gab1}\\
Let $A$ be a ring and $G$ be a group. If the group ring $R=AG$ is centrally essential, then $A$ also is a centrally essential ring and the group $G$ is an $FC$-group.

$\lhd$ Let $0\neq a\in A$. Since $A\subseteq R$ and $R$ is centrally essential, there exists an element $c\in Z(R)$ such that $0\neq ca\in Z(R)$. We have $c=\sum_{g\in G}c_g\cdot g$ and $ca=\sum_{g\in G}c_ga\cdot g$. It follows from relations $0=[c,b]=\sum_{g\in G}[c_g,b]\cdot g$ for any $b\in A$ that $c_g\in Z(A)$ for all $g\in G$. Similarly, we have $c_ga\in Z(A)$ for any $g\in G$. Since there exists at least one element $g\in G$ with $c_ga\neq 0$, we obtain our assertion on the ring $A$.

Now let $g$ be an arbitrary element of the group $G$. It is well known (e.g., see \cite[Lemma 4.1.1]{Pas77}) that $Z(AG)$ is a free $Z(A)$-module with basis
$$
\{\Sigma_K\,|\,K\mbox{ is a finite class of conjugated elements in }G\}.\eqno{(3.2.1)}
$$
In particular,
$$
r\in Z(AG)\Rightarrow |g^G|<\infty \text{ for any } g\in \text{supp}\,(r).\eqno{3.2.1(2)}
$$
Since $AG$ is centrally essential, $0\neq cg=d$ for some $c,d\in Z(AG)$. By comparing coefficients in the left part and the right part of the relation $cg=d$, we obtain that for any $y\in\text{supp}\,(d)$, there exists an element $x\in\text{supp}\,(c)$ such that $xg=y$. For any $h\in G$, we have $hgh^{-1}=(hxh^{-1})^{-1}hyh^{-1}$, whence we have $g^G\subseteq (x^{-1})^G\cdot y^G$. Since $|(x^{-1})^G|=|x^G|$, we have that
$|g^G|\leq |x^G|\cdot |y^G|<\infty$, by 3.2.1(2).~$\rhd$

\textbf{3.2.2. Lemma.}\label{subgroups}\\
Let $G$ be a group, $F$ be a field of characteristic $p>0$, and
let $q$ be a prime integer which is not equal to $p$. If the ring
$FG$ is centrally essential, then every $q$-subgroup in $G$ is a
normal commutative subgroup.

$\lhd$ First, let $H$ be a finite $q$-subgroup 
of the group $G$. Then $|H|=n=q^k$ is a non-zero element of the field $F$ and the element $e_H=\frac 1 n \Sigma_H$ is an idempotent of the ring $FG$. By 1.1.4, $e_H$ is a central idempotent. Consequently, $ge_Hg^{-1}=\frac 1 n \sum_{h\in H}ghg^{-1}=\frac 1n \sum_{h\in H}h$ for any $g\in G$. By comparing coefficients in the both parts of the last relation, we see that $ghg^{-1}\in H$, i.e., the subgroup $H$ is normal.

Let $F_0$ be a prime subfield of the field $F$. We consider the
finite ring $F_0H$. By the Maschke theorem, it is isomorphic to some finite direct product of matrix rings over division rings; in addition, any finite division ring is a field by the Wedderburn theorem. We assume that the group $H$ is not commutative. Then one of the summands of the ring $F_0H$ is the matrix ring of order $k>1$ over some field; this is impossible, since such a matrix ring contains a non-central idempotent.

Now let $H$ be an arbitrary $q$-subgroup in $G$. We take any
element $h\in H$ and an arbitrary element $g\in G$. Since $h$
generates a cyclic $q$-subgroup $H_0=\langle h \rangle$, we have $ghg^{-1}\in H_0\subseteq H$ for any $g\in G$, i.e., the 
subgroup $H$ normal.

If $x,y\in H$, then the subgroup $H_1=\langle x,y \rangle$ is finite by Proposition 3.2.1 and the following Dicman's lemma:

\textsf{If $x_1,\ldots, x_n$ are elements of finite order of an arbitrary group $G$ and each of the elements $x_1,\ldots, x_n$ has only finite the number of conjugated elements, then there exists a finite normal subgroup $N$ of the group $G$ containing $x_1,\ldots,x_n$}.

By the first part of the proof, the subgroup $H_1$ is commutative, therefore, $xy=yx$.~$\rhd$

In the case of finite groups, we have a more strong assertion which reduces the study of centrally essential group algebras of finite groups to the study of centrally essential group algebras of finite $p$-groups.

\textbf{3.2.3. Proposition.}\\
Let $|G|=n<\infty$ and let $F$ be a field of characteristic $p>0$. Then the following conditions are equivalent.

\textbf{1)} The ring $FG$ is centrally essential.

\textbf{2)} $G=P\times H$, where $P$ is the unique Sylow $p$-subgroup of the group $G$, the group $H$ is commutative, and the ring $FP$ is centrally essential.

$\lhd$ Let $FG$ is centrally essential. By Lemma 3.2.2, every Sylow $q$-subgroup for $q\neq p$ is normal in $G$ and it is commutative; consequently, the product $H$ of all such subgroups is a commutative normal subgroup. Let $m=|H|$. We note that $(m,p)=1$, whence we have that the element $m$ is invertible in $F$.

We prove that the Sylow $p$-subgroup $P$ is normal in $G$.

We consider the following linear mapping $f\colon R\rightarrow R$: 
$$
f(r)=\frac 1 m \sum_{h\in H}hrh^{-1}.
$$
It is clear that $f(1)=1$ and $f(yry^{-1})=f(r)$ for any $y\in H$, since the left part and right part of the relation contain equal summands. Now we assume that $xy\neq yx$ for some $x\in P$ and $y\in H$. We set $r=x-yxy^{-1}$. It is directly verified that $r\neq 0$ but $f(r)=f(x)-f(yxy^{-1})=0$; this contradicts to 1.1.5. Therefore, elements of $P$ and $H$ commute, $G=PH$ and $P\cap H=\{1\}$; consequently, $G=P\times H$. By considering $FG$
as the group ring $(FP)H$, we obtain from Proposition 3.2.1 that $FP$ is centrally essential.

The converse assertion directly follows from 3.1.1 and the isomorphism $FG\cong (FP)H$.~$\rhd$

\textbf{3.2.4. Proposition.}\label{nilp2}\\
Let $G$ be a finite $p$-group and let $F$ be a field of characteristic $p$. If $\text{NC}(G)\le 2$, then the ring $FG$ is centrally essential.

$\lhd$ We recall that for any subgroup $H$ of $G$, we denote by $\omega H$ the right ideal of the ring $FG$ generated by the set $\{1-h|h\in H\}$; we also recall that this right ideal is an ideal if and only if the subgroup $H$ is normal. It is well known (e.g., see \cite[Lemma 3.1.6]{Pas77}) that the ideal $\omega G$ is nilpotent in our case.

Let $0 \neq x \in FG$. We consider all products $x(1-z)$, where $z\in Z=Z(G)$. If at least one of them (say, $x_1=x(1-z_1)$) is non-zero, then we consider the product $x_1(1-z)$ and so on. This process terminates at some step, i.e. there exists an integer $k\geq 0$ such that $x_k\neq 0$ but $x_k\omega Z=0$ (we assume that $x_0=x$). Then $x_k\in FG\Sigma_Z$ (see \cite[Lemma 3.1.2]{Pas77}). We note that $FG\Sigma_Z\subseteq Z(FG)$. Indeed, if $g,h\in G$, then 
$$
[g,h\Sigma_Z]=[g,h]\Sigma_Z=gh(1-h^{-1}g^{-1}hg)\Sigma_Z=0,
$$
since $h^{-1}g^{-1}hg\in G'\subseteq Z$. Therefore, by setting $c=(1-z_1)\ldots(1-z_k)$ (or $c=1$ for $k=0$), we obtain $c\in Z(FG)$ and $xc=x_k\in Z(FG)\setminus\{0\}$, which is required.~$\rhd$

\textbf{3.2.5. Lemma.}\label{conj}\\
Let $F$ be a field of characteristic $p$ and let $G$ be a finite $p$-group which satisfies the following condition:

$(*)$ for any element $g\in G\setminus Z(G)$, there exists a non-trivial subgroup $H\subseteq Z(G)$ such that $Hg\subseteq g^G$.

If $\text{NC}(G)>2$, then the ring $R=FG$ is not centrally essential.

$\lhd$ Let $K=g^G$ be a class of conjugated elements of the group $G$ such that $|K|>1$, and let $H$ be a subgroup which satisfies $(*)$. We note that $Hg'\subseteq K$ for any $g'\in K$, since $g'=a^{-1}ga$ for some $a\in G$ and $Hg'=Ha^{-1}ga=a^{-1}Hga\subseteq a^{-1}Ka=K$. Let $Hx_1,\ldots,Hx_t$ be all distinct cosets for $G$ with respect to $H$ contained in $K$. Then $K$ is a disjoint union of these cosets, whence we have $\Sigma_K=\sum_{i=1}^t \Sigma_{Hx_i}$. Now we note that $(\Sigma_Z)h=\Sigma_Z$ for any $h\in H$ since $H\subseteq Z$; therefore, $\Sigma_Z\cdot\Sigma_H=|H|\Sigma_Z=0$. Then we obtain that
$$
\Sigma_Z\Sigma_K=\sum_{i=1}^t\Sigma_Z\cdot\Sigma H\cdot x_i=0.\eqno{3.2.5(1)}
$$
Next, if $\text{NC}(G)>2$, then there exists an element $g\in G\setminus Z_2(G)$. This means that there exists an element $a\in G$ such that $(g,a)\in Z$. We consider element $x=g\Sigma_Z\neq 0$. We have
$$
[a,x]=[a,g\Sigma_Z]=(ag-ga)\Sigma_Z=ag(1-(g,a))\Sigma_Z\neq 0,
$$
since $1-(g,a)\in \omega Z$. Consequently, $x\in C=Z(R)$. With the use of the basis $(3.2.1)$, an arbitrary element $c\in C$ can be represented in the form $c=c_0+c_1$, where $c_0\in FZ$, $c_1=\sum_{i=0}^s\alpha_i\Sigma_{K_i}$, $K_1,\ldots,K_s$ are classes of conjugated elements of the group $G$, $|K_i|>1$ and $\alpha_i\in F$ for all $i=1,\ldots,s$. We assume that $xc\in Z(R)$. By $3.2.5(1)$, we have $xc_1=0$, whence we have $xc=xc_0$. Since $(\Sigma_Z)z=\Sigma_Z$ for any $z\in Z$, we obtain that $xc_0=\alpha x$ for some $\alpha\in F$. If $\alpha\neq 0$, then $x\in C$; this is a contradiction. We obtain that $xC\cap C=0$.~$\rhd$

\textbf{3.2.6. Lemma.}\label{centralizer}\\ 
Let $G$ be a group. If the centralizer $C_G(Z_2(G))$ of the subgroup $Z_2(G)$ is contained in $Z_2(G)$, then $G$ satisfies condition $(*)$ of Lemma 3.2.5.

$\lhd$ Let $g$ be an element of $G\setminus Z(G)$. We assume that there exists an element $a\in G$ such that
$$
(g,a)\in Z(G)\{1\}.\eqno{(*)}
$$ Let $z=(g,a)$. Then $gz=a^{-1}ga\in g^G$, whence we have $gz^k=a^{-k}ga^k\in g^G$ for any $k\geq 1$. Therefore, the subgroup $H$, generated by $z$, satisfies $(*)$.

Now we consider two cases. If $g\in Z_2(G)\setminus Z(G)$, then there exists an element $a\in G$ such that $(g,a)\neq 1$. However, it follows from the definition of $Z_2(G)$ that $(g,a)\in Z(G)$, whence $(*)$ is true.

It remains to consider the case $g\in Z_2(G)$. In this case, $g\in Z(Z_2(G))$, whence we have that there exists an element $a\in Z_2(G)$ such that $z=(g,a)\neq 1$. However $z\in Z_1(G)$, since $a\in Z_2(G)$, and we obtain $(*)$.~$\rhd$

\textbf{3.2.7. Remark; A.Yu. Olshansky.}\\
There exists another series of groups which satisfy conditions of Lemma 3.2.6. Namely, let $p$ be a prime integer and let $G$ be a free 3-generated group of the variety defined by identities $x^p=1$ and $(x_1,x_2,x_3,x_4)=1$. Then $G/G'$ is an elementary Abelian $p$-group; therefore, $G'$ is the Frattini subgroup of the group $G$. If $g\in G'$, then we can include $gG'$ in a system of free generators of the group $G/G'$; consequently, the element $g$ can be included in a system consisting of three generators of the group $G$. Since the group $G$ is finite, this generator system is free. Therefore, if $g\in C_G(G')$, then $G$ satisfies the identity $(x_1,x_2,x_3)=1$; this is impossible, since the group $G$ can be homomorphically mapped onto the group of upper uni-triangular matrices of order 4
over $GF(p)$ which does not satisfy this identity. Therefore, $Z_2(G)\supseteq G'\supseteq C_G(G')\supseteq C_G(Z_2(G))$.

\textbf{3.2.8. Proposition.}\\ 
If $F$ is a field of characteristic $p>0$, then there exists a group $G$ of order $p^5$ such that the group algebra $FG$ is not centrally essential.

$\lhd$ We construct of the group which satisfy conditions of Lemma 3.2.6. We consider cases $p=2$ and $p\neq 2$ separately.

Let $p=2$. We consider the direct product $N$ of the quaternion group $Q_8=\{\pm 1,\pm i,\pm j,\pm k\}$ and the cyclic group $\langle a \rangle$ of order 2 with generator $a$, we also consider the automorphism $\alpha$ of the group $N$ defined on generators by relations $\alpha(i)=j$, $\alpha(j)=i$, $\alpha(a)=(-1)a$. We set $\Gamma=\langle\alpha\rangle$. We have the semidirect product $G=N\rtimes\Gamma$ whose elements are considered as products $x\gamma$, where $x\in N$ and $\gamma \in \langle\alpha\rangle$ and the operation is defined by the relation $x\gamma x'\gamma'=x\gamma(x')\gamma\gamma'$. Elements of the form
$x\cdot 1$ are naturally identified with elements $x\in N$ and elements of the form $1\cdot\gamma$ identified with elements $\gamma\in \Gamma$. It is directly verified that $Z_1(G)=\langle-1\rangle$, $Z_2(G)=\langle k, a\rangle=C_G(Z_2(G))$.

Now we assume that $p>2$. We consider the semidirect product $N=A\rtimes \Gamma$ of the elementary Abelian group $A$ of order $p^3$ with generators $a,b,c$ and the cyclic group $\Gamma=\langle\gamma\rangle$, where $\gamma$ is an automorphism of the group $A$ defined on generators by relations 
$$
\gamma(a)=a,\ \gamma(b)=b,\ \gamma(c)=bc.
$$
It is directly verified that $|N|=p^4$ and any element of the group $N$ can be uniquely represented as the product $a^kb^lc^m\gamma^r$, where {$k,l,m,r\in\{0,\ldots,p-1\}$.} We prove that mapping 
$\beta\colon \{a,b,c,\gamma\}\rightarrow N$ defined by relations
$$
\beta(a)=a,\ \beta(b)=b,\ \beta(c)=ac,\ \beta(\gamma)=abc\gamma,
$$
can be extended to an automorphism $\hat\beta$ of the group $N$. Indeed, for any $k,l,m,r\in\mathbb{Z}$, we set
$$
\hat\beta(a^kb^lc^m\gamma^r)=a^kb^la^mc^ma^rb^r(c\gamma)^r=a^{k+m+r}b^{l+\frac{r (r+1)} 2}c^{m+r}\gamma^r.$$
This definition is correct, since $p|\frac{r(r+1)}2$ if $p|r$. It is directly verified that for any
$k,l,m,r,k',l',m',r'\in\{0,\ldots,p-1\}$, relations
$$
a^kb^lc^m\gamma^r\cdot a^{k'}b^{l'}c^{m'}\gamma^{r'}=a^{k+k'}b^{l+l'+rm'}c^{m+m'}\gamma^{r+r'}
$$
hold. Consequently,
$$
\hat\beta(a^kb^lc^m\gamma^r\cdot a^{k'}b^{l'}c^{m'}\gamma^{r'})=\\
$$
$$
=a^{k+k'+m+m'+r+r'}b^{l+l'+rm'+\frac{(r+r')(r+r'+1)}2}c^{m+m'+r+r'}\gamma^{r+r'}.
$$
On the other hand,
$$
\hat\beta(a^kb^lc^m\gamma^r)\cdot\hat\beta(a^{k'}b^{l'}c^{m'}\gamma^{r'})=
$$
$$
=(a^{k+m+r}b^{l+\frac{r(r+1)}2}c^{m+r}\gamma^r)\cdot(a^{k'+m'+r'}b^{l'+\frac{r'(r'+1)}2}
c^{m'+r'}\gamma^{r'})=
$$
$$
= a^{k+m+r+k'+m'+r'}b^{l+\frac{r(r+1)}2+l'+\frac{r'(r'+1)}2+r(m'+r')}c^{m+r+m'+r'}\gamma^{r+r'}.
$$
It remains to note that we have the following identity
$$
l+\frac {r(r+1)}2+l'+\frac{r'(r'+1)}2+r(m'+r')=
$$
$$
=l+l'+rm'+rr'+\frac{r^2+r+{r'}^2+r'} 2= l+l'+rm'+$$
$$
+\frac{r^2+2rr'+r'^2+r+r'}2=
l+l'+rm'+\frac{(r+r')(r+r'+1)} 2
$$
Now we set $G=N\rtimes \langle\beta\rangle$. It is directly verified that $Z_1(G)=\langle a,b \rangle$ and $Z_2(G)=\langle a,b,c \rangle=C_G(Z_2(G))$.~$\rhd$

\textbf{3.2.9. Theorem.}\\
Let $F$ be a field of characteristic $p>0$.

\textbf{a.} If $G$ is an arbitrary finite group, then the group algebra $FG$ is a centrally essential ring if and only if $G=P\times H$, where $P$ is the unique Sylow $p$-subgroup of the group $G$, the group $H$ is commutative, and the ring $FP$ is centrally essential.

\textbf{b.} If $G$ is a finite $p$-group and the nilpotence class\footnote{It is well known that every finite $p$-group is nilpotent, e.g., see \cite[Theorem 10.3.4]{Hal59}.} of the group $G$ does not exceed $2$, then the group algebra $FG$ is a centrally essential ring.

\textbf{c.} There exists a group $G$ of order $p^5$ such that 
the group algebra $FG$ is centrally essential.

$\lhd$ Theorem 3.2.9 follows from Proposition 3.2.3, Proposition 3.2.4 and Proposition 3.2.8.~$\rhd$

\textbf{3.2.10. Remarks.}\\
\textbf{a.} It is well known that a group ring $AG$ is semiprime if and only if the ring $A$ is semiprime and orders of finite normal subgroups of the group $G$ are not zero-divisors in $A$.

\textbf{b.} If $A$ is a semiprime ring such that
its additive group is torsion-free and $G$ is an arbitrary
group, then the group ring $AG$ is centrally essential if and only if the ring $A$ and the group $G$ are commutative. Indeed, by Theorem 1.2.2, any centrally essential semiprime ring is commutative. Therefore, Remark \textbf{b} follows from Remark \textbf{a}.

\textbf{c.} 
Let $F$ be an arbitrary field of zero characteristic and let $G$ be a group $G$. In connection to Theorem 3.2.9, we note that the group algebra $FG$ is centrally essential if and only if the algebra $FG$ is commutative; see Remark \textbf{b}. Therefore, only the case of fields of positive characteristic is of interest under the study of centrally essential group algebras over fields.

\textbf{d.} 
Let $G$ be a finite $p$-group of nilpotence class 3.
In connection to Theorem 3.2.9(c), we note that group rings of $G$ can be centrally essential and can be not centrally essential. More precisely, we used computer algebraic system GAP \cite{GAP17} to verify the property that for any group of order 16 and nilpotence class 3, its group algebra over a field $\text{GF}\,(2)$ is centrally essential.

\textbf{e.} There exists a finite $2$-group $G$ such that the group algebra $R=FG$ over the field $F$ of order 2 is centrally essential and contains an element $x$ such that $x^2=0$ but $xRx\ne 0$.\\ 
$\lhd$ Let $G=D_4$ be the dihedral group of order 8 defined by generators $a,b$ and defining relations $a^4=b^2=(ab)^2=1$. It is easy to verify that 
$$
G=\{1,a,a^2,a^3,b,ab,a^2b,a^3b\},\;
G'=Z(G)=\langle a^2 \rangle.
$$
Therefore, the group algebra $FG$ is centrally essential.\\
In the same time, $(1+b)^2=1+b^2=1+1=0$ and
$$(1+b)a(1+b)=a+ba+ab+bab=1+a^3+ab+a^3b\neq 0.$$

\subsection{Rings of Fractions, Group and Semigroup Rings}\label{subsection3.3}

This subsection is based on \cite{LT21b} and \cite{LT23}.

\subsubsection{Rings of Fractions and Group Rings}\label{subsubsection3.3.1}

\textbf{3.3.1. Notation and used information.}\\
For a fixed group $G$, we denote by $P(G)$ and $G_p$ the torsion part $G$ and the set of elements of the group $G$ whose orders are degrees of prime integer $p$, respectively. In addition, $Z(G)$ is the center of the group $G$, $K$ is a field of characteristic $p > 0$, and $KG$ is the group algebra of the group $G$ over $K$.

For an $FC$-group $G$, it is known that $G_p$ is a characteristic subgroup in $G$; see, for example \cite[Lemma 8.1.6]{Pas77}.
In Theorem 3.2.9(a), it is proved that if a group $G$ is finite, then $KG$ is a centrally essential ring if and only if $G=G_p\times H$, where $G_p$ is the unique Sylow $p$-subgroup in $G$ and the ring $KG_p$ is centrally essential. 

If $G$ is a finitely generated $FC$-group, then the torsion part $P(G)$ is a finite normal subgroup in $G$. In addition, $G/P(G)$ is a finitely generated Abelian torsion-free group; see, for example \cite[Lemma 4.1.5]{Pas77}. Let 
$$
G/P(G) = \langle\overline{a_1}\rangle\times\ldots \times \langle\overline{a_n}\rangle, \text{ where } \langle\overline{a_i}\rangle = a_iP(G),\;i = 1, 2, \ldots, n.
$$
Then $G$ is a semidirect product $G = P(G)\rtimes F$, where $F = A_1 \times \ldots \times A_n$, $A_i = <a_i>$ is an infinite cyclic group, $i = 1, 2, \ldots, n$; see \cite[Theorem 4]{Nis57}.

\textbf{3.3.2. Remark.}\label{remark-1.1}\\
Let $K$ be a field of zero characteristic, $G$ be any group, and let the group algebra $KG$ be centrally essential. In 3.2.10(c), it is proved that $KG$ is commutative.\\
Let $R$ be a commutative ring and let $G$ be a torsion-free group such that the group ring $RG$ is centrally essential. Then the ring $RG$ is commutative. Indeed, by Proposition 3.2.1, the ring $R$ also is centrally essential and all classes of conjugated elements in $G$ are finite. By \cite[Lemma 4.1.6]{Pas77}, the group $G$ is Abelian, the group ring $RG$ is commutative and, consequently, $RG$ has commutative classical ring of fractions. 

\textbf{3.3.3. Remark.}\label{proposition-2.1}\\
If the group $G$ does not have elements of order $p$ and the group algebra $KG$ is centrally essential, then $KG$ is commutative.
 
$\lhd$ It follows from \cite[Theorem 4.2.13]{Pas77} that $KG$ is a semiprime algebra. Consequently, the centrally essential semiprime algebra $KG$ is commutative by Theorem 1.2.2.~$\rhd$

\textbf{3.3.4. Proposition.}\label{proposition-3.1}\\
Let $R$ be a ring. If for any non-zero-divisor $b$, there exists a non-zero-divisor $x$ such that $bx = y \in Z(R)$ (resp., $xb = y\in Z(R)$), then $R$ has the right (resp., left) classical ring of fractions.

$\lhd$ 
Let $a, b\in R$, where $b$ is a non-zero-divisor in $R$. Then $b(xa) = a(bx) = ay$. Therefore, the ring $R$ satisfies the right \"Ore condition, $(ax)b = a(xb) = (xb)a$, and the ring $R$ satisfies the left \"Ore condition.~$\rhd$

\textbf{3.3.5. Corollary.}\label{corollary-3.2}\\
Any centrally essential group algebra has the two-sided classical ring of fractions.

$\lhd$ Since the group $G$ is an $FC$-group, it follows from \cite[Lemma 4.4.4]{Pas77} that for any non-zero-divisor $b$, there exists a non-zero-divisor $x\in KG$ such that $xb = y\in Z(KG)$ ($bx = y\in Z(KG)$) and $y$ is a non-zero-divisor in $KG$.~$\rhd$

\textbf{3.3.6. Remark.}\label{remark-3.3}\\
Corollary 3.3.5 also follows from \cite{HerS71}, since all classes of conjugated elements in $G$ are finite. 

\textbf{3.3.7. Proposition.}\label{proposition-2.2}\\
Let $G$ be a finitely generated group. If $KG$ is a centrally essential ring, then $G$ is an $FC$-group and $KG_p$ is a centrally essential ring. If $F\subseteq Z(G)$ (under the above notations), then the converse is true, as well.

$\lhd$ By Proposition 3.2.1, the group $G$ is an $FC$-group. As was mentioned above, $KG_p$ is a centrally essential ring provided $KP(G)$ is a centrally essential ring. Consequently, without loss of generality, it is sufficient to prove that $KP(G)$ is a centrally essential ring. 

By assumption, for $0\neq \alpha\in KP(G)$, there exist non-zero central elements $\beta, \gamma\in KG$ such that $\alpha\beta = \gamma$. 
Let $\pi(\gamma)\neq 0$, where $\pi\colon KG \to KP(G) $ is a natural projection defined by the relation $\pi(\sum_{x\in G}a_xx) = \sum_{x\in P(G)}a_xx$. Let $\mu\in KP(G)$. Then it follows from \cite[Lemma 1.1.2]{Pas77} that
$$
\mu\pi(\beta) = \pi(\mu\beta) = \pi(\beta \mu) = \pi(\beta)\mu.
$$
Therefore, $\pi(\beta)\in Z(KP(G))$. Next,
$$
\mu\alpha\pi(\beta) = \pi(\mu\alpha\beta) = \pi(\alpha\beta \mu) = \alpha\pi(\beta)\mu.
$$
Therefore, $\alpha\pi(\beta)\in Z(KP(G))$ and $\alpha\pi(\beta) = \pi(\alpha\beta) = \pi(\gamma)\neq 0$.

Let $\pi(\gamma) = 0$. We verify that $\gamma$ is a non-zero-divisor in $KG$. It follows from \cite[Lemma 6]{Bro78} that $\gamma$ is a non-zero-divisor in $KG$ if and only if $\gamma$ is a non-zero-divisor in the ring $Z(KG)\cap KH$, where $H = \langle\text{supp}\,\gamma\rangle$. Since the group $H$ is a finitely generated normal torsion-free subgroup in $G$, we have that $H$ has a finitely generated central torsion-free subgroup of finite index; see \cite[Lemma 4.1.8]{Pas77}. It follows from \cite[Corollary 2]{Bro76} that $KH$ does not have zero-divisors. Consequently, elements $\gamma$ and $\beta$ are non-zero-divisors in $KG$. Since $KG$ has the classical ring of fractions $Q_{\text{cl}}(KG)$, there exists a $\beta^{-1}\in Q_{\text{cl}}(KG)$. Next, $\beta\in Z(KG)$ and $Z(Q_{\text{cl}}(KG)) = Q_{\text{cl}}(Z(KG))$; see \cite[Theorem 4.4.5]{Pas77}. Then $\beta^{-1}\in Z(Q_{\text{cl}}(KG))$. Since $\alpha\beta = \gamma$, we have $\alpha = \gamma\beta^{-1}\in Z(KP(G))$.

Conversely, we take the set $A_1\times \ldots \times A_n = F$ as a transversal for the subgroup $P(G)$. It follows from \cite[Lemma 1.1.4]{Pas77} that $F$ 
is a basis of the algebra $KG$ over $KP(G)$. Since $F\subseteq Z(G)$ and the algebra $KP(G)$ is a centrally essential ring, the algebra $KG$ also is a centrally essential ring by Proposition 1.1.7.~$\rhd$

\textbf{3.3.8. Example.}\label{example-3.4}\\
Let $K$ be a field. We consider the ring $\mathcal{R}$ of all $3\times 3$ matrices of the form
$$
A = 
\left(\begin{matrix}
k & a & b\\
0 & k & a\\
0 & 0 & k\\
\end{matrix}\right),
$$ 
where $k\in K$, $a$ and $b$ are contained in the polynomial ring $K\langle x, y\rangle$ in two non-commuting variables $x$ and $y$ over the field $K$ with relations $xk = kx$ and $ky = yk$, where $k\in K$, and $yx - xy = x$; for example, see \cite{Eli69}. 

We note that the ring $\mathcal{R}$ is not commutative. Indeed, 
$$
\left(\begin{matrix}
0 & x & 0\\
0 & 0 & x\\
0 & 0 & 0\\
\end{matrix}\right)
\left(\begin{matrix}
0 & y & 0\\
0 & 0 & y\\
0 & 0 & 0\\
\end{matrix}\right) = \left(\begin{matrix}
0 & 0 & xy\\
0 & 0 & 0\\
0 & 0 & 0\\
\end{matrix}\right)\neq 
\left(\begin{matrix}
0 & 0 & yx\\
0 & 0 & 0\\
0 & 0 & 0\\
\end{matrix}\right) =
$$
$$ 
=\left(\begin{matrix}
0 & y & 0\\
0 & 0 & y\\
0 & 0 & 0\\
\end{matrix}\right)
\left(\begin{matrix}
0 & x & 0\\
0 & 0 & x\\
0 & 0 & 0\\
\end{matrix}\right).
$$
Next, 
$$
Z(\mathcal{R}) = \left\lbrace\left.\left(\begin{matrix}
k & k' & h\\
0 & k & k'\\
0 & 0 & k \\
\end{matrix}\right) \right| k, k'\in K; h\in K\langle x, y\rangle\right\rbrace.
$$ 
In addition, if $f\in K\langle x, y\rangle$, then
$$
\left(\begin{matrix}
0 & f & 0\\
0 & 0 & f\\
0 & 0 & 0\\
\end{matrix}\right)
\left(\begin{matrix}
0 & 1 & 0\\
0 & 0 & 1\\
0 & 0 & 0\\
\end{matrix}\right) = \left(\begin{matrix}
0 & 0 & f\\
0 & 0 & 0\\
0 & 0 & 0\\
\end{matrix}\right)\in Z(\mathcal{R}).
$$
Consequently, $\mathcal{R}$ is a centrally essential ring. 
For any regular matrix 
$$
A = \left(\begin{matrix}
k & a & b\\
0 & k & a\\
0 & 0 & k\\
\end{matrix}\right),
$$
there exists a regular matrix 
$$
A' = \left(\begin{matrix}
k' & a' & 0\\
0 & k' & a'\\
0 & 0 & k'\\
\end{matrix}\right),
$$
where $k'\neq 0$, $a' = \cfrac{1}{k}(k'' - ak')$ for some $0\neq k''\in K$ such that $AA'\in Z(\mathcal{R})$. It follows from Proposition 3.3.4 that $\mathcal{R}$ is a non-commutative centrally essential ring which has the two-sided classical ring of fractions.

\textbf{3.3.9. Proposition.}\label{proposition-3.5}\\
Let $R$ be a centrally essential ring and let $\alpha$ be a non-zero-divisor in $Z(R)$. Then $\alpha$ is a non-zero-divisor in $R$.

$\lhd$
Let $\alpha\beta = 0$ for some $0\neq \beta\in R$. Then $\beta\notin Z(R)$ and there exist elements $c, d\in Z(R)$ such that $0\neq \beta c = d$. Since $\alpha$ is a non-zero-divisor in $Z(R)$, we have that $\alpha d\neq 0$ and $\alpha\beta c\neq 0$. This is a contradiction.~$\rhd$

\textbf{3.3.10. Proposition.}\label{proposition-3.6}\\
Let $R$ be a centrally essential ring and $R$ has the classical ring of fractions. Then 
$Q_{\text{cl}}(Z(R))\subseteq Z(Q_{\text{cl}}(R))$.

$\lhd$ Let $\alpha\in Z(R)$ be a non-zero-divisor in $Z(R)$. It follows from Proposition 3.3.9 that $\alpha$ is a non-zero-divisor in $R$. Consequently, there exists an element $\alpha^{-1}\in Q_{\text{cl}}(R)$. We verify that $\alpha^{-1}\in Z(Q_{\text{cl}}(R))$. 

Let $\beta = \gamma\delta^{-1}\in Q_{\text{cl}}(R)$. Since $\alpha\gamma = \gamma\alpha$, we have that $\gamma = \alpha^{-1}\gamma\alpha$ in the ring $Q_{\text{cl}}(R)$. Then
$$
(\gamma\delta^{-1})\alpha^{-1} = \gamma(\alpha\delta)^{-1} = \gamma(\delta\alpha)^{-1} = 
\gamma\alpha^{-1}\delta^{-1} = 
$$
$$
=\alpha^{-1}\gamma\alpha\alpha^{-1}\delta^{-1} = \alpha^{-1}(\gamma\delta^{-1}).
$$
Therefore, $\alpha^{-1}\in Z(Q_{\text{cl}}(R))$.
If $\alpha\in Z(R)$ is a zero-divisor, then it follows from relations $\alpha\delta = \delta\alpha$ and $\alpha = \delta^{-1}\alpha\delta$ that
$$
\alpha(\gamma\delta^{-1}) = (\gamma\alpha)\delta^{-1} 
= \gamma\delta^{-1}\alpha\delta\delta^{-1} = (\gamma\delta^{-1})\alpha.
$$
Therefore, $Q_{\text{cl}}(Z(R))\subseteq Z(Q_{\text{cl}}(R))$. 
The left-sided analogue is similarly verified.~$\rhd$

\textbf{3.3.11. Theorem.}\label{theorem-1.2}\\
Every centrally essential group algebra over any field has two-sided classical ring of fractions. In addition, the group algebra over a field is centrally essential if and only if it has right classical ring of fractions which is a centrally essential ring.

$\lhd$ The first assertion of the theorem follows from Corollary 3.3.5. Let $0\neq as^{-1}\in Q_{\text{cl}}(KG)$, where $a, s\in KG$ and $s$ is a non-zero-divisor. Since $G$ is an $FC$-group, it follows from \cite[Lemma 4.4.4]{Pas77} that there exists a non-zero-divisor $\gamma\in KG$ such that $s\gamma = t\in Z(KG)$ and $t$ is a non-zero-divisor in $KG$. Therefore, we have that $s^{-1} = \gamma t^{-1}$, in the ring $Q_{\text{cl}}(KG)$. 
By assumption, for the non-zero element $a\gamma\in KG$, there exist non-zero elements $c, d\in Z(KG)$ such that $(a\gamma)c = d$.
By Proposition 3.3.10 (also see \cite[Theorem 4.4.5]{Pas77}), any element of $Z(KG)$ is central in $Q_{\text{cl}}(KG)$, i.e., $c, d, t^{-1}\in Z(Q_{\text{cl}}(KG))$. Then
$$
0\neq (as^{-1})c = (a\gamma t^{-1})c = (a\gamma c)t^{-1} = dt^{-1}\in Z(Q_{\text{cl}}(KG)).
$$
Conversely, let $0\neq r\in KG$. By assumption, there exist elements $t, s\in Z(Q_{\text{cl}}(KG))$ such that $0\neq rt = s$. Since $Z(Q_{\text{cl}}(KG)) = Q_{\text{cl}}(Z(KG))$, we have that $t = cd^{-1}$ and $s = mn^{-1}$ for some $c, d, m, n\in Z(KG)$. Then the relation 
$rcd^{-1} = rt = s = mn^{-1}$ implies that $rc = mn^{-1}d$ and 
$$
r(cn) = (rc)n = md\in Z(KG).
$$
In addition, $md\neq 0$, since $d$ is a non-zero-divisor in $KG$.

\textbf{3.3.12. Remarks.}
Let $G$ be a group,
$$
\Delta(G) = \{x\in G: |G : C_G(x)| < \infty\}
$$
(i.e., $\Delta(G)$ is an $FC$-subgroup in $G$), and let
$$
\Delta^{+}(G) = \{x\in \Delta(G)\,|\, o(x) < \infty\}.
$$
It is well known that $\Delta(G)$ and $\Delta^{+}(G)$ are characteristic subgroups in $G$; see details in \cite{Pas77}. If the group $\Delta^{+}(G)$ is finite, then the ring $Q_{\text{cl}}(K\Delta(G))$ and is a quasi-Frobenius ring; in particular, it coincides with maximal ring of fractions $Q_{\text{max}}(K\Delta(G))$; see \cite{Bro78}. It follows from these facts and theorem 3.3.11 we obtain the following remark.

\textbf{a.} \label{remark-5.1}
If the subgroup $\Delta^{+}(G)$ of the group $G$ is finite, then the following conditions are equivalent.
\begin{itemize}
\item
$KG$ is a centrally essential ring.
\item
$Q_{\text{cl}}(KG)$ is a centrally essential ring.
\item
$Q_{\text{max}}(KG)$ is a centrally essential ring.
\end{itemize}

\textbf{b.}\label{remark-5.2}
A ring $A$ is called a \textsf{ring with large center}\label{larcen} if any non-zero ideal of the ring $A$ has the non-zero intersection with the center of the ring $A$. In \cite[Theorem 2]{Zly14} is proved that if $R$ is a ring with large center, then $Q_{\text{max}}(Z(R))\subseteq Z(Q_{\text{max}}(R))$. Since it is clear that any centrally essential ring is a ring with large center, the assertion of Proposition 3.3.10 remains true for maximal rings of fractions, as well.

\subsubsection{Rings of Fractions and Semigroup Rings}\label{subsubsection3.3.2}

In this subsection, $F$ is a field, $S$ is a semigroup, $F$ is a field, and $FS$ is the semigroup $F$-algebra of the semigroup $S$. The centers of the semigroup $S$ and the semigroup algebra $FS$ are denoted by $Z(S)$ and $Z(FS)$, respectively. If $a = \sum\alpha_ss\in FS$, then $\text{supp }(a) = \{s\in S \, | \, \alpha_s\neq 0\}$. 

\textbf{3.3.13. Remarks.}\\
\textbf{a.} A semigroup $S$ is called a \textsf{left cancellative semigroup} if $a = b$ for any $a, b, c\in S$ with $ca = cb$. A \textsf{right cancellative semigroup} is defined dually. A right and left cancellative semigroup is called a \textsf{cancellative semigroup}. It is well known that a torsion cancellative semigroup is a group; e.g., see \cite{CliP61}. A cancellative semigroup is embedded in the right group of fractions if and only if the intersection of any two principal right ideals of the semigroup $S$ is non-empty, i.e., $sS\cap tS\neq \emptyset$ for all $s, t\in S$ (the right \"Ore condition). If $S$ also satisfies the left \"Ore condition which is symmetrically defined, then the group $G_S = SS^{-1} = S^{-1}S$ is called the \textsf{group of fractions} of the semigroup $S$. Any element of the group $G_S$ can be written in the form $a^{-1}b$ and in the form $cd^{-1}$, where $a, b, c, d\in S$. 

\textbf{b.} We recall that the subgroup $\Delta(G)$ of the group $G$ and properties of $\Delta(G)$ are considered in Remarks 3.3.12.

\textbf{c.} Let $S$ be a cancellative semigroup and $s\in S$. If for some $x\in S$, there exists an element $t\in S$ such that $xs = tx$, then the element $t$ is uniquely defined; it is denoted by $s^x$. Then $\Delta(S)$ is the set of elements $s\in S$ such that elements $s^x$ are defined for all $x\in S$ and the set 
$\{s^x \, | \, x\in S\}$ is finite. If $s\in \Delta(S)$, then we set $D_S(s) = \{s^x \, | \, x\in S\}$. It is clear that if $S$ is embedded in the group of fractions $G_S$, then for $s\in \Delta(S)$, the set $D_S(s)$ is embedded in the set of conjugated elements for $s$ in $G_S$. If $S$ is a cancellative semigroup, then 
$Z(FS)$ is an $F$-subspace in $FS$ generated by elements of the form $\sum_{t\in D_S(s)}t$, where $s\in \Delta(S)$; see \cite[Theorem 9.10]{Okn91}.

\textbf{3.3.14. Proposition.} Let $S$ be a cancellative semigroup. If the semigroup algebra $FS$ is a centrally essential ring, then $S = \Delta(S)$. 

$\lhd$ Let $s\in S$. By assumption, $0\neq cs = d$ for some $c, d\in Z(FS)$. For any $y\in \text{supp }(d)$, there exists an element $x\in \text{supp }(c)$ such that $xs = y$. It follows from \cite[Proposition 9.2(iii)]{Okn91} that $x, y\in \Delta(S)$. In addition, $\Delta(S)$ is a right and left \"Ore set in $S$ and $G_{\Delta(S)} = \Delta(S)^{-1}\Delta(S) = \Delta(S)\Delta(S)^{-1}$ is an $FC$-group; see \cite[Corollary 9.6 and Proposition 9.8(iii)]{Okn91}. Consequently, $s = x^{-1}y$, where $x^{-1}\in \Delta(S)^{-1}$, $y\in \Delta(S)$. For any $t\in S$, we have $x^t\in \Delta(S)$. Therefore, it follows from $tx = x^tt$ that $(x^t)^{-1}t = tx^{-1}$, i.e., $(x^t)^{-1} = (x^{-1})^t$ in the group $G_{\Delta(S)}$. Then the element $s^t = (x^{-1}y)^t = (x^{-1})^ty^t$ there exists a for any $t\in S$; see \cite[Basic property (a), p.108]{Okn91}. Next,
$$
\{s^t\, | \, t\in S\} = \{(x^{-1}y)^t\, | \, t\in S\} = \{(x^{-1})^t\, | \, t\in S\}\cdot \{y^t\, | \, t\in S\} = 
$$
$$
\{(x^t)^{-1}\, | \, t\in S\}\cdot \{y^t\, | \, t\in S\}.
$$
The first set is finite, since is finite the set $\{x^t\, | \, t\in S\}$. It follows from $y\in \Delta(S)$ that the second set is finite. Therefore, the set $D_S(s)$ is finite and $s\in\Delta(S)$.~$\rhd$

\textbf{3.3.15. Corollary.} If $FS$ is a centrally essential semigroup algebra of the cancellative semigroup $S$, then $S$ has the group of fractions $G_S$. 

$\lhd$ By Proposition 3.3.14, we have $S = \Delta(S)$. Since $\Delta(S)$ is a right and left \"Ore set, $S$ has the group of fractions $G_S$.~$\rhd$

By Corollary 3.3.15, under the study of centrally essential semigroup algebras of cancellative semigroups, it is sufficient to consider semigroups $S$ which have the group of fractions $G_S$.

\textbf{3.3.16. Corollary.} Let $F$ be a field with $\text{char }F = 0$. Then any centrally essential semigroup $F$-algebra of a cancellative semigroup is commutative.

$\lhd$ The algebra $FS$ is semiprime if and only if the algebra $FG_S$ is semiprime; see \cite[Theorem 7.19]{Okn91}. It is well known that the group algebra over a field of characteristic $0$ is semiprime; e.g., see \cite[Theorem 4.2.12]{Pas77}. By Theorem 1.2.2, all centrally essential semiprime rings are commutative.~$\rhd$

\textbf{3.3.17. Example.} We consider the subring $\mathcal{R}$ of the ring $M_{7}(F)$ of all matrices of order 7 over a field $F$ of characteristic $0$ consisting of matrices of the form 
$$
\left(\begin{matrix}
\alpha & a & b & c & d & e & f\\
0 & \alpha & 0 & b & 0 & 0 & d\\
0 & 0 & \alpha & 0 & 0 & 0 & e\\
0 & 0 & 0 & \alpha & 0 & 0 & 0\\
0 & 0 & 0 & 0 & \alpha & 0 & a\\
0 & 0 & 0 & 0 & 0 & \alpha & b\\
0 & 0 & 0 & 0 & 0 & 0 & \alpha\\
\end{matrix}\right).
$$
Then $\mathcal{R}$ is a non-commutative centrally essential ring by 3.6.8  below.  
Let $e_{\alpha} = E_7$, $e_a$, $e_b$, $e_c$, $e_d$, $e_e$, $e_f$ be matrices, in which only non-null entry with value $1$ in places $\alpha$, $a$, $b$, $c$, $d$, $e$, $f$, respectively. We consider the semigroup $S = \{e_{\alpha}, e_a, e_b, e_c, e_d, e_e, e_f\}\cup \{\theta\}$, where $\{\theta\}$ acts as zero. Then $\mathcal{R}\cong F_0S$, where $F_0S$ is the compressed semigroup algebra of the semigroup $S$ over the field $F$.
Since $FS\cong F\bigoplus F_0S$ (e.g., see \cite[Corollary 4.9]{Okn91}), then $FS$ is a centrally essential semigroup algebra as a direct sum of centrally essential algebras.
 
\textbf{3.3.18. Theorem.}\\
\textbf{a.} Let $S$ be a cancellative semigroup and let $F$ be a field. The semigroup $F$-algebra $FS$ is centrally essential if and only if the group of fractions $G_S$ of the semigroup $S$ exists and the group algebra $FG_S$ of the group $G_S$ is a centrally essential group algebra.

\textbf{b.} There exist non-commutative centrally essential semigroup algebras over fields of characteristic zero (in addition, it is known that centrally essential group algebras over fields of characteristic $0$ are commutative).

$\lhd$
\textbf{a.} Let $FS$ be a centrally essential ring and let $0\neq a\in FG_S$, $a = \sum_{i=1}^n\alpha_ig_i$, where $\alpha_i\in F$, $g_i\in G_S$. It follows known that $G_S = SZ(S)^{-1}$; see \cite[Proposition 9.8(iv)]{Okn91}. Then $a = \sum_{i=1}^n\alpha_is_it_i^{-1}$ for some $s_i\in S$, $t_i\in Z(S)$, $i = 1,\ldots, n$. We set $a' = \alpha_1s_1t_2\ldots t_n +\ldots + \alpha_ns_nt_1\ldots t_{n-1}\in FS$. We note that $a'\neq 0$; it is sufficient to verify that $s_1t_2\ldots t_n,\ldots, s_nt_1\ldots t_{n-1}$ are distinct elements in $FS$. Indeed, if $s_it_1\ldots\widehat{t_i}\ldots t_n = s_jt_1\ldots\widehat{t_j}\ldots t_n$ for $i\neq j$, then we multiply this relation by $(t_1\ldots t_n)^{-1}$ and obtain $s_it_i^{-1} = s_jt_j^{-1}$, i.e., $g_i = g_j$. This is a contradiction. By assumption, $0\neq a'c' = d'$ for some $c', d'\in Z(FS)$. Then $0\neq ac'' = d'$, where $c'' = t_1\ldots t_nc'$ and $d'$ are central elements in $FS$ which remain central in $FG_S$; see \cite[Corollary 9.11(i)]{Okn91}.

Conversely, let $0\neq a\in FS$, $a = \sum\alpha_is_i$, where $\alpha_i\in F$, $s_i\in S$. By assumption, $0\neq ac = d$ for some $c, d\in Z(FG_S)$, $c = \sum_{i=1}^n\beta_ig_i$, $d = \sum_{j=1}^m\gamma_jh_j$, where $g_i, h_j\in G_S$. Let $g_i = x_iy_i^{-1}$, $h_j = z_jt_j^{-1}$, $x_i, z_j\in S$, 
$y_i, t_j\in Z(S)$, $i = 1,\ldots, n$, $j = 1,\ldots, m$. We set $y = y_1\ldots y_n$, $t = t_1\ldots t_m$, $c' = cyt\in Z(FS)$. Then
$$
ac' = (ac)yt = dyt\in Z(FS).
$$ 
It remains to verify that $ac'\neq 0$. We have 
$$
dyt = \gamma_1z_1yt_2\ldots t_m + \ldots + \gamma_mz_myt_1\ldots t_{m-1}.
$$
If $i\neq j$ and $z_iyt_1\ldots\widehat{t_i}\ldots t_m = z_jyt_1\ldots\widehat{t_j}\ldots t_m$, then 
$z_it_i^{-1} = z_jt_j^{-1}$ and $g_i = g_j$; this is a contradiction. 

\textbf{b.} The assertion follows from Example 3.3.17.~$\rhd$

\textbf{3.3.19. Example.} Let $S = \langle x,y,z \,|\, z\in Z(S), z^2 = e, xy = zyx\rangle$. It is directly verified that $S$ is a cancellative semigroup which has the group of fractions $G_S=\langle x, y, z | \, z\in Z(G_S), z^2 = e, x^{-1}y^{-1}xy = z\rangle$. Since $z$ is a central involution and $x^2, y^2\in Z(G_S)$, we have that the unique non-trivial commutator in $G_S$ is $x^{-1}y^{-1}xy$. Therefore, the commutant $G'_S = <z>$. We have 
$Z(G_S) = <x^2, y^2, z>$. Let $H = G'_S = \{1, z\}$, $\hat{H} = 1 + z$, $\text{char }F = 2$ and $0\neq \alpha\in FG_S$. If $\alpha(1 + z) = 0$, then 
$\alpha\in FG_S\hat{H}$; see \cite[Lemma 3.1.2]{Pas77}.
Since $G_S$ is the class 2 nilpotent group, then $H\subseteq Z(G_S)$ and for $g_1, g_2\in G$ we have:
$$
[g_1, g_2\hat{H}] = [g_1, g_2]\hat{H} = g_1g_2(1 - g_2^{-1}g_1^{-1}g_2g_1)\hat{H} = 0.
$$
Thus, $\alpha\in Z(FG_S)$. If $\alpha_1 = \alpha(1 + z)\neq 0$, then $\alpha_1(1 + z) = 0$ and 
$\alpha(1 + z) = \alpha_1\in Z(FG_S)\backslash\{0\}$. Consequently, the group algebra $FG_S$ is centrally essential. By Theorem 3.3.18., 
 the semigroup algebra $FS$ is centrally essential, as well. 

\textbf{3.3.20. Example.} Let $S = \langle x, y, z\, | \, z\in Z(S),\, xy = zyx\rangle$. The semigroup $S$ has the group of fractions $G_S$ which is a free nilpotent group of nilpotence class 2; see \cite[Example 21]{Okn91}. It follows known that if the group does not contain of elements of order $p$, then centrally essential group algebra is commutative; see \cite[Proposition 1]{LT21b}. Therefore, the group algebra $FG_S$ is not centrally essential. By Theorem 3.3.18, the semigroup algebra $FS$ also is not a centrally essential.

\textbf{3.3.21. Lemma.} Let $FS$ be a centrally essential semigroup algebra of the cancellative semigroup $S$. Then for every non-zero-divisor $b\in FS$, there exists a non-zero-divisor $z\in FS$ such that $bz\in Z(FS)$. 

$\lhd$ It follows from \cite[Lemma 4.4.4]{Pas77} that there exists a non-zero-divisor $x\in FG_S$ such that $bx = y\in Z(FG_S)$.
If $x = \sum_{i=1}^n\alpha_is_it_i^{-1}$, where $\alpha_i\in F$, $s_i\in FS$, $t_i\in Z(FS)$ ($i = 1,2, \ldots, n$), then the element $z = xt_1\ldots t_n$ is a non-zero-divisor in $FS$ and $bz\in Z(FS)$.~$\rhd$

\textbf{3.3.22. Proposition.} If $FS$ is a centrally essential semigroup algebra of the cancellative semigroup $S$, then $FS$ has the classical ring of fractions.

$\lhd$ The assertion follows from Proposition 3.3.4, Lemma 3.3.21 and the property that their left-side analogues are true.~$\rhd$

The following theorem extends Theorem 3.3.11 to semigroup algebras of cancellative semigroups.

\textbf{3.3.23. Theorem.} A semigroup algebra of a cancellative semigroup is centrally essential if and only if it has the right classical ring of fractions which is a centrally essential ring.

$\lhd$ Let $FS$ be a centrally essential ring and $0\neq as^{-1}\in Q_{\text{cl}}(FS)$, where $s$ is a non-zero-divisor in $FS$. Let a non-zero-divisor $\gamma\in FS$ be such that $s\gamma = t\in Z(FS)$. Then $s^{-1} = \gamma t^{-1}$ in the ring $Q_{\text{cl}}(FS)$. By assumption, for the element $a\gamma\in FS$, there exist non-zero elements $c, d\in Z(FS)$ such that $0\neq (a\gamma)c = d\in Z(FS)$. Then
$$
(as^{-1})c = (a\gamma t^{-1})c = (a\gamma c)t^{-1} = dt^{-1}\neq 0,
$$
where $dt^{-1}\in Z(Q_{\text{cl}}(FS))$. Consequently, $Q_{\text{cl}}(FS)$ is a centrally essential ring.

Conversely, let $0\neq s\in FS$. By assumption, there exist elements $t, r\in Z(Q_{\text{cl}}(FS))$ such that $0\neq st = r$. We note that
$Z(Q_{\text{cl}}(FS))\subseteq Q_{\text{cl}}(Z(FS))$; cf. \cite[Theorem 4.4.5]{Pas77}. Indeed, let $\rho\in Z(Q_{\text{cl}}(FS))$, 
$\rho = \alpha\beta^{-1}$, where $\alpha, \beta\in FS$ and $\beta$ is a non-zero-divisor . Then $\alpha\beta = \beta\alpha$ and $\alpha\beta^{-1} = \beta^{-1}\alpha$. By Lemma 3.3.21, there exists a non-zero-divisor $\gamma\in FS$ such that $\beta\gamma\in Z(FS)$. By denoting 
$\epsilon = \beta\gamma$ and $\eta = \alpha\gamma$, we obtain:
$$
\eta\epsilon^{-1} = \alpha\gamma\gamma^{-1}\beta^{-1} = \alpha\beta^{-1} = \rho.
$$
In addition, $\epsilon, \eta\in Z(Q_{\text{cl}}(FS))$. By considering the above, we have $t = cd^{-1}$, $r = mn^{-1}$ for some $c, d, m, n\in Z(FS)$. Then
$$
s(cn) = (sc)n = (mn^{-1}d)n = md\in Z(FS),
$$
and $md\neq 0$, since $d$ is a non-zero-divisor in $FS$.~$\rhd$

\textbf{3.3.24. Open questions.}

\textbf{a.}\label{q-5.3} Is it true the assertion of Remark 3.3.12(a) provided the subgroup $\Delta^{+}(G)$ of the group $G$ is infinite?

\textbf{b.}\label{q-5.4} Is it true that every centrally essential ring has the right classical ring of fractions?

\textbf{c.}\label{q-5.5} Is it true that a centrally essential ring with right classical ring of fractions also has the left classical ring of fractions? 

\subsection{Construction of One Centrally Essential Ring}\label{subsection3.4}

The main results of this subsection are proved in \cite{MT19d}

\textbf{3.4.1. Rings with polynomial identity.}\\ 
Let $X$ be a countable set and let $F=\mathbb{Z}\langle X\rangle$ be the free ring with free generator set $X$. A \textsf{classical identity}\label{cliden} (in the sense of Rowen) is an identity with integral coefficients, i.e. element of the free ring $F$ contained in the kernel of any homomorphism from $F$ into the ring $R$. 
A classical identity is called a \textsf{polynomial identity}\label{poliden} if it is multilinear and has $1$ as one its coefficients; a ring with polynomial identity is called a \textsf{PI ring}\label{pirin}\footnote{See \cite[Definitions 1.1.12, 1.1.17]{Row80}}. 

\textbf{3.4.2. Rings which are algebraic or integral over their centers.}\\
Let $R$ be a ring with center $C=Z(A)$. An element $r\in R$ is said to be \textsf{algebraic}\label{algcenele} (resp., \textsf{integral})\label{intcenele} over the center if for some $n\in \mathbb{N}$, there exist $c_0,\ldots c_n\in C$ such that $c_n$ is a non-zero-divisor in $R$ (resp., invertible element in $R$) and
$$\label{algeq}
c_nr^n+c_{n-1}r^{n-1}+\ldots+c_{1}r+c_0=0.\eqno (3.4.2(1))
$$ 
We denote by $n_{1}(r)$ (resp., $n_{2}(r)$) the least integer $n$ which satisfies this condition. A ring $R$ is said to be \textsf{algebraic}\label{algcenRing} (resp., \textsf{integral})\label{celcenRing} over its center if any element $r\in R$ is algebraic (resp., integral) over its center. We set $m_{1}(R)=\max\{n_{1}(r)\,|\,r\in R\}$ and $m_{2}(R)=\max\{n_{2}(r)\,|\,r\in R\}$; it is possible that $m_{1}(R)=\infty$, $m_{2}(R)=\infty$. 

Finite rings and finite-dimensional algebras over fields are examples of rings $R$ such that $m_1(R)=m_2(R)<\infty$. 

\textbf{3.4.3. Example of a ring which is algebraic over \\ center of the ring and is not integral over the center.}\\
Let
$R=\left\lbrace \begin{pmatrix}
a&b\\0&z
\end{pmatrix}\colon a,b\in \mathbb{Q},\, z\in\mathbb{Z}
\right\rbrace.$
It is clear that the center of the ring $R$ is of the form $\mathbb{Z} E$, where $E$ is the identity matrix.

We note that the ring $R$ is not integral over its center. Indeed, if $r=\begin{pmatrix}
\frac 1 2&0\\0&0
\end{pmatrix}$, then the relation 
$$r^n+c_{n-1}r^{n-1}+\ldots+c_0E=0,$$ where $n\in\mathbb{N}$ and $c_0,\ldots c_{n-1}\in \mathbb{Z}$, implies relations
$c_0=0$ and 
$$
-\dfrac 1 {2^n}=
\dfrac{c_{n-1}} {2^{n-1}}+\ldots+\dfrac{c_{1}} 2;\;\text{this is impossible.}
$$

On the other hand, if $ r=\begin{pmatrix}
a&b\\0&z
\end{pmatrix}\in R$, then $na\in \mathbb{Z}$ for some $n\in \mathbb{N}$, $\text{trace}\,(nr)=na+nz$ and $\det(nr)=na\cdot nz$ are integers. Therefore, $nr$ is a root of the polynomial $x^2-\text{trace}\,(nr)x+\det(nr)\in\mathbb{Z}[x]$ by the Hamilton--Cayley theorem. 

We note that classes of centrally essential rings, PI rings, and rings, which are algebraic or integral over its center, properly contain all commutative rings.

The main result of this subsection is Theorem 3.4.5.
For the proof of Theorem 1.3, we need the following familiar result.

\textbf{3.4.4. Theorem; \cite[Theorem 5.3.9(ii)]{Pas77}.}\\
Let $F$ be a field of characteristic $p>0$. If the group algebra $FG$ satisfies a polynomial identity of degree $d$, then there exists a subgroup $H$ in $G$ such that $[G:H]\cdot|H'|<g(d)$, where $g(d)$ is some function integer $d$.

\textbf{3.4.5. Theorem.}\\ 
For any prime integer $p$ and every field $F$ of characteristic $p$, there exists a centrally essential $F$-algebra which is not a PI ring and is not algebraic over its center.

$\lhd$ We fix a prime integer $p$ and the field $F$ of characteristic $p$. We denote by $Z(G)$ the center of the group~$G$.

For any positive integer $n$, we construct the group $G=G(n)$, see below. 
Let $A=\langle a\rangle$, $B=\langle b\rangle$ and $C=\langle c\rangle$ be three cyclic groups such that 
$|A|=|B|=|C|=p^n$. We consider the automorphism $\alpha\in \text{Aut}\,(B\times C)$ defined on generators by relations $\alpha(b)=bc$ and $\alpha(c)=c$. It is clear that $\alpha^n$ the identity automorphism; 
whence we have we have such homomorphism $\varphi\colon A\rightarrow\text{Aut}\,(B\times C)$ that $\varphi(a)=\alpha$. This 
homomorphism corresponds to the semidirect product $G=(B\times C)\ltimes A$ which can be considered as the group generated by elements $a,b,c$ which satisfy relations
$a^{p^n}=a^{p^n}=c^{p^n}=1$, $bc=cb$, $ac=ca$ and $aba^{-1}=bc$. 
It follows from these relations that $c\in Z(G)$. 
It is directly verified that for any integers $x,y,z, x',y',z'$, we have
$$\label{comm}
[b^yc^za^x,b^{y'}c^{z'}a^{x'}]=
b^ya^xb^{y'}a^{x'}a^{-x}b^{-y}a^{-x'}b^{-y'}
$$
$$
=b^y(a^xb^{y'}a^{-x})(a^{x'}b^{-y}a^{-x'})b^{-y'}= \eqno (3.5.4(1))
$$
$$
=b^y(b^{y'}c^{xy'})(b^{-y}c^{-yx'})b^{-y'}=c^{xy'-yx'}.
$$
Therefore, $Z(G)=G'=\langle c\rangle$ and $G$ is a group of nilpotence class 2.

Now let $H$ be any subgroup of the group $G$. We prove that
$$
\label{ineq}
[G:H]\cdot|H'|\geq p^n. \eqno (3.5.4(2))
$$
We note that $[G:HZ(G)]\leq [G:H]$ and $(HZ(G))'=H'$; consequently, it is sufficient to prove the inequality $(3.5.4(2))$ in the case, where $H\supseteq Z(G)$. We set $\bar G=G/Z(G)$ and denote by $\bar a, \bar b, \bar H$ the images $a,b,H$ under the canonical homomorphism $G$ onto the group $\bar G$. We also set $\bar B=\langle \bar b\rangle$. We have $[G:H]=[\bar G:\bar H]$. It follows from the standard isomorphism $(\bar H\bar B)/\bar B\cong \bar H/(\bar H\cap \bar B)$ that $\bar H/(\bar H\cap\bar B)$ is a cyclic group which is isomorphic to some subgroup of the group $\langle \bar a\rangle$. The group $\bar H\cap\bar B$ is cyclic, as well. Consequently, the group $\bar H$ is generated by two elements
of the form $\bar b^{p^m}$ and $\bar a^{p^k}\bar b^l$ for some non-negative integers $k,l,m$. Therefore,
$$
[\bar G:\bar H]=[\bar G:\bar H\bar B][\bar H\bar B:\bar H]=
[\langle \bar a\rangle:
\langle \bar a^{p^k}\rangle] [\langle \bar b\rangle: \langle \bar b^{p^m} \rangle] =
$$
$$p^kp^m=p^{m+k}.
$$
If $m+k\geq n$, then the inequality $(3.5.4(2))$ holds. If $m+k<n$, then it follows from $(3.5.4(1))$ and the property that elements $a^{p^k}b^l$ and $b^{p^m}$ are contained in the subgroup $H$, that $[a^{p^k}b^l,b^{p^m}]=c^{p^{m+k}}\in H'$. Therefore, $|H'|\geq |\langle c^{p^{m+k}}\rangle|=p^{n-m-k}$ and we have
$$
[G:H]\cdot|H'|\geq p^{m+k}\cdot p^{n-m-k}=p^n,
$$
i.e., $(3.5.4(2))$ holds in this case.

Now it is sufficient to take the direct product of the group algebras $FG(n)$, $n\in \mathbb{N}$, as the ring $R$. We note that the direct product of any set of rings is centrally essential if and only if every factor is a centrally essential ring. Therefore, the ring $R$ is centrally essential by Theorem 3.2.9(b). However, if the algebra $R$ satisfies some polynomial identity of degree $d$, then for any $n\in\mathbb{N}$, the inequality $p^n<g(d)$ follows from $(3.5.4(2))$ and Theorem 3.4.4; this is impossible.

Now we prove that the constructed ring is not algebraic over its center.

It is well known (e.g., see \cite[Proposition 1.1.47]{Row80} or \cite[Lemma 5.2.6]{Zhe82}) that $R$ satisfies polynomial identity of degree $\displaystyle{d(m)=\dfrac{m(m+1)}{2}+m}$ provided $m_1(R)=m<\infty$.

We note that for any $m\in \mathbb{N}$, there exists an integer $n_m$ such that $p^{n_m}>g(d(m))$; in addition, we can take integers $n_1,n_2,\ldots$ such that these integers form an ascending sequence. By the definition of $d(m)$, there exists an element $r'_m\in FG(n_m)$ which does not satisfy any relation of the form $(3.4.2(1))$ of degree $m$. Now we consider the element
$r=\prod_{n=1}^\infty r_n\in \prod_{n=1}^\infty FG(n)$,
where $r_n\in FG(n)$, $r_n=r'_m$ provided $n=n_m$ for some $m\in\mathbb{N}$ and $r_n=0$ otherwise. It is clear that if $r$ satisfies to some relation
of the form $(3.4.2(1))$ of degree $m$, then every element $r_n$ satisfies relation of the same degree; this is impossible by the choice of the element 
$r'_m$.~$\rhd$

\subsection{Centrally Essential Rings $R$ with Non-Commutative $R/J(R)$}\label{subsection3.5}

\textbf{3.5.1. Proposition.}\label{prod}\\
Let $\{R_\alpha\}_{\alpha\in A}$ be an arbitrary set of rings, $R=\prod_{\alpha\in A}R_\alpha$, and let $f(x_1,\ldots,x_n)$ belong to the free ring $\mathbb{Z}\langle X\rangle$ with countable set of free generators. If for any $m\in\mathbb{N}$, there exists infinitely many subscripts $\alpha\in A$ such that the ring $R_\alpha$ does not satisfy the identity $f(x_1,\ldots,x_n)^m$, then the ring $R/K(R)$ does not satisfy the identity $f(x_1,\ldots,x_n)$.

$\lhd$ By assumption, for any $m\in\mathbb{N}$, there exists a subscript $\alpha=\alpha_m\in A$ such that $f(r^{(m)}_1,\ldots,r^{(m)}_n)^m\neq 0$ for some 
$r^{(m)}_1,\ldots,r^{(m)}_n\in R_\alpha$ and all subscripts $\alpha_m, m\in\mathbb{N}$ can be chosen pairwise distinct.

For any $i=1,\ldots,n$, we set 
$s_i=\prod_{\alpha\in A}s^{(i)}_\alpha$, where
$$\label{elements}
s^{(i)}_\alpha=\left\{\begin{array}{l}
r^{(m)}_i\mbox{ for }\alpha=\alpha_m\mbox{ for some }m\in\mathbb{N},\\
0,\mbox{ otherwise,.}\ 
\end{array}\right.\eqno (3.5.1)
$$
Then $f(s_1,\ldots,s_n)=\prod_{\alpha\in A}f(s_\alpha^{(1)},\ldots,s_\alpha^{(n)})$. It follows from $(3.5.1)$ that for any $m\in \mathbb{N}$, there exists a subscript $\alpha\in A$ such that
$f(s_\alpha^{(1)},\ldots,s_\alpha^{(n)})^m\neq 0$; therefore, $f(s_1,\ldots,s_n)^m\neq 0$ and, consequently,
$f(s_1,\ldots,s_n)\notin K(R)$.~$\rhd$

\textbf{3.5.2. Proposition.}\label{polynomials}\\
Under the conditions of Proposition 3.5.1, the ring $R[t]/J(R[t])$ does not satisfy the identity $f(x_1,\ldots,x_n)$.

$\lhd$ For an arbitrary ring $R$, it follows from the Amitsur theorem \cite{Ami56} that $J(R[t])=I[t]$ for some nil-ideal $I$ of the ring $R$. Therefore, we have $J(R[t])\subseteq K(R)[t]$. However, for our purposes, it is sufficient to use the following elementary remark: if $r\in J(R[t])\cap R$, then the element $1-rt$ is invertible. We write $(1-rt)(a_0+a_1t+a_2t^2+\ldots+a_mt^m)=1$ and compare coefficients of degrees $t$. We obtain that $a_0=1$, $a_i=r^i$ for $i=1,\ldots,m$ and $r^{m+1}=0$, i.e., $r$ is a nilpotent element. Therefore, $R\cap J(R[t])\subseteq K(R)$. However, the ring $R/K(R)$ does not satisfy the identity $f(x_1,\ldots,x_n)$, Therefore, the ring $R[t]/J(R[t])$ does not satisfy this identity.~$\rhd$

\textbf{3.5.5. Theorem.}\label{main} 
There exists a centrally essential ring $R$ such that the ring $R/J(R)$ is not a PI ring.\\
Consequently, the ring $R/N(R)$ is not a PI ring, as well; in particular, the rings $R/J(R)$ and $R/N(R)$ are not commutative. 

$\lhd$ We use the sequence of rings from Theorem 3.4.5.

We note that a direct product of any set of rings is centrally essential if and only if every direct factors of the product is a centrally essential ring. Therefore, the ring $R=\prod_{n\in\mathbb{N}}FG(n)$ is centrally essential by the first assertion of Theorem 3.4.5. 

For any polynomial identity $f(x_1,\ldots,x_n)$ of degree $d$ and every integer $m\in \mathbb{N}$, there exists an infinite set of integers $k\in \mathbb{N}$ such that the identity $f(x_1,\ldots,x_n)^m$ does not hold in the ring $FG(k)$. If the ring $FG(k)$ satisfies the identity $f(x_1,\ldots,x_n)^m$, then it also satisfies polynomial identity of degree $dm$ obtained by linearization of this identity. By the second assertion of Theorem 3.4.5, this is impossible for infinite sets of integers $k$. It follows from Proposition 3.5.2 that $R[t]/J(R[t])$ is not a PI ring.

It remains to note that the polynomial ring in one variable over a centrally essential ring is centrally essential by Remark 3.2.10(c).~$\rhd$

\subsection[Local Subalgebras of Triangular Algebras]{Local Sublgebras\\ of Triangular Algebras}\label{subsection3.6}

This subsection is based on \cite{LT22b}. In this subsection, we consider not necessarily unital rings and study local centrally essential subalgebras of the algebra $T_n(\mathbb{F})$ of all upper of triangular matrices, where $\mathbb{F}$ is a field of characteristic $\ne 2$. Such subalgebras are of interest, since, for $\mathbb{F} = \mathbb{Q}$, they are quasi-endomor\-phism algebras of strongly indecomposable torsion-free Abelian groups of finite rank $n$. Quasi-endomorphism algebras of all such groups are local matrix subalgebras of the algebra $M_n(\mathbb{Q})$ of all matrices of order $n$ over a field $\mathbb{Q}$; e.g., see \cite[Chapter I, \S 5]{KMT03}.

We note that the algebra $\mathbb{Q}E$ is the quasi-endomorphism algebra of a strongly indecomposable Abelian torsion-free group of prime rank $p$ if and only if $\mathbb{Q}E$ is isomorphic to to a local subalgebra of the algebra $T_p(\mathbb{\mathbb{Q}})$. Indeed,
$\mathbb{Q}E/J(\mathbb{Q}E)\cong \mathbb{Q}$ in this case; see \cite[Theorem 4.4.12]{Fat07}, where $J(\mathbb{Q}E)$ is the Jacobson radical; it is nilpotent, since $\mathbb{Q}E$ Artinian. It follows from the Weddenburn-Malcev theorem\footnote{For example, see \cite [Theorem 6.2.1]{Dro94}.} that
$\mathbb{Q}E\cong \mathbb{Q}E_p\bigoplus J(\mathbb{Q}E)$, where $E_p$ is the identity matrix.
It is known that every nilpotent subalgebra of the matrix algebra $M_n(\mathbb{F})$ over an arbitrary field $\mathbb{F}$ is transformed by conjugation into nil-triangular subalgebra; see \cite[Chapter 2, Theorem 6]{SupT68}. Since diagonal matrices of a local matrix algebra have equal elements on the main diagonal, they are transformed to itself under conjugation. Consequently, quasi-endomorphism algebras of such Abelian groups are realized as matrix subalgebras if and only if these subalgebras are conjugated with some local subalgebra of the algebra $T_p(\mathbb{\mathbb{Q}})$. Necessary information on Abelian groups is contained in \cite{Fuc15} and \cite{KMT03}.

Let $\mathbb{F}$ be a field and let $\mathcal{A}$ be a finite-dimensional $\mathbb{F}$-algebra. For the nil-algebra $\mathcal{A}$, the maximal nilpotence index $\nu(\mathcal{A})$ of its elements is called the \textsf{nil-index}.\label{nilInd} If $\mathcal{A}^k = (0)$ and $\mathcal{A}^{k-1}\neq (0)$, then $k$ is the \textsf{nilpotence index}\label{nilpInd} of the algebra $\mathcal{A}$ and this algebra is called an algebra of \textsf{nilpotence index} $k$.

Next, $\mathcal{A}$ denotes a local subalgebra of the algebra $T_n(\mathbb{F})$ and $N_n(\mathbb{F})$ denotes the subalgebra of nilpotent matrices in $\mathcal{A}$ (i.e. algebra of properly upper triangular matrices). We note that any matrix $A\in \mathcal{A}$ is of the form
$$
A = \begin{pmatrix}
\lambda& a_{12} &\ldots & a_{1n}\\
0& \lambda &\ldots & a_{2n}\\
\vdots& \vdots &\ddots & \vdots\\
0& 0 &\ldots & \lambda
\end{pmatrix}.
$$
We denote by $E_{ij}$ the matrix unit, i.e. the matrix with $1$ on the position $(i, j)$ and zeros on remaining positions; $E_k$ denotes the identity $k\times k$ matrix. For a subset $S$ of a vector space, the linear hull of $S$ is denoted by $\langle S\rangle$.

\textbf{3.6.1. Proposition.}\\
Let $\mathcal{A}$ be a local subalgebra of the algebra $T_n(\mathbb{F})$ with Jacobson radical $J(\mathcal{A})$. 
The algebra $\mathcal{A}$ is centrally essential if and only if $J(\mathcal{A})$ is a centrally essential algebra.

$\lhd$ Let we have a matrix $A\in J(\mathcal{A})$ with $A\notin Z(J(\mathcal{A}))$. Since $\mathcal{A}$ is a centrally essential algebra, there exists a matrix $B\in Z(\mathcal{A})$ such that $0\neq AB = C\in Z(\mathcal{A})$. Since $J(\mathcal{A})$ is an ideal, we have $C\in Z(J(\mathcal{A}))$. If $B\notin J(\mathcal{A})$, then $A = CB^{-1}\in Z(J(\mathcal{A}))$; this contradicts to the choice of the matrix $A$.

Conversely, let $\mathcal{A} = \mathbb{F}E_n\bigoplus J(\mathcal{A})$.
Since $\mathbb{F}E_n\subset Z(\mathcal{A})$, we have 
$$
Z(J(\mathcal{A}))\subset Z(\mathcal{A}). \eqno(*)
$$ 
If $0\neq A\in \mathcal{A}$ and $A\in Z(\mathcal{A})$, then $0\neq AE_n\in Z(\mathcal{A})$. Let $A\notin Z(\mathcal{A})$ and $A\in J(\mathcal{A})$. Then there exists a matrix $B\in Z(J(\mathcal{A}))$ such that $0\neq AB = C\in Z(J(\mathcal{A}))$.
It follows from relation $(*)$ that $B\in Z(\mathcal{A})$ and $C\in Z(\mathcal{A})$.

Let $A\notin J(\mathcal{A})$. Then $A = A' + A''$, where $0\neq A'\in \mathbb{F}E_n$, $A''\in J(\mathcal{A})$. If $A'' = 0$, then 
$A\in Z(\mathcal{A})$. Otherwise, $0\neq A''B\in Z(J(\mathcal{A}))$ for some $B\in Z(J(\mathcal{A}))$. Then
$$
AB = A'B + A''B = BA' + BA'' = BA.
$$
Since $A'B, A''B\in Z(J(\mathcal{A}))$, we have $AB\in Z(J(\mathcal{A}))\subset Z(\mathcal{A})$. We note that $AB\neq 0$, since the matrix $A$ is invertible.~$\rhd$

It follows from Proposition 3.6.1 that the problem of constructing local centrally essential subalgebras of the algebra $T_n(\mathbb{F})$ is equivalent to the problem of constructing centrally essential subalgebras of the algebra $N_n(\mathbb{F})$.

Let $\mathcal{A}$ be a subalgebra of the algebra $N_n(\mathbb{F})$ of nilpotence index $n$. We assume that $\nu(\mathcal{A}) = n$. There exists a matrix $A\in \mathcal{A}$ such that $A^{n-1}\neq 0$. We transform $A$ to the Jordan normal form,
$$
A = E_{12} + E_{23} +\ldots +E_{(n-1)n},
$$
and pass to the corresponding conjugated subalgebra $\mathcal{A}_c$. We denote by $\text{Cen}(A)$ the centralizer of the matrix $A$ in $\mathcal{A}_c$. 
Since the minimal polynomial of the matrix $A$ coincides with its characteristic polynomial, we have
$\text{Cen}(A) = \mathbb{F}[A]$, where $\mathbb{F}[A]$ is the ring of all matrices which can be presented in the form $f(A)$, $f(x)\in \mathbb{F}[x]$; 
see \cite[Chapter 1, Theorem 5]{SupT68}. For $B\in \text{Cen}(A)$, we have
$$
B = f(A) = \alpha_0E_n + \alpha_1A +\ldots +\alpha_{n-1}A^{n-1}. 
$$
In addition, $\alpha_0 = 0$, since the matrix $B$ is nilpotent. 

\textbf{3.6.2. Remark.} \\
If $Z(\mathcal{A}_c) = \text{Cen}(A)$, then the algebra $\mathcal{A}_c$ is commutative.

$\lhd$ Indeed, if $A'\notin \text{Cen}(A)$, then $AA'\neq A'A$. 
However, $A\in \text{Cen}(A) = Z(\mathcal{A}_c)$. This is a contradiction.~$\rhd$

\textbf{3.6.3. Remark.}\\ 
Let $\mathcal{A}_c$ be a centrally essential algebra and\\ $Z(\mathcal{A}_c) = \langle A^{n-1}\rangle$. Then the algebra $\mathcal{A}_c$ is commutative.

$\lhd$ Indeed, if $\mathcal{A}_c$ is not commutative and the matrix $A'$ is not in $Z(\mathcal{A}_c)$, then we have $BA' = 0$ for any matrix $B\in Z(\mathcal{A}_c)$.~$\rhd$

\textbf{3.6.4. Remark.}\\ 
If $\mathbb{F}$ is a field of characteristic $\ne 2$, then
every centrally essential subalgebra of the algebra $N_3(\mathbb{F})$ is commutative.

$\lhd$ Every matrix $A\in N_3(\mathbb{F})$ is of the form
$$
A = \begin{pmatrix}
0& a & b\\
0& 0 &c\\
0& 0 & 0
\end{pmatrix}.
$$
Let $\mathcal{A}$ be a non-commutative centrally essential subalgebra of the algebra $N_3(\mathbb{F})$ of nilpotence index 3.
Then $\nu(\mathcal{A}) = 3$. Let the matrix $A\in \mathcal{A}$ have the nilpotence index 3. We transform $A$ to the Jordan normal form:
$A = E_{12} + E_{23}$. Now if $B\in \text{Cen}(A)$, then $B = \alpha_1A +\alpha_2A^2$. We note that $Z(\mathcal{A}_c)\subseteq \text{Cen}(A)$; in addition, $\nu(Z(\mathcal{A}_c)) = 3$ by Remark 3.6.3. However, $Z(\mathcal{A}_c) = \text{Cen}(A)$ and algebra $\mathcal{A}_c$ is commutative by Remark 3.6.2. This is a contradiction.~$\rhd$

\textbf{3.6.5. Remark.}\\ 
It follows from Remark 3.6.4 that all centrally essential endomorphism rings of strongly indecomposable Abelian torsion-free groups of rank 3 are commutative.

\textbf{3.6.6. Proposition.}\\ 
Any centrally essential subalgebra $\mathcal{A}$ of the algebra $N_4(\mathbb{F})$ is commutative.

$\lhd$ If the algebra $\mathcal{A}$ is of nilpotence index index $2$, then it is commutative. Let the nilpotence index of the algebra $\mathcal{A}$ be equal to $4$. There exists a matrix $A\in \mathcal{A}$ such that $A^3\neq 0$. Indeed, the algebra $\mathcal{A}$ contains three matrices $S = (s_{ij})$, $T = (t_{ij})$, $P = (p_{ij})$ with $s_{12}\neq 0$, $t_{23}\neq 0$, $p_{34}\neq 0$. Otherwise, the nilpotence index $\mathcal{A}$ is less than $4$. 
As the the required matrix, we can take a matrix $A = (a_{ij})$ such that $a_{i(i+1)}\neq 0$, $i = 1, 2, 3$.
We transform $A$ to the Jordan normal form,
$$
A = E_{12} + E_{23} + E_{34},
$$
and pass to the corresponding conjugated subalgebra $\mathcal{A}_c$. For the matrix $B\in \text{Cen}(A)$, we have
$$
B = \alpha_1A + \alpha_2A^2 + \alpha_3A^3,
$$
where $\alpha_1, \alpha_2, \alpha_3\in \mathbb{F}$.
It follows from Remarks 3.6.2 and 3.6.3 that $Z(\mathcal{A}_c)\neq \text{Cen}(A)$ and $Z(\mathcal{A}_c)\neq <A^3>$ if algebra $\mathcal{A}_c$ is not commutative. Then any matrix $C\in Z(\mathcal{A}_c)$ is of the form
$$
C = \begin{pmatrix}
0 & 0 & c_{13} & c_{14}\\
0 & 0 & 0 & c_{13}\\
0 & 0 & 0 &0\\
0 & 0 & 0 & 0\\
\end{pmatrix}.
$$
Since $\mathcal{A}_c$ is a centrally essential algebra, we have that for a non-zero matrix $D\notin Z(\mathcal{A}_c)$, there exists a matrix $C\in Z(\mathcal{A}_c)$ such that $0\neq DC\in Z(\mathcal{A}_c)$. Since the matrix $D$ is nilpotent, we have $\textbf{tr} D = 0$. In addition, $\mathcal{A}_c$ is local; therefore, all elements on the main diagonal of the matrix $D$ are equal to zero. In this case, it is directly verified that 
$$
D = \begin{pmatrix}
0 & d_{12} & d_{13} & d_{14}\\
0 & 0 & d_{23} & d_{24}\\
0 & 0 & 0 & d_{12}\\
0 & 0 & 0 & 0\\
\end{pmatrix}.
$$ 
If $d_{12} = 0$ and $D\notin Z(\mathcal{A}_c)$, then $DC = 0$ for any matrix $C\in Z(\mathcal{A}_c)$; this is a contradiction. 
Let $d_{12}\neq 0$ and $DF \neq FD$ for some matrix $F = (f_{ij})\in \mathcal{A}_c$. We take an element $\lambda\in \mathbb{F}$ such that $f_{12} = \lambda d_{12}$. We set $G = \lambda D - F$, $G = (g_{ij})$. Then 
$$
FG = F(\lambda D - F) = \lambda FD - F^2,
$$
$$
GF = (\lambda D - F)F = \lambda DF - F^2.
$$
Therefore, $G\notin Z(\mathcal{A}_c)$ and $g_{12} = 0$. It follows from obtained contradiction that the algebra $\mathcal{A}_c$ is commutative.

Let the nilpotence index of the algebra $\mathcal{A}$ be equal to 3. Then $\nu(\mathcal{A}) = 3$, i.e., $\mathcal{A}$ contains a matrix $A$ such that $A^2\neq 0$. Indeed, we assume the contrary,
$A^2 = 0$ for all $A\in \mathcal{A}$. If $A\notin Z(\mathcal{A})$, then $0\neq AB\in Z(\mathcal{A})$ for some matrix $B\in Z(\mathcal{A})$.
Then
$$
(A + B)^2 = A^2 + 2AB + B^2 = 2AB = 0.
$$
Therefore, $AB = 0$. This is a contradiction.

We transform the matrix $A$ to a Jordan normal form,
$$
A = E_{12} + E_{23}.
$$
In the corresponding conjugated subalgebra $\mathcal{A}_c$, the centralizer $\text{Cen}(A)$ consists of matrices $B$ of the form
$$
B = \begin{pmatrix}
0 & b_{12} & b_{13} & b_{14}\\
0 & 0 & b_{12} & 0\\
0 & 0 & 0 & 0\\
0 & 0 & b_{43} & 0\\
\end{pmatrix}; \eqno (3.6.6(1))
$$ 
see \cite[Chapter 3, \S 1]{SupT68}. In addition, if $C\in Z(\text{Cen}(A))$, then we have
$$
C = \begin{pmatrix}
0 & c_{12} & c_{13} & 0\\
0 & 0 & c_{12} & 0\\
0 & 0 & 0 & 0\\
0 & 0 & 0 & 0\\
\end{pmatrix}. \eqno (3.6.6(2))
$$ 
Let $Z(\mathcal{A}_c)$ have the nilpotence index $3$. Then we can take a matrix in $Z(\mathcal{A}_c)$ as the matrix $A$; see 
\cite[Chapter 1, Proposition 5, Corollary]{SupT68}.
In this case, all matrices in $\mathcal{A}_c$ are contained in $\text{Cen}(A)$. Then $\mathcal{A}_c$ consists of matrices of the form $(3.6.6(1))$
and matrices in $Z(\mathcal{A}_c)$ of the form $(3.6.6(2))$.
If $B = (b_{ij})\notin Z(\mathcal{A}_c)$ and $b_{12} = 0$, then $BC = 0$ for all $C\in Z(\mathcal{A}_c)$. Then $\mathcal{A}_c$ is not a centrally essential algebra.
Let $b_{12}\neq 0$ and $BD\neq DB$ for some matrix $D = (d_{ij})\in\mathcal{A}_c$. Let $d_{12} = \lambda b_{12}$ and $F = \lambda B - D$, 
$F = (f_{ij})$. Then $f_{12} = 0$ and $F\notin Z(\mathcal{A}_c)$. This is a contradiction.

Let $Z(\mathcal{A}_c)$ be of nilpotence index 2. Then for $C\in Z(\mathcal{A}_c)$, we obtain
$$
C = \begin{pmatrix}
0 & 0 & c_{13} & 0\\
0 & 0 & 0 & 0\\
0 & 0 & 0 & 0\\
0 & 0 & 0 & 0\\
\end{pmatrix}. 
$$
It follows from relation $AC = CA$ for $A\in \mathcal{A}_c$ that
$$
A = \begin{pmatrix}
0 & a_{12} & a_{13} & a_{14}\\
0 & 0 & a_{23} & a_{24}\\
0 & 0 & 0 & 0\\
0 & a_{42} & a_{43} & 0\\
\end{pmatrix}. 
$$
However, $AC = 0$ for any matrix $C\in Z(\mathcal{A}_c)$.
Consequently, if $\mathcal{A}_c$ is a centrally essential algebra, then $\mathcal{A}_c$ is commutative.~$\rhd$

\textbf{3.6.7. Remark.}\\ 
In Theorem 2.2.3, it is proved that the Grassmann algebra $\Lambda(V)$ over a field $\mathbb{F}$ of characteristic $\ne 2$ is a centrally essential algebra if and only if the dimension of the space $V$ is odd. By considering the regular matrix representation of the algebra $\Lambda(V)$,
we obtain that for an odd positive integer $n > 1$, there exists a non-commutative centrally essential subalgebra of the algebra $N_{2^n}(\mathbb{F})$; also see Example 3.7.10 below. Therefore, the minimal order of matrices of non-commutative centrally essential Grassmann of the algebra is equal to $8$.

We recall that for a right ideal $I$ of the ring $R$, any right ideal $J$ in $R$ which is maximal with respect to the property $I\cap J = 0$, is said to be \textsf{$\cap$-complement} for $I$.

\textbf{Example 3.6.8.} There exists a non-commutative centrally essential algebra of $7\times 7$ matrices which has a closed right ideal which is not an ideal.

We consider the subalgebra $\mathcal{A}$ in $N_{7}(\mathbb{F})$ consisting of matrices $A$ of the form 
$$
A = 
\left(\begin{matrix}
0 & a & b & c & d & e & f\\
0 & 0 & 0 & b & 0 & 0 & d\\
0 & 0 & 0 & 0 & 0 & 0 & e\\
0 & 0 & 0 & 0 & 0 & 0 & 0\\
0 & 0 & 0 & 0 & 0 & 0 & a\\
0 & 0 & 0 & 0 & 0 & 0 & b\\
0 & 0 & 0 & 0 & 0 & 0 & 0\\
\end{matrix}\right).
$$
Let for $A'\in \mathcal{A}$, we have $a' = a + 1$ and the remaining components coincide with the corresponding components of the matrix $A$. Then $AA'\neq A'A$ if $a\neq 0$ and $b\neq 0$. Therefore, the algebra $\mathcal{A}$ is not commutative. It is easy to see that $Z(\mathcal{A})$ consists of matrices $C$ of the form
$$
C = 
\left(\begin{matrix}
0 & 0 & 0 & c & d & e & f\\
0 & 0 & 0 & 0 & 0 & 0 & d\\
0 & 0 & 0 & 0 & 0 & 0 & e\\
0 & 0 & 0 & 0 & 0 & 0 & 0\\
0 & 0 & 0 & 0 & 0 & 0 & 0\\
0 & 0 & 0 & 0 & 0 & 0 & 0\\
0 & 0 & 0 & 0 & 0 & 0 & 0\\
\end{matrix}\right).
$$
If $0\neq A\notin Z(\mathcal{A})$, then $0\neq AC\in Z(\mathcal{A})$ for some matrix $C\in Z(\mathcal{A})$. Consequently, $\mathcal{A}$ is a centrally essential algebra. 

We consider the right ideal $I$ in $\mathcal{A}$ consisting of matrices of the form 
$$
B = 
\left(\begin{matrix}
0 & 0 & b & 0 & 0 & 0 & f\\
0 & 0 & 0 & b & 0 & 0 & 0\\
0 & 0 & 0 & 0 & 0 & 0 & 0\\
0 & 0 & 0 & 0 & 0 & 0 & 0\\
0 & 0 & 0 & 0 & 0 & 0 & 0 \\
0 & 0 & 0 & 0 & 0 & 0 & b \\
0 & 0 & 0 & 0 & 0 & 0 & 0\\
\end{matrix}\right).
$$
It is directly verified that $I$ is not an ideal of $\mathcal{A}$. In addition, $I$ is a closed right ideal. Indeed, the ideal of $\mathcal{A}$ which has only the element $c$ as a non-zero component, is a $\cap$-complement for $I$. 

In the same time, the closed left ideal $J$ in $\mathcal{A}$ whose elements are matrices 
$$
D = 
\left(\begin{matrix}
0 & a & 0 & 0 & 0 & 0 & f\\
0 & 0 & 0 & 0 & 0 & 0 & 0\\
0 & 0 & 0 & 0 & 0 & 0 & 0\\
0 & 0 & 0 & 0 & 0 & 0 & 0\\
0 & 0 & 0 & 0 & 0 & 0 & a \\
0 & 0 & 0 & 0 & 0 & 0 & 0 \\
0 & 0 & 0 & 0 & 0 & 0 & 0\\
\end{matrix}\right),
$$
is not an ideal. The ideal, which has only element $c$ as the non-zero component, is a $\cap$-complement also and for $J$. 

\textbf{3.6.9. Theorem.}\\
For any field $\mathbb{F}$ of characteristic $\ne 2$ and an arbitrary positive integer $n\ge 7$, there exists a local non-commutative centrally essential subalgebra of the algebra $T_n(\mathbb{F})$ upper of triangular $n\times n$ matrices. 

$\lhd$ In $N_n(\mathbb{F})$, we consider the subalgebra $\mathcal{A}$ of matrices $A$ of the form 
$$
A = 
\left(\begin{matrix}
0 & a_{12} & a_{13} & a_{14} & a_{15} & \ldots & a_{1n-2} & a_{1n-1} & a_{1n}\\
0 & 0 & 0 & a_{13} & 0 & \ldots & 0 & 0 & a_{1\;n-2}\\
0 & 0 & 0 & 0 & 0 & \ldots & 0 & 0 & a_{1\;n-1}\\
0 & 0 & 0 & 0 & 0 & \ldots & 0 & 0 & 0\\
\ldots & \ldots & \ldots & \ldots & \ldots & \ldots & \ldots & \ldots &\ldots\\
0 & 0 & 0 & 0 & 0 & \ldots & 0 & 0 & 0\\
0 & 0 & 0 & 0 & 0 & \ldots & 0 & 0 & a_{12}\\
0 & 0 & 0 & 0 & 0 & \ldots & 0 & 0 & a_{13}\\
0 & 0 & 0 & 0 & 0 & \ldots & 0 & 0 & 0\\
\end{matrix}\right).
$$
We note that the algebra $\mathcal{A}$ is not commutative; also see Example 3.6.8. If $B\in Z(\mathcal{A})$, then
$$
B = 
\left(\begin{matrix}
0 & 0 & 0 & b_{14} & b_{15} & \ldots & b_{1n-2} & b_{1n-1} & b_{1n}\\
0 & 0 & 0 & 0 & 0 & \ldots & 0 & 0 & b_{1n-2}\\
0 & 0 & 0 & 0 & 0 & \ldots & 0 & 0 & b_{1n-1}\\
0 & 0 & 0 & 0 & 0 & \ldots & 0 & 0 & 0\\
\ldots & \ldots & \ldots & \ldots \ldots & \ldots & \ldots & \ldots &\ldots & \ldots\\
0 & 0 & 0 & 0 & 0 & \ldots & 0 & 0 & 0\\
\end{matrix}\right).
$$
For $A = (a_{ij})\notin Z(\mathcal{A})$, we have $a_{12}\neq 0$, $a_{13}\neq 0$. 
Let $B = (b_{ij})\in Z(\mathcal{A})$ and $b_{1n-2} = a_{12}$, $b_{1n-1} = a_{13}$. Then $0\neq AB\in Z(\mathcal{A})$.
Indeed, let $AB = C = (c_{ij})$, $BA = D = (d_{ij})$. Then $c_{ij} = d_{ij} = 0$ for all $i\neq 1$, $j\neq n$.
In addition, $c_{1n} = d_{1n} = a_{12}^2 + a_{13}^2$.
Therefore, $\mathcal{A}$ is a centrally essential algebra.~$\rhd$

\subsection{Endomorphism Rings of Abelian groups}\label{subsection3.7}

In this subsection, we study Abelian groups $A$ with centrally essential rings endomorphism ring $\text{End}\,A$. The subsection is based on \cite{LT20}. 

We denote by $\text{End}\,A$ the endomorphism ring of an Abelian group $A$. If $A = \bigoplus_{p\in P} A_p$ is a decomposition of a torsion Abelian group $A$ into a direct sum of $p$-components, then $\text{supp} \,A = \{p\in P \mid A_p\neq 0\}$. We use the following notation: $\mathbb {Z}_{p^k}$ (resp., $Z_{p^k}$) is the residue ring (resp., the additive group modulo $p^k$);~ $\mathbb Q$ (resp., $Q$) is the ring (resp., the additive group) of rational numbers; $Z_{p^{\infty}}$ is a quasi-cyclic Abelian $p$-group; $\hat{\mathbb Z}_p$ is the ring of $p$-adic integers. 

An Abelian group $A$ is said to be \textsf{divisible}\label{divGr} if $nA = A$ for any positive integer $n$. An Abelian group is said to be \textsf{reduced}\label{redGr} if it does not contain non-zero divisible subgroups and \textsf{non-reduced}\label{nonRed} otherwise. 

A subgroup $B$ of an Abelian group $A$ is said to be \textsf{pure}\label{purSub} if the equation $nx = b\in B$ which is solvable in the group $A$, also is solvable in $B$.

\textbf{3.7.1. Remark.}\\ 
Let $A$ be an Abelian group which is either a torsion group or a non-reduced group and let the endomorphism ring $\text{End}\,A$ is centrally essential. We will prove below that the ring $\text{End}\,A$ is commutative. Therefore, only reduced torsion-free groups and reduced mixed groups are of interest under the study Abelian groups with non-commutative centrally essential endomorphism rings.

Theorem 3.7.13(c) below contains an example of an Abelian torsion-free group of finite rank with centrally essential non-commutative ring endomorphism ring. In Example 3.7.15 below, we give other examples of non-commutative centrally essential endomorphism rings of some Abelian torsion-free groups of infinite rank.

Let $A$ be an Abelian torsion-free group. A \textsf{pseudo-socle}\label{pseSoc} $\text{PSoc}\,A$ of the group $A$ is the pure subgroup of the group $A$ generated by all its minimal pure fully invariant subgroups.

\textbf{3.7.2. Lemma.} 
Let $A$ be a module and let $A = \bigoplus_{i\in I}A_i$ be the direct decomposition of the module $A$. The endomorphism ring $\text{End} A$ is centrally essential if and only if for every $i\in I$, all rings $\text{End} A_i$ are centrally essential and all $A_i$ are fully invariant submodules in $A$. 

$\lhd$ Let $\text{End}\,A = E$ be a centrally essential ring. If $A_i$ is not a fully invariant submodule for some $i\in I$, then there exists a subscript $j\in I$, $j\neq i$, such that $\text{Hom}\,(A_i, A_j) = e_jEe_i\neq 0$, where $e_i$ and $e_j$ are projections from the module $A$ onto the submodules $A_i$ and $A_j$, respectively. In addition,
$$
e_i\cdot e_jEe_i = 0\neq e_jEe_i = e_jEe_i\cdot e_i,
$$
i.e. the idempotent $e_i$ is not central; this contradicts to 1.1.4.

If $A_i$ is a fully invariant submodule in $A$, $i\in I$, then $\text{End}\,A\cong \text{End}\,A_i\times \text{End}\,\overline{A}_i$, where $\overline{A}_i$ is a complement direct summand for $A_i$. It is clear that if the ring $\text{End}\,A_i$ is not centrally essential, then and $\text{End}\,A$ is not centrally essential.

We assume that for any $i\in I$, the ring $\text{End} A_i$ is centrally essential and $A_i$ is a fully invariant submodule in $A$. Then $\text{End}\,A\cong \prod_{i\in I}\text{End}\,A_i$ and all rings $\text{End}\,A_i$ are centrally essential. It is clear that $\text{End}\,A$ is centrally essential, as well.~$\rhd$

\textbf{3.7.3. Lemma.}\\ 
For a divisible Abelian group $A$, the endomorphism ring of $A$ is centrally essential if and only if $A\cong Q$ or $A\cong Z_{p^{\infty}}$.

$\lhd$ Let $A = F(A)\bigoplus T(A)$, where $0\neq F(A)$ is the torsion-free part and $0\neq T(A)$ is the torsion part of the group $A$. Then the subgroup $F(A)$ in $A$ is not fully invariant (see \cite[Theorem 7.2.3]{Fuc15}) and, by Lemma 3.7.2, the ring $\text{End}\,A$ is not centrally essential. Hypothetically, $F(A)$ or $T(A)$ is the direct sum of groups ${\mathbb{Z}}_p^{\infty}$ or $\mathbb{Q}$. It is clear that if the number of summands exceeds $1$, then $\text{End}\,A$ has a non-central idempotent; this is a contradiction.~$\rhd$
 
Let $A = \bigoplus_{p\in P} A_p$ be the decomposition of the torsion Abelian group $A$ into a direct sum of primary components of $A$. It follows from of Lemma 3.7.2 that $\text{End}\,A$ is a centrally essential ring if and only if every ring $\text{End}\,A_p$ is centrally essential.

\textbf{3.7.4. Lemma.}
The endomorphism ring of a primary Abelian group $A_p$ is centrally essential if and only if $A_p\cong Z_{p^k}$ or $A_p\cong Z_{p^{\infty}}$.

$\lhd$ If the group $A_p$ is not indecomposable, then it has a cocyclic direct summand; see \cite[Corollary 5.2.3]{Fuc15}. By considering \cite[Theorem 7.1.7, Example 7.1.3 and Theorem 7.2.3]{Fuc15}, we see that this summand and complement summands of it are fully invariant in $A$. Consequently, $A_p\cong Z_{p^k}$ or $A_p\cong Z_{p^{\infty}}$. The converse is obvious, since the rings $\mathbb{Z}_{p^k}$ and $\hat{\mathbb{Z}}_p$ are commutative.~$\rhd$

\textbf{3.7.5. Theorem.} 
Let $A = D(A)\bigoplus R(A)$ be a non-reduced Abelian group, where $0\neq D(A)$ and $0\neq R(A)$ is the divisible part and the reduced part of the group $A$, respectively. The endomorphism ring of the group $A$ is centrally essential if and only if 
$A = D(A)\bigoplus R(A)$, where $R(A) = \bigoplus_{p\in P'}Z_{p^k}$ and $D(A)\cong Q$ or $D(A)\cong \bigoplus_{p\in P''} Z_{p^{\infty}}$;
$P', P''$ are subsets of distinct prime numbers with $P' \cap P'' = \emptyset$.

$\lhd$ Let $\text{End}\,A$ be a centrally essential ring. We verify that $D(A)$ and $R(A)$ are fully invariant of the subgroup in $A$. Indeed, it is well known that $\text{Hom}\,(D(A), R(A)) = 0$. If $R(A)$ is a torsion-free group, then $\text{Hom}\,(R(A), D(A)) \neq 0$ (see \cite[Theorem 7.2.3]{Fuc15}); this contradicts to Lemma 3.7.2. In addition, it is clear that $\text{Hom}\,(R(A), D(A)) \neq 0$ if $R(A)$, $D(A)$ are torsion groups and $\text{supp} \,R(A) \cap \text{supp} \,D(A) \neq \emptyset$. It follows from Lemma 3.7.4 that $R(A)$ is the direct sum of its cyclic $p$-components and it follows from Lemma 3.7.3 that $D(A)\cong Q$ or $D(A)\cong \bigoplus_{p\in P} Z_{p^{\infty}}$.

The converse assertion follows from Lemmas 3.7.2--3.7.4.~$\rhd$

\textbf{3.7.6. Corollary.}
The endomorphism ring of a non-reduced Abelian group is centrally essential if and only if the ring is commutative. In other words, only reduced Abelian groups can have non-commutative centrally essential endomorphism rings.

$\lhd$ Indeed, it follows from Theorem 3.7.5 that the centrally essential endomorphism ring of an arbitrary non-reduced Abelian group is a direct product of rings which can be only of the ring $\mathbb{Z}_{p^k}$, $\mathbb{Q}$ and $\hat{\mathbb{Z}}_p$.~$\rhd$

\textbf{3.7.7. Quasi-decompositions and strong indecomposability.}\\ 
Let $A$ and $B$ be two Abelian torsion-free groups. One says that $A$ is \textsf{quasi-contained}\label{quaCon} in $B$ if $nA\subseteq B$ for some positive integer $n$. If $A$ is quasi- contained in $B$ and $B$ is quasi-contained in $A$ (i.e. if $nA\subseteq B$ and $mB\subseteq A$ for some
$n,m\in \mathbb{N}$), then one says that $A$ is \textsf{quasi-equal}\label{quaEq} to $B$ (one writes $A\doteq B$). A quasi-relation $A\doteq\bigoplus_{i\in I}A_i$ is called a \textsf{quasi-decomposition}\label{quaDec} (or a \textsf{quasi-direct decomposition}) of the Abelian group $A$; these subgroups $A_i$ are called \textsf{quasi-summands} of the group $A$. If the group $A$ does not have non-trivial quasi-decompositions, then $A$ is said to be \textsf{strongly indecomposable}.\label{strInd} The ring $\mathbb{Q}\otimes \text{End}\,A$ is called the \textsf{quasi-endomorphism ring}\label{quasEndRing} of the group $A$; it is denoted by $\mathbb{Q}\text{End}\,A$; see details in \cite[Chapter I, \S 5]{KMT03}. Elements of the ring $\mathbb{Q}\otimes \text{End}\,A$ are called \textsf{quasi-endomorphisms}\label{quasEnd} of the group $A$. We note that
$$
\mathbb{Q}\text{End}\,A = \{\alpha\in \text{End}_{\mathbb{Q}}(Q\otimes A) \mid (\exists n\in \mathbb{N}) (n\alpha\in \text{End}\,A)\}.
$$
It is well known (\cite [Proposition 5.2]{KMT03}) that the correspondence
$$
A \doteq e_1A\bigoplus \ldots \bigoplus e_kA \to \mathbb{Q}\text{End}\,A = 
$$
$$
=\mathbb{Q}\text{End}\,A e_1\bigoplus \ldots \bigoplus\mathbb{Q}\text{End}\,A e_k
$$
between finite quasi-decompositions of the torsion-free group $A$ and finite decompositions of the ring $\mathbb{Q}\text{End}\,A$ in direct sum of left ideals, where
$\{e_i \mid i = 1,\ldots,k\}$ is a complete orthogonal system idempotents of the ring $\mathbb{Q}\text{End}\,A$, is one-to-one.

\textbf{3.7.8. Proposition.}\\
For an Abelian torsion-free group $A$, the endomorphism ring $E$ of $A$ is centrally essential if and only if the quasi-endomorphism ring $\mathbb{Q}E$ of $A$ is centrally essential.

$\lhd$ Let $0\neq \tilde{a}\in \mathbb{Q}E$. For some $n\in \mathbb{N}$, we have $n\tilde{a} = a\in E$ and there exist $x,y\in Z(E)$ with $ax = y\neq 0$. In this case, $\tilde{a}\tilde{x} = \tilde{y}$, where $\tilde{x} = x$, $\tilde{y} = \cfrac{1}{n}\cdot y\in Z(\mathbb{Q}E)$, i.e., $\mathbb{Q}E$ is a centrally essential ring.

Conversely, for every $0\neq a\in E$, there exist non-zero $\tilde{x},\tilde{y}\in Z(\mathbb{Q}E)$ with $a\tilde{x} = \tilde{y}$. In addition, there exist $n, m\in \mathbb{N}$ such that $n\tilde{x}\in Z(E)$ and $m\tilde{y}\in Z(E)$. Then $ax = y$, where $x = mn\tilde{x}, y = mn\tilde{y}\in Z(E)$.~$\rhd$

Let $A\doteq \bigoplus_{i = 1}^n A_i = A'$ be a decomposition of the Abelian torsion-free group $A$ of finite rank into a quasi-direct sum of strongly indecomposable groups (e.g., see \cite[Theorem 5.5]{KMT03}). By using Lemma 3.7.2 and Proposition 3.7.8, we obtain that the ring $\text{End}\,A$ is centrally essential if and only if all subgroups $A_i$ are fully invariant in $A'$, and every ring $\text{End}\,A_i$ is centrally essential. 
Therefore, the problem of describing Abelian torsion-free groups of finite rank with centrally essential endomorphism rings is reduced to a similar problem for strongly indecomposable groups.

\textbf{3.7.9. Proposition.}\\
Let $A$ be a strongly indecomposable Abelian group and $A = \text{PSoc}\,A$. If the ring $\text{End}\,A$ is centrally essential, then the ring $\text{End}\,A$ is commutative.

$\lhd$ If $A = \text{PSoc}\,A$, then $\text{End}\,A$ is a semiprime ring (e.g., see \cite [Theorem 5.11]{KMT03}). By Theorem 1.2.2, the ring $\text{End}\,A$ is commutative.~$\rhd$

\textbf{3.7.10. Example.}
We take centrally essential endomorphism rings of strongly indecomposable Abelian torsion-free groups of rank 2 and 3.

If $A$ is a strongly indecomposable group of rank $2$, then the ring $\text{End}\,A$ is commutative (e.g., see \cite[Theorem 4.4.2]{Fat07}). Consequently, $\text{End}\,A$ is a centrally essential ring. Let $A$ be a strongly indecomposable group of rank $3$. Then the algebra $\mathbb{Q}\text{End}\,A$ is isomorphic to one of the following $\mathbb{Q}$-algebras (\cite[Theorem 2]{Che98}):
$$
K\cong\left\{\left. \left(\begin{matrix}
x & 0 & z\\0 & x & 0\\0 & 0 & x 
\end{matrix}\right)\right| x, z\in \mathbb{Q}\right\},\, 
R\cong\left\{\left.\left(\begin{matrix}
x & y & z\\0 & x & 0\\0 & 0 & x 
\end{matrix}\right)\right| x, y, z\in \mathbb{Q}\right\},
$$
$$
S\cong\left\{\left. \left(\begin{matrix}
x & y & z\\0 & x & ky\\0 & 0 & x 
\end{matrix}\right)\right| x, y, z\in \mathbb{Q},\,0\neq k\in \mathbb{Q},\,k = \text{const}\right\},
$$
$$
T\cong\left\{\left. \left(\begin{matrix}
x & y & z\\0 & x & t\\0 & 0 & x 
\end{matrix}\right)\right| x, y, z, t\in \mathbb{Q}\right\}.
$$
The rings $K$, $R$, $S$ are commutative; consequently, they are centrally essential. The ring $T$ is not commutative (in addition, $\text{PSoc}\,A$ has the rank $1$). We have 
$$
J(T) = \left\{\left. \left(\begin{matrix}
0 & y & z\\0 & 0 & t\\0 & 0 & 0 
\end{matrix}\right) \right| y, z, t\in \mathbb{Q}\right\},
$$
$$
Z(T) = \left\{\left. \left(\begin{matrix}
x & 0 & z\\0 & x & 0\\0 & 0 & x 
\end{matrix}\right) \right| x, z \in \mathbb{Q}\right\},
$$ 
$$
M = \left\{\left. \left(\begin{matrix}
0 & 0 & 0\\0 & 0 & t\\0 & 0 & 0 
\end{matrix}\right) \right| t\in \mathbb{Q}\right\},
$$ 
where $M$ is a minimal right ideal of the ring $T$. We note that the ring $T/J(T)$ is commutative but $Z(T) \cap M = 0$. 
It follows from Remark 3.7.11(a) that the ring $T$ is not centrally essential.
As a result, we obtain that endomorphism rings of strongly indecomposable groups of rank 2 or 3 are centrally essential if and only if they are commutative.

\textbf{3.7.11. Remarks.}\\
Let $R$ be a local Artinian ring with center $Z(R) = C$ and $R$ is not a division ring. 

\textbf{a.} If $R$ is centrally essential, then $R/J(R)$ is commutative and $C\cap M\neq 0$ for every minimal ideal $M$.

\textbf{b.} If $R/J(R)$ is commutative, $\text{Soc}\, (R_C) = \text{Soc}\, (R_R)$ and $C\cap M\neq 0$ for every minimal ideal $M$, then $R$ is centrally essential.

$\lhd$. \textbf{a.} Let $R$ be centrally essential. By Theorem 1.3.2, the ring $R/J(R)$ is commutative. Since $R$ is Artinian, the ideal $J(R)$ is nilpotent; let $k$ be the nilpotence index $J(R)$. 
We note that if $M$ is a minimal ideal of $R$, then $MJ(R) = 0$.

We assume that $C\cap M = 0$ for some minimal ideal $M$. By assumption, for $0\neq a\in M$, there exist $x, y\in Z(R)$ such that $ax = y\neq 0$. Since $x\notin J(R)$ (otherwise $ax = 0$), the element $x$ is invertible in $R$ and $a = x^{-1}y\in C$; this is a contradiction. 

\textbf{b.} Let $C\cap M\neq 0$ for every minimal ideal $M$. We verify that $M\subseteq C$. Let $C\cap M = K$. By assumption, $R/J(R)$ is commutative. Therefore, $rs - sr\in J(R)$ for all $r,s\in R$. Then $k(rs - sr) = 0$ for every $k\in K$. 
In addition, since $k\in C$, we have $(kr)s = ksr = s(kr)$ and $kr\in C$. Similarly, $rk\in C$. In addition, $kr\in M$ and $rk\in M$. Therefore, 
$K$ is an ideal.
Since the ideal $M$ is minimal, we have that $K = M$ or $K = 0$. However, $K\neq 0$, whence we have $K = M$ and $M\subset C$. 
Therefore, $\text{Soc}\,R_C=\text{Soc}\,R_R\subseteq C$. By Theorem 4.4.1(b), the ring $R$ is centrally essential.~$\rhd$

\textbf{3.7.12. Example.}
Let $V$ be a vector $\mathbb{Q}$-space with basis $e_1, e_2, e_3$ and let $\Lambda(V)$ be the Grassmann algebra of the space $V$, i.e., $\Lambda(V)$ is an algebra with operation $\wedge$, generators $e_1, e_2, e_3$ and defining relations
$$
e_i \wedge e_j + e_j \wedge e_i = 0 \quad \mbox {for all} \quad i, j = 1, 2, 3.
$$
Then $\Lambda(V)$ is a $\mathbb{Q}$-algebra of dimension 8 with basis $$
\{1, e_1, e_2, e_3, e_1 \wedge e_2, e_2 \wedge e_3, e_1 \wedge e_3, e_1 \wedge e_2 \wedge e_3\}
$$
and $\Lambda(V)$ is a non-commutative centrally essential ring; see Introduction, Example 2.
We consider the regular representation $\Lambda(V)$. If $x\in \Lambda(V)$ and
$$
x = q_0\cdot 1 + q_1e_1 + q_2e_2 + q_3e_3 + q_4e_1 \wedge e_2 + q_5e_2 \wedge e_3 + q_6e_1 \wedge e_3 + q_7e_1 \wedge e_2 \wedge e_3,
$$
then the matrix $A_x\in Mat_8(\mathbb{Q})$ is of the form
$$
\left(\begin{matrix}
q_0 & q_1 & q_2 & q_3 & q_4 & q_5 & q_6 & q_7\\
0 & q_0 & 0 & 0 & -q_2 & 0 & -q_3 & q_5\\
0 & 0 & q_0 & 0 & q_1 & -q_3 & 0 & -q_6\\
0 & 0 & 0 & q_0 & 0 & q_2 & q_1 & q_4\\
0 & 0 & 0 & 0 & q_0 & 0 & 0 & q_3\\
0 & 0 & 0 & 0 & 0 & q_0 & 0 & q_1\\
0 & 0 & 0 & 0 & 0 & 0 & q_0 & -q_2\\
0 & 0 & 0 & 0 & 0 & 0 & 0 & q_0 \end{matrix}\right).
$$
We denote by $R$ the corresponding subalgebra in $\text{Mat}_8(\mathbb{Q})$. It is clear that the radical $J(R)$ consists of properly upper triangular matrices in $R$ and $A_x\in Z(R)$ if and only if $q_1 = q_2 = q_3 = 0$. In addition, 
$\text{Soc}\,R_R = \{A_x = (a_{ij})\in R \mid a_{ij} = 0, i\neq 1, j\neq 8\}$ and
$\text{Soc}\,R_C = \{A_x = (a_{ij})\in Z(R) \mid a_{ii} = 0\}$.
Since $\text{Soc}\,(R_C)\neq \text{Soc}\,R_R$, the corresponding condition of Remark 3.7.11(b) is not necessary.

\textbf{3.7.13. Theorem.}
Let $A$ be a strongly indecomposable Abelian torsion-free group of finite rank, $\mathbb{Q}\text{End}\,A$ be the quasi-endomorphism algebra, and let $A\neq \text{PSoc}\,A$. 

\textbf{a.} If $\mathbb{Q}\text{End}\,A$ is a centrally essential ring, then the ring $\mathbb{Q}\text{End}\,A/J(\mathbb{Q}\text{End}\,A)$ is commutative and $Z(\mathbb{Q}\text{End}\,A) \cap M\neq 0$ for every minimal right ideal $M$ of the ring $\mathbb{Q}\text{End}\,A$.

\textbf{b.} Let the ring $\mathbb{Q}\text{End}\,A/J(\mathbb{Q}\text{End}\,A)$ be commutative, 
$$
\text{Soc}\,(\mathbb{Q}\text{End}\,A_{\mathbb{Q}\text{End}\,A}) = \text{Soc}\,(\mathbb{Q}\text{End}\,A_{Z(\mathbb{Q}\text{End}\,A)})$$
and $Z(\mathbb{Q}\text{End}\,A) \cap M\neq 0$ for every minimal right ideal $M$ of the ring $\mathbb{Q}\text{End}\,A$. Then the ring $\mathbb{Q}\text{End}\,A$ is centrally essential.

\textbf{c.} Let $n > 1$ be an odd integer. There exists a strongly indecomposable Abelian torsion-free group $A(n)$ of rank $2^n$ such that its endomorphism ring is a non-commutative centrally essential ring.
 
$\lhd$ \textbf{a, b.} It is known that the ring $\mathbb{Q}\text{End}\,A$ is a local Artinian ring (e.g., see \cite[Corollary 5.3]{KMT03}). It remains to use 3.7.11.

\textbf{c.} By Theorem 2.2.3, the Grassmann algebra $\Lambda(V)$ over a field $F$ of characteristic $0$ or $p\neq 2$ is a centrally essential ring if and only if the dimension of the space $V$ is odd. We set $F = \mathbb{Q}$. It is known (e.g., see \cite{PieV83}) that every $\mathbb{Q}$-algebra of dimension $n$ can be realized as the quasi-endomorphism ring of an Abelian torsion-free group of rank $n$. Therefore, we consider Example 3.7.11 and Proposition 3.7.8 and obtain the required property.~$\rhd$
 
Under the conditions of Theorem 3.7.13, if the rank of the group $A$ is square-free, then the ring $\mathbb{Q}\text{End}\,A/J(\mathbb{Q}\text{End}\,A)$ is commutative \cite[Lemma 4.2.1]{Fat07}. By considering Proposition 3.7.8, we obtain Corollary 3.7.14.

\textbf{3.7.14. Corollary.}
Let $A$ be a strongly indecomposable Abelian torsion-free group of finite rank, $A\neq \text{PSoc}\,A$, and let the rank of the group $A$ be square-free. 

\textbf{a.} If the endomorphism ring $\text{End}\,A$ of the group $A$ is centrally essential, then $Z(\mathbb{Q}\text{End}\,A) \cap M\neq 0$ for every minimal right ideal $M$ of the ring $\mathbb{Q}\text{End}\,A$.

\textbf{b.} If for every minimal right ideal $M$ of the ring $\mathbb{Q}\text{End}\,A$, we have 
$\text{Soc}\,(\mathbb{Q}\text{End}\,A_{\mathbb{Q}\text{End}\,A}) = \text{Soc}\,(\mathbb{Q}\text{End}\,A_{Z(\mathbb{Q}\text{End}\,A)})$ and $Z(\mathbb{Q}\text{End}\,A) \cap M\neq 0$, then the ring $\text{End}\,A$ is centrally essential.

\textbf{3.7.15. Example.}\\ 
Let $R = \mathbb{Z}[x,y]$ be the polynomial ring in two variables $x$ and $y$. We use the construction described in \cite[Proposition 7]{Jel16}. We consider the ring
$$
T(R) = \left\{\left. \left(\begin{matrix}
f & d_1(f) & g\\0 & f & d_2(f)\\0 & 0 & f 
\end{matrix}\right) \right| f, g\in \mathbb{Z}[x,y]\right\},
$$
where $d_1$, $d_2$ are two derivations of the ring $\mathbb{Z}[x,y]$, $d_1 = \cfrac{\partial}{\partial x}$, $d_2 = \cfrac{\partial}{\partial y}$.
The ring $T(R)$ is not commutative and $J(R) = e_{13}R\subseteq Z(T(R))$, where $e_{13}$ is the matrix unit; \cite[Corollary 8]{Jel16}. If $0\neq a\in T(R)\setminus Z(T(R))$, then $0\neq ae_{13}\in Z(T(R))$. Therefore, $T(R)$ is centrally essential.
Since $T(R)$ is a countable ring with reduced additive torsion-free group, it follows from the familiar Corner theorem (e.g., see \cite[Theorem 29.2]{KMT03}) that the ring $T(R)$ contains $\mathfrak{M}$ of Abelian groups $A_i$ such that $\text{End}\,A_i\cong T(R)$ 
and $\text{Hom}\,(A_i, A_j) = 0$ for all $i\neq j$, where $\mathfrak{M}$ is an arbitrary predetermined cardinal number; see \cite{DugG82}, \cite{CorG85}. We note that the endomorphism ring of a direct sum of such groups is a non-commutative centrally essential ring, as well.

\section[Distributive and Uniserial Rings]{Distributive and\\ Uniserial Rings}\label{chapter4}

\subsection{Uniserial Artinian Rings}\label{section4.1}

This subsection is based on \cite{MT20}.

The following fact is well known and is directly verified.

\textbf{4.1.1. Uniserial Artinian rings.}\label{simpleseries}\\
Let $R$ be a ring with Jacobson radical $J=J(R)$. If $J$ is nilpotent of nilpotence index $n$, then the following conditions are equivalent.

\textbf{a)} $J^{k-1}/J^{k}$ is a simple left $R$-module for all $k=1,\ldots,n$ (we assume that $J^0=R$).

\textbf{b)} $R$ is a left uniserial, left Artinian ring.

\textbf{c)} $R$ is a local ring and $J$ is a principal left ideal of $R$.

\textbf{4.1.2. Lemma.}\label{sigma}
Let $R$ be a left uniserial, left Artinian ring, $J=J(R)=R\pi$, $D=\overline{R}$ be a division ring, and let $\sigma\colon D\rightarrow D$ be the mapping defined by the relation
$$
\sigma(\overline{r})=\overline{a}, \mbox{ where } a\pi = \pi r.
\eqno (*)
$$ 
Then $\sigma$ is a homomor\-phism from the division ring $D$ into itself.

$\lhd$ First, we note that the mapping $\sigma$ is well defined. Indeed, the existence of the element $a$ from $(*)$ follows from the property that $R\pi$
is a two-sided ideal. If $r,r'\in R$, $\pi r=a\pi$, $\pi r'=a'\pi$, and $\overline{r}=\overline{r'}$, then $(a-a')\pi=\pi(r-r')\in J^2$. However, $J/J^2$ is an one-dimensional linear space over the division ring $\overline{R}$ generated by the element $\pi+J^2$; therefore, $\overline{a-a'}=0$ and $\overline{a}=\overline{a'}$.

Second, for any two elements $\text{r}_1,\text{r}_2\in R$, we have
$$
\begin{array}{l}
\sigma(\text{r}_1+\text{r}_2)(\pi+J^2)=
\pi(\text{r}_1+\text{r}_2)+J^2=\\
\qquad\qquad=\pi \text{r}_1+\pi
\text{r}_2+J^2=(\sigma(\text{r}_1)+\sigma(\text{r}_2))(\pi+J^2),\\
\sigma(\text{r}_1\text{r}_2)(\pi+J^2)=\pi
\text{r}_1\text{r}_2+J^2=\\
\qquad\qquad=\sigma(\text{r}_1)(\pi \text{r}_2+J^2)
=\sigma(\text{r}_1)\sigma(\text{r}_2)(\pi+J^2).
\end{array}
$$
Therefore, $\sigma$ is a ring homomor\-phism.~$\rhd$

\textbf{4.1.3. Remark.} Without using special links, we often use the property that for any centrally essential local ring $R$, the division ring $\overline{R}$ is a field by Theorem 1.3.2. In particular, this is the case if $R$ is a left or right uniserial centrally essential ring.

\textbf{4.1.4. Proposition.}\label{sigma-id}\\
Let $R$ be a left uniserial, left Artinian, centrally essential ring, $C=Z(R)$, and let $J=J(R)$. Then the homomor\-phism $\sigma$ from Lemma 4.1.2 is the identity automor\-phism.

$\lhd$ Let $c$ be an element of the center of the ring $R$ with $c\pi\in C\setminus \{0\}$. Let $r$ be any
element of the ring $R$. We have $c=a\pi^k+b$, where $b\in J^{k+1}$ and $\overline{a}\neq 0$. It follows from the relation $rc=cr$ that $\overline{a}(\overline{r}-\sigma^k(\overline{r}))=0$; then we have $c\pi= a\pi^{k+1}+b\pi$ and it follows from the relation $r(c\pi)=(c\pi)r$ that $\overline{a}(\overline{r}-\sigma^{k+1}(\overline{r}))=0$. Since $\text{Ker}\,(\sigma)=0$, it follows from the relation $\sigma^k(\overline{r})=\sigma^{k+1}(\overline{r})$ that $\overline{r}=\sigma(\overline{r})$.~$\rhd$

\textbf{4.1.5. Corollary.}\\ 
If a left uniserial, left Artinian ring $R$ is centrally essential, then $R$ is a right uniserial, right Artinian ring.

$\lhd$ If $J(R)=R\pi$ for some element $\pi\in J(R)$, it follows from Proposition 4.1.4 that $J(R)=\pi R+J(R)^2$ and $J(R)=\pi R$ by the Nakayama lemma.
It remains to note that the right-sided analogue of condition \textbf{c)} of 4.1.1 holds.~$\rhd$

We recall that for a ring $R$, we denote by $\ell_R(S)=\{r\in R\,|\, rS=0\}$ the \textsf{left annihilator} of an subset $S$ of $R$. The right annihilator $\text{r}_R(S)$ is similarly defined.

\textbf{4.1.6. Proposition.}\label{suff}\\
Let $R$ be a left Artinian, left uniserial ring with center $C=Z(A)$ and Jacobson radical $J$ and let $n$ be the nilpotence index of the ideal $J$. If $\displaystyle{J^{[n/2]}\subseteq C}$, then the ring $R$ is centrally essential.

$\lhd$
First, we note that $\ell_R(J^k)=\text{r}_R(J^k)=J^{n-k}$ for any $k=0,1,\ldots,n$. In particular, $$\ell_R(J^{[n/2]})\subseteq J^{[n/2]}.$$

Let $0\neq r\in R$. If $r\in J^{[\dfrac n 2]}$, then $r\in C$. If $r\not\in J^{[\dfrac n 2]}$, then it remarked above that $r \not \in \ell_R(J^{[\dfrac n 2]})$. Therefore, $rJ^{[n/2]}\neq 0$ and
$rJ(R)^{[\dfrac n 2]}\subseteq J(R)^{[n/2]}\subseteq C$.
In the both cases, we have $rC\cap C\neq 0$.~$\rhd$

\textbf{4.1.7. Open question.}\\ Is it true that the assertion, which is converse to Proposition 4.1.6, holds?

Now we prove that there exists a non-commutative uniserial centrally essential ring. For this purpose, we use the construction which is similar to the one described in \cite{Jel16}.

For a field $F$, we recall that a \textsf{derivation}\label{deriv} of $F$ is an arbitrary endomor\-phism of the additive the group $(F,+)$ which satisfies the relation $\delta(ab)=a\delta(b)+\delta(a)b$ for any two elements $f,b\in F$. General properties of derivations are given, e.g., in \cite[\S II.17]{ZarS75}. Any field
has the \textsf{trivial} derivation $F\rightarrow 0$. An example of a non-trivial derivation is the ordinary derivation on the field of rational functions.

For a ring $R$, we denote by $[a,b]$ the \textsf{commutator}\label{comm} $ab-ba$ of two elements $a,b$ of the ring $R$ and we denote by $[A,b]$ the ideal of $R$ generated by the set $\{[a,b]\,|\,a\in A,b\in B\}$, where $A,B$ are any two subsets of $R$. For any three elements $a,b,c\in R$, we have the following well known properties of commutators: $[a,b]=-[b,a]$, $[ab,c]=a[b,c]+[a,c]b$.

\textbf{4.1.8. Example.}\label{matr}\\
Let $F$ be a field with non-trivial derivation $\delta$. Then there exists a non-commutative Artinian uniserial centrally essential ring $R$ with $R/J(R)\cong F$.

$\lhd$ We consider a mapping $f\colon F\rightarrow M_3(F)$ from the field $F$ into the ring of $4\times 4$ matrices over $F$ defined by the the relation
$$
\forall a\in F,\quad f(a)=\begin{pmatrix}a&0&0&0\\
0&a&0&0\\\delta(a)&0&a&0\\0&0&0&a\end{pmatrix}.
$$
It is directly verified that $f$ is a ring homomor\-phism.
We set $$x=\begin{pmatrix}0&0&0&0\\
1&0&0&0\\0&1&0&0\\0&0&1&0\end{pmatrix}.$$
We consider the subring $R$ of the ring $M_3(F)$ generated by the set $f(F)\cup\{x\}$. It is easy to see that $xf(a)=f(a)x+f(\delta(a))x^3$ for any $a\in A$. Therefore, $Rx=xR$, $(Rx)^4=0$ and $R/Rx\cong F$. It follows from 
4.1.1 that $R$ is a uniserial Artinian ring.

Since $Rx^2\subseteq Z(R)$, it follows from Proposition 4.1.6 that the ring $R$ is centrally essential.

Finally, if $a\in R$ and $\delta(a)\neq 0$, then $[x,f(a)]=f(\delta(a))x^3\neq 0$,
therefore, the ring $R$ is not commutative.~$\rhd$

\textbf{4.1.9. Proposition.}\label{noncomm}\\
Let $R$ be a uniserial Artinian ring with radical $J=R\pi$ such that $\overline{R}$ is a field and $[r,\pi]\in J^2$ for any $r\in R$. If the ring $R$ is not commutative, then the field $F$ has a non-trivial derivation.

$\lhd$ Let $\gamma\colon \overline{R}\rightarrow R$
be an arbitrary mapping such that $\overline{\gamma{\overline{r}}}=\overline{r}$ for any $r\in R$; in other words, $\gamma(\overline{r})$ is a fixed representative of the coset $r+J$. Without loss of generality, we can assume that $\gamma(\overline{0})=0$. We
set $\Gamma=\gamma(\overline{R})$. Then any element $r\in R$ is uniquely represented as the sum
$$
\label{gamma-rep}
r=\sum_{i=0}^{n-1}g_i\pi^i, \eqno (*)
$$
where $g_i\in \Gamma$ for all
$i=0,1,\ldots,n-1$ (we assume that $\pi^0=1$). Indeed, $r-g_0\in J$ for the unique element $g_0=\overline{r}$. Next, if $g_0,\ldots,g_{k-1}$ are already defined with
$s=r-\sum_{i=0}^{k-1}g_i\pi^i\in{J^{k}}$, where $0 < k< n$, then the next coefficient is uniquely defined as $\gamma(\lambda)$ from the relation
$$
s+J^{k+1}=\lambda (\pi^{k}+J^{k+1}),\quad \lambda\in\overline{R}.
$$
For $k=n-1$, we obtain the required representation.

First, we assume that $\pi\not\in Z(R)$. By $(*)$, we have 
$[\Gamma,\pi]\neq 0$. Let $n(R)$ be the nilpotence index of $J(R)$. 
Then $[\Gamma,\pi]=J^k$ for some integer $k$ with $2\leq k < n(R)$. By the induction on the positive integer $m$, we prove that $[\Gamma,\pi^m]\subseteq J^{m-1+k}$ for any $m>0$.
Indeed, this is true for $m=1$ by the choice of $k$. Further for any $g\in\Gamma$, we have
$$[g,\pi^{m+1}]=\pi[g,\pi^m]+[g,\pi]\pi^m\in JJ^{m-1+k}+J^k\pi^m\subseteq J^{m+k}.$$
Thus, for any $g,h\in \Gamma$, we have $[g,[h,\pi]]=[g,x\pi^k]$ for some $x\in R$. Therefore,
$$
\label{three-comm}
[g,[h,\pi]]=[g,x]\pi^k+x[g,\pi^k]\in J\pi^k+J^{2k-1}\subseteq J^{k+1}, 
$$
since $k\geq 2$. It follows that $[hg,\pi]=h[g,\pi]+[h,\pi]g=$
$$\label{comm-prod}
=h[g,\pi]+g[h,\pi]-[g,[h,\pi]]\in h[g,\pi]+g[h,\pi]+J^{k+1}.
\eqno (**)
$$

Since $J^k/J^{k+1}$ is a simple left module, it is an one-dimensional left vector space over the field $\overline{R}$ with basis consisting of one element $v=\pi^k+J^{k+1}$. We define a mapping
$\delta_\pi\colon \overline{R}\rightarrow \overline{R}$ by the rule 
$$[\gamma(\overline{r}),\pi]+{J^{k+1}}=\delta_\pi(\overline{r})v\mbox{ for
any }r\in R.
$$
If $r\in R\pi^m$ for some $m>0$, we have that the coefficients
$g_0,\ldots,g_{m-1}$ of the representation $(*)$ are equal to $0$; therefore, $[r,\pi]=\sum_{i=m}^{n-1}[g_i,\pi]\pi^i\in J^{k+m}$. Therefore, for any $r,s\in R$, we have
$$[\gamma(\overline{r})\gamma(\overline{s}),\pi]+J^{k+1}=[\gamma(\overline{r}\overline{s}),\pi]+J^{k+1}=\delta_\pi(\overline{r}\overline{s})v.$$

On the other hand, it follows from $(**)$ that
$$
[\gamma(\overline{r})\gamma(\overline{s}),\pi]=\gamma(\overline{r})[\gamma(\overline{s}),\pi]+[\gamma(\overline{r}),\pi]\gamma(\overline{s})+{J^{k+1}}=
$$
$$
=(\overline{r}\delta_\pi(\overline{s})+\overline{s}\delta_\pi(\overline{r}))v;
$$
therefore,
$\delta_\pi(\overline{r}\overline{s})=\overline{r}\delta_\pi(\overline{s})+
\overline{s}\delta_\pi(\overline{r})$.
The relation
$\delta_\pi(\overline{r}+\overline{s})= \delta_\pi(\overline{r})+\delta_\pi(\overline{s})$ is similarly verified.
Consequently, $\delta_\pi$ is a non-trivial derivation of the field $\overline{R}$.

Now we assume that $\pi\in Z(R)$. Then $[\Gamma,\Gamma]=J^k$ where $1\leq k<n(R)$. We choose an element $f\in\Gamma$ such that $[\Gamma,f]\not\subseteq
J^{k+1}$ and define a new mapping
$\delta_f\colon\overline{R}\rightarrow \overline{R}$ by the rule
$$
[\gamma(\overline{r}),f]+{J^{k+1}}=\delta_f(\overline{r})v\mbox{ for any
}r \in R,
$$
where $v=\pi^k+J^{k+1}$. We verify that $\delta_f$ is a derivation. First, we note that $[J,\Gamma]\subseteq
J^{k+1}$, since, for $r\in J$, we have $g_0=0$ in the representation $(*)$ and
$$[r,g]=\sum_{i=1}^{n-1}[g_i,g]\pi^i\in\sum_{i=1}^{n-1}J^k\pi^i=J^{k+1}$$
for any $g\in \Gamma$. It follows that for any $r,s\in R$, we have
$$[\gamma(\overline{r})\gamma(\overline{s}),f]+J^{k+1}=[\gamma(\overline{r}\overline{s}),f]+J^{k+1}=\delta_f(\overline{r}\overline{s})v.$$
On the other hand, for any $r,s\in R$, we have
$$\renewcommand{\arraycolsep}{0pt}
\begin{array}{rcl} [\gamma(\overline{r})\gamma(\overline{s}),f]&=&\gamma(\overline{r})[\gamma(\overline{s}),f]+[\gamma(\overline{r}),f]\gamma(\overline{s})\\
&=&\gamma(\overline{r})[\gamma(\overline{s}),f]+\gamma(\overline{s})[\gamma(\overline{r}),f]-[\gamma(\overline{s}), [\gamma(\overline{r}),f]\\&\in& \gamma(\overline{r})[\gamma(\overline{s}),f]+\gamma(\overline{s})[\gamma(\overline{r}),f]+[\Gamma,J]\\ &\subseteq& \gamma(\overline{r})[\gamma(\overline{s}),f]+\gamma(\overline{s})[\gamma(\overline{r})+J^{k+1}\\&=&(\overline{r}\delta_f(\overline{s})+\overline{s}\delta_f(\overline{r}))v.
\end{array}$$

Consequently $\delta_f(\overline{r}\overline{s})=\overline{r}\delta_f(\overline{s})+
\overline{s}\delta_f(\overline{r})$. The relation
$\delta_f(\overline{r}+\overline{s})= \delta_f(\overline{r})+\delta_f(\overline{s})$ is similarly verified. It follows that $\delta_f$ is a non-trivial derivation of the field $\overline{R}$.~$\rhd$

For a field $F$, it is well known (e.g., see \cite[\S II.17]{ZarS75}) that $F$ does not have a non-trivial derivation provided $F$ is a separable algebraic extension of its prime subfield (all finite fields and all
fields of algebraic numbers are such fields) or $F$ is a perfect field (i.e., ${\rm char} F=p>0$ and $F^p=F$).

\textbf{4.1.10. Theorem.}\\
A field $F$ does not have a non-trivial derivation if and only if
any left uniserial, left Artinian, centrally essential ring $R$ with \hbox{$R/J(R)\cong F$} is commutative.

$\lhd$ Theorem 4.1.10 follows from Proposition 4.1.9 and Proposition 4.1.4.~$\rhd$

For a ring $R$ and any element $r$ (resp., a subset $S$) of the
ring $R$, we set $\overline{r}=r+J(R)\in R/J(R)$ (resp., $\overline{S}=\{\overline{s}\,|\,s\in S\}$). In particular, $\overline{R}=R/J(R)$.

\textbf{4.1.11. Theorem.}\\
\textbf{a.} Every finite, left uniserial, centrally essential ring $R$ is commutative.

\textbf{b.} There exists a non-commutative uniserial Artinian centrally essential ring.

$\lhd$ \textbf{a.} For the finite local ring $R$, the division ring $\overline{R}$ is a field by the Wedderburn theorem \cite[Theorem 3.1.1]{Her05}. In addition, this field does not have a non-zero derivation. Then the ring $R$ is commutative by Theorem 4.1.10. 

\textbf{b.} The assertion follows from 4.1.8.~$\rhd$

\subsection{Uniserial Noetherian Rings}\label{section4.2}

This subsection is based on \cite{MT20} and \cite{MT20b}.

\textbf{4.2.1. Remarks.}\\
\textbf{a.} It is directly verified that the ring $A$ is a right (resp., left) uniserial, right (resp., left) Noetherian ring if and only if $R$ is a local principal right (resp., left) ideal ring. 

\textbf{b.} It follows from \textbf{a} and Theorem 1.2.2 that centrally essential uniserial Noetherian semiprime rings coincide with commutative local principal right ideal domains.

\textbf{c.} There exist right uniserial, right Noetherian rings which are neither prime rings no right Artinian rings; e.g., see \cite[Example 9.10(3)]{Tug98}.

For convenience, we give brief proofs of the following two well known assertions.

\textbf{4.2.2. Remarks.}\label{div-fg-mod}\\
\textbf{a.}\label{N=aN}
Let $A$ be a right uniserial ring and $B$ a completely prime ideal of $A$. Then $B=aB$ for every $a\in A\setminus B$.

\textbf{b.} Let $A$ be a commutative domain which has a non-zero finitely generated divisible torsion-free $A$-module $M$. Then $A$ is a field.

$\lhd$ \textbf{a.} Let $a\in A\setminus B$. Since $aA\not\subseteq B$, we have $B\subseteq aA$. Therefore, for every $x\in B$, there exists an element $b\in A$ with $x=ab$. Since $B$ is a completely prime ideal, $b\in B$ and $x\in aB$.

\textbf{b.} Let's assume the contrary. Then $A$ has a non-zero maximal ideal $\mathfrak m$ and $M$ naturally turns into a non-zero finitely generated module over the local ring $R_{\mathfrak m}$ with radical $J={\mathfrak m}R_{\mathfrak m}$. Since the module $M$ is divisible, we have that $MJ\supseteq M{\mathfrak m}=M$ and $M=0$ by the Nakayama lemma. This is a contradiction.~$\rhd$

\textbf{4.2.3. Lemma.}\label{local-pri}
Let $A$ be a local ring and let $J(A)=\pi A$ for some element $\pi\in A$ of nilpotence index $n$ (maybe, $n=\infty$). For any two integers $k,\ell $ and each $a,b\in A$ such that $k,\ell \geq 0$, $k+\ell <n$, $a\in\pi^kA\setminus\pi^{k+1}A$ and $b\in\pi^lA\setminus\pi^{\ell +1}A$, we have $ab\in\pi^{k+\ell }A\setminus\pi^{k+\ell +1}A$.

$\lhd$ It follows from the inclusion $A\pi\subseteq \pi{}A$ that $ab\in \pi^{k+\ell }A$. If $\pi^m\in \pi^{m+1}A$ for some $m\geq 0$, then it is clear that $\pi^m(1-\pi t)=0$ for some $t\in A$ and $\pi^m=0$, since $1-\pi t\in A^*$. We set $a=\pi^{k}r$ and $b=\pi^{\ell }s$ for some $r,s\in A\setminus J(A)$. Then $r,s\in A^*$, since the ring $A$ is local and $r\pi^\ell \in\pi^lA\setminus\pi^{\ell +1}A$. Consequently, $r\pi=\pi{}r'$ for some $r'\in A^*$ and $ab=\pi^{k+\ell }r's$. It remains to remark that $ab\not\in\pi^{k+\ell +1}A$, since $r's\in A^*$ and
$\pi^{k+\ell }\neq 0$.~$\rhd$

\textbf{4.2.4. Lemma.}\label{left}\\
A right uniserial, right Artinian, centrally essential ring is a left uniserial, left Artinian ring.

$\lhd$ Let $A$ be a right uniserial, right Artinian, centrally essential ring, $N=J(A)$, and let $n$ be the nilpotence index of the ideal $N$. If $n=1$, then the ring $A$ is commutative by Lemma 2.2(3,4); it is nothing to prove in this case. Any right uniserial ring is a local ring; therefore, every element of $A\setminus N$ is invertible. Let $n>1$, i.e., $N\neq 0$. Since a right Noetherian (e.g., a right Artinian) right uniserial ring is a principal right ideal ring, $N=\pi A$ for some element $\pi\in N$. There exist two elements $x,y\in Z(A)$ with $\pi x=y\neq 0$.
Let $x\in N^k \setminus N^{k+1}$ for some $k$, $0\leq k<n$. Then $y\in N^{k+1}$, whence $k+1<n$. If $[a,\pi]\not \in N^2$, then it follows from Lemma 4.2.3 that $[a,\pi]x\not\in N^{k+2}$; consequently, $[a,\pi]x\neq 0$. However, $[a,\pi]x=[a,\pi x]=[a,y]=0$. This is a contradiction; therefore, $[a,\pi]\in N^2$ for every $a\in A$. Consequently, $N=A\pi+N^2$, whence $N/A\pi=N(N/a\pi)=\ldots=N^n(N/A\pi)=0$, i.e., $N=A\pi$.
It follows from the left-side analogue of Lemma 4.2.3 that every left ideal of the ring $A$ coincides with one of the ideals $A,N,N^2, \ldots, N^{n-1},\{0\}$, i.e., $A$ is a left uniserial, left Artinian ring.~$\rhd$

\textbf{4.2.5. Theorem.}\label{main}
For a ring $A$, the following conditions are equivalent.

\textbf{a)} $A$ is a right uniserial, right Noetherian, centrally essential ring.

\textbf{b)} $A$ is a left uniserial, left Noetherian, centrally essential ring.

\textbf{c)} $A$ is a commutative local principal ideal domain or a uniserial Artinian ring.

$\lhd$ It is sufficient to prove the equivalence of conditions \textbf{a} and \textbf{c}. 

\textbf{c)}\,$\Rightarrow$\,\textbf{a)}. The implication is directly verified.

\textbf{a)}\,$\Rightarrow$\,\textbf{c)}.
We set $N=\text{Sing}\,A_A$. The ideal $N$ is nilpotent; e.g., see \cite[9.2]{Tug98}. It follows from 1.1.2(c,d) that the ideal $N$ is completely prime and contains all zero-divisors of the ring $R$ and the ring $A/N$ is a commutative domain. 

Therefore, the proposition is true for $N=0$. Now let $N\neq 0$. We denote by $n$ the nilpotence index of the ideal $N$. Then $0\neq N^{n-1}\subseteq{}\ell_A(N)$. It follows from 1.1.2(e) that $N^{n-1}\subseteq Z(A)$. Next, for every $a\in A\setminus N$, we have $N=aN$ by Remark 4.2.2(a), whence $N^{n-1}=aN^{n-1}=N^{n-1}a$. Consequently, $N^{n-1}$ is a divisible right $(A/N)$-module and $N^{n-1}$ is a torsion-free $(A/N)$-module, since all zero-divisors of the ring $A$ are contained in $N$. By Remark 4.2.2(b), the ring $A/N$ is a field and each of the cyclic $(A/N)$-modules $(N^{k-1}/N^k)$ for $k=1,\ldots,n$ is a simple module. Consequently, the right uniserial ring $A$ is right Artinian. By Lemma 4.2.4, $A$ is a uniserial Artinian ring.~$\rhd$

\textbf{4.2.6. Example.}\label{exists}\\
Let $F$ be a field and let $D_1,D_2\colon F\rightarrow F$ be two derivations of the field $F$ with incomparable kernels (for example, we can take the field of rational functions $\mathbb{Q}(x,y)$ in two independent variables as $F$ and set $D_1=\partial/\partial x$, $D_2=\partial/\partial y$).

Then for every positive integer $n\geq 2$, there exists a non-com\-mu\-ta\-tive uniserial, Artinian, centrally essential ring $A$ such that $A/J(A)\cong F$ and the nilpotence index of $J(A)$ is equal to $n$.

$\lhd$ We use a construction which is similar to the one described in \cite{Jel16}. Let $N=2n-1$, $R=M_N(F)$ be the matrix ring of order $N$ over the field $F$, $e_{i,j}$ denote the matrix unit for any $i,j\in\{1,\ldots,N\}$, and let $f\colon F\rightarrow R$ be the mapping defined by the rule
$$
f(\alpha)=\alpha E+D_1(\alpha)e_{1,N-1}+D_2(\alpha)e_{N-1,N}
$$
for every $\alpha\in F$, where $E$ is the identity matrix. Let $A$ be the subring of the ring $R$ generated by the set $f(F)$ and the matrix
$\pi=\sum_{i=1}^{n-1}e_{2i-1,2i+1}$. It is directly verified that 
$\pi^n=0$, $\pi^{n-1}=e_{1,N}$, $f(\alpha)\pi=\pi f(\alpha)=\alpha \pi$ and 
$$[f(\alpha),f(\beta)]=(D_1(\alpha)D_2(\beta)-D_1(\beta)D_2(\alpha))\pi^{n-1}$$ 
for any $\alpha,\beta\in F$. It follows from these relations that $\pi A=A\pi=J(A)$, $J(A)^k=\pi^kA=A\pi^k$ for all $k=1,\ldots,n-1$ and $\pi A\subseteq Z(A)$. It is clear that $A$ is a uniserial Artinian ring. If $a\in A\setminus \{0\}$ and $a\in \pi A$, then $a\in Z(A)$; otherwise, $a\pi^{n-1}\in Z(A)\setminus \{0\}$ and $X^{n-1}\in Z(A)$. Consequently, the ring $A$ is centrally essential.

Finally, if $\alpha\in \text{Ker}\,D_2\setminus \text{Ker}\,D_1$ and $\beta\in \text{Ker}\,D_1\setminus \text{Ker}\,D_2$, then $$[f(\alpha),f(\beta)]=D_1(\alpha)D_2(\beta)X^{n-1}\neq 0,$$ i.e., the ring $A$ is not commutative.~$\rhd$

\textbf{4.2.7. Remark.} Until the end of Section 4.2, we assume that $A$ is a non-simple ring and $\varphi\colon A\rightarrow A$ is an injective homomor\-phism from the ring $A$ into itself such that $\varphi(A\setminus \{0\})\subseteq A^*$.

\textbf{4.2.8. Right skew power series rings.}\\
We denote by $A_r[[x,\varphi]]$ the right skew power series ring in the sense of \cite[9.8]{Tug98}; this ring
consists of all formal series $\sum_{k=0}^{+\infty}x^ka_k$, $a_k\in A$, the addition of series is component-wise and the multiplication is naturally defined with the use of the rule $ax^k=x^k\varphi^k(a)$.

\textbf{4.2.9. Lemma.}\label{two-sided}
Let $A$ be a non-simple ring, $R=A_r[[x,\varphi]]$, and let $I$ be a non-zero two-sided ideal of $R$. Then $I=x^mB+x^{m+1}R$ for some $n\geq 0$ and some non-zero right ideal $B$ of $A$, which is a left $\varphi^m(A)$-module (we assume that $\varphi^0$ is the identity mapping).

$\lhd$ Let $0\neq I\lhd R$. Since $\bigcap_{i=0}^\infty x^mR=0$, there exists an integer $m\geq 0$ such that $I\subseteq x^mR$ and $I\not\subseteq x^{m+1}R$. Let $f\in
I\setminus x^{m+1}R$, then it follows from \cite[9.9(3)]{Tug98}\footnote{There is a misprint in the text of \cite[9.9(3)]{Tug98}: the correct relation is $M\subseteq N\Leftrightarrow \mbox{ either } m>n \mbox{ or } m=n \mbox{ and } D\subseteq E$.} that $fR=x^mDR$ for some non-zero principal right ideal $D$ of the ring $A$; in addition, we set $E=R$, $n=m+1$ and obtain $fR\supseteq x^{n}R$. It remains to set $B=\{b\in A\,|\,x^mb+x^{m+1}R\subseteq I\}$. We multiply the elements of $I$ from the left and right by elements of the ring $A$ and obtain that $B$ is a $(\varphi^m(A),A)$-sub-bimodule in $A$.~$\rhd$

\textbf{4.2.10. Lemma.}\label{center}\\
Let $A$ be a non-simple ring, $R=A_r[[x,\varphi]]$, and let $Z(A)\not \subseteq A^*\cup \{0\}$. Then $xR\cap Z(R)=0$.

$\lhd$ We choose $a\in Z(A)\setminus (A^*\cup \{0\})$. Let $f=\sum_{i=1}^\infty x^if_i \in xR\cap Z(R)$. Then it follows from the relation $[a,f]=0$ that $\varphi^i(a)f_i=f_ia=af_i$ for any $i>0$. The relation $a=\varphi^i(a)$
is impossible for $i>0$, since $a\not\in A^*\cup\{0\}$ and the ring $A$ is a domain; therefore, $f_i=0$ for all $i>0$, i.e., $f=0$.~$\rhd$

\textbf{4.2.11. Proposition.}
\label{PI} Let $A$ be a non-simple PI ring, $R=A_r[[x,\varphi]]$, and let $I$ be an ideal of $R$. Then $R/I$ is a PI ring if and only if $I\neq 0$.

$\lhd$
Let $I\neq 0$. By Lemma 4.2.9, we have that $I\supseteq x^nR=(xR)^n$ for some $n>0$. If $f(x_1,\ldots,x_t)$ is an admissible identity of the ring $A=R/xR$, then the admissible identity
$$
f(x_1,\ldots,x_t)f(x_{t+1},\ldots,x_{2t})\ldots f(x_{(n-1)t+1},\ldots,x_{nt})
$$
holds in the ring $R/I$.

Let $R/I$ is a PI ring.
We have to prove that $R$ is not a PI ring under the conditions of the proposition. The rings $A$ and $R$ are domains; see \cite[9.9(1)]{Tug98}. 

We need the following well known fact $(*)$, see \cite[Theorem 2]{Row73}.\\
If $S$ is a semiprime PI ring, and $I$ is a non-zero two-sided ideal of the ring $S$, then $Z(S)\cap I\neq 0$.

By applying $(*)$ to the proper non-zero ideal $B$ of the ring $A$, we obtain that $0\neq Z(A)\cap B\not\subseteq A^*$. By Lemma 4.2.10, $xR\cap Z(R)=0$; in addition, $xR$ is a non-zero two-sided ideal of the semiprime ring $R$. We again use $(*)$ and see that $R$ cannot be a PI ring.~$\rhd$

\textbf{4.2.12. Proposition.}\label{non-ce}\\ Let $A$ be a commutative non-simple ring, $R=A_r[[x,\varphi]]$, and let $I$ be an ideal of the ring $R$. Then the following conditions are equivalent.

\textbf{a)} $I\supseteq xR$.

\textbf{b)} $R/I$ is a commutative ring.

\textbf{c)} The ring $R/I$ is centrally essential.

$\lhd$ The implications a)\,$\Rightarrow$\,b)$\Rightarrow$\,c) are directly verified. Let's assume that \textbf{c)} holds. Since the ring $A$ is not simple (in the commutative case, this means that $A$ is not a field), we have $A=Z(A)\not \subseteq A^*\cup
\{0\}$, whence $xR\cap Z(R)=0$ by Lemma 4.2.10; this is impossible in a centrally essential ring. Consequently, $I\neq 0$.

By Lemma 4.2.9, we have $I=x^mB+x^{m+1}R$, where $B$ is a non-zero
ideal of the ring $A$. If $m=0$, then \textbf{a)} holds. We assume that $m>0$. For any element $r\in R$, we set $\widehat r=r+I\in \widehat R=R/I$. Since $I\subseteq xR$, we can identify the elements $a$ and $\widehat a$ for any $a\in A$. Then any element of the ring $\widehat R$ can be considered as the
sum $f=f_0+\widehat x_1f_1+\ldots+\widehat x^mf_m$, where $f_i\in A$, $i=0,1,\ldots,m$, the coefficients $f_0,\ldots,f_{m-1}$ are uniquely defined and $f_m$ is determined up to a summand which is an arbitrary element of $B$. Let $f\in Z(\widehat R)$. It follows from the relation $[a,f]=0$ (where $a\in A$) that $\varphi^i(a)f_i=f_ia=af_i$ for any $i=0,1,\ldots,m-1$. If $a$ is a non-zero non-invertible element of the ring $A$, we obtain $a\neq \varphi^i(a)$ for $i>0$, whence $f_1=\ldots=f_{m-1}=0$. By \textbf{c)}, there exist two elements $c,d\in Z(\widehat R)$ such that $\widehat x c=d\neq 0$. We set $c=c_0+\widehat x^mc_m$. Then $d=\widehat xc=\widehat x c_0$, since $\widehat
x^{m+1}=0$. First, we assume that $m>1$. Then it follows from the inclusion $d\in Z(\widehat R)$ that $c_0=0$ and $d=0$; this is a contradiction. Thus,
$m=1$.
Let $c=c_0+\widehat xc_1\in Z(\widehat R)$. For
$b\in B\setminus \{0\}$, it follows from the relation $[c,b]=0$ that $\widehat x(c_1b-c_1\varphi(b))=0$; therefore,
$c_1b-c_1\varphi(b)\in B$ and $c_1\varphi(b)\in B$ and we have $c_1\in B$, since $\varphi(b)$ is invertible, i.e., $c=c_0$.
We assume that $B$ is a proper ideal of the ring $A$. Then $\widehat x\neq 0$ and for some $c_0\in Z(\widehat R)$, we have $d=\widehat x c_0\in Z(\widehat R)\setminus \{0\}$.
Similar to the above case with $c$, it follows from the relation $[d,b]=0$ that $c_0\in B$, i.e., $d=0$. This contradiction shows that $B=A$ and $I=xA+x^2R=xR$.~$\rhd$

\textbf{4.2.13. Example.}\\
There exists a right uniserial, right Noetherian, non-semiprimary PI ring $\widehat R$ with prime radical $\widehat N$ and Jacobson radical $\widehat M$ such that $\widehat R$ is not centrally essential, and $\widehat R$ is not left Noetherian or left uniserial, $\widehat R/\widehat N$ is a commutative discrete valuation domain, $\widehat N$ is a minimal right ideal, and $\widehat N=\widehat M\widehat N\ne \widehat N\widehat M=0$.

$\lhd$ Let $A$ be a non-simple ring and $R=A_r[[x,\varphi]]$. We use the example \cite[9.10]{Tug98}, where $N=xR$ and $0\neq NM\neq N$. Then $\widehat R=R/(NM)$ is a PI ring by Proposition 4.2.11; however, it is not a centrally essential ring by Proposition 4.2.12.~$\rhd$

\textbf{4.2.14. Remark.} With the use of the above results of Section 4.2, it is easy verified that $J^{n-1}\subseteq C$ under the conditions of 4.2.13. On the
other hand, the following example shows that the inclusion
$\displaystyle{J^{\left[\dfrac n 2\right]+1}} \subseteq C$ does not necessarily imply that the left uniserial, left Artinian ring $R$ is centrally essential.

\textbf{4.2.15. Example.}\\ 
Let $F=GF(4)$, $F_0=GF(2)\subseteq F$, and let $\sigma\colon x\mapsto x^2$ be the Frobenius automor\-phism of the field $F$. We consider the skew polynomial ring
$S=F[X,\sigma]$ and its factor ring $R=S/(X^3)$. Then $R$ is a left and right uniserial ring, left and right Artinian ring, $J(R)$ is a nilpotent ideal of nilpotence index 3, and $\displaystyle{J(R)^{\left[\dfrac 3 2\right]+1}}\subseteq Z(R)$; however, $R$ is not centrally essential.

$\lhd$ It is clear that $F=F_0[\theta]$, where $\theta$ is a root of the irreducible polynomial $t^2+t+1\in F_0[t]$. We denote by $x$ the image of the variable $X$ under the canonical homomor\-phism from the ring $S$ onto $R$ and identify the
elements of the field $F$ with their images in $R$. It is directly verified that $J(R)=(x)$, $n=3$ is the nilpotence index of the ideal $J=J(R)$ and the left (and right) modules $J/J^2$ and $J^2$ are one-dimensional vector spaces over $F=R/J$. Consequently, $R$ is a left and right uniserial, left and right
Artinian ring by 4.1.1. We consider the element
$r=a_0+a_1x+a_2x^2$. It follows from the relation $x^3=0$ that
$$
[r,x]=rx-xr=(a_0-\sigma(a_0))x+(a_1-\sigma(a_1)x^2 \text{ and}
$$
$$
[r,\theta]=(a_1\sigma(\theta)-a_1\theta)x+(a_2\sigma^2(\theta)-\theta
a_2)x^2=a_1x,
$$
since $\sigma^2$ is the the identity automor\-phism and $\sigma(\theta)=\theta+1$.
Therefore, $Z(R)=F_0+Fx^2$, since $x$ and $\theta$ generate the ring $R$ (as a ring). It remains to note that $Z(R)x=F_0x$ and $F_0x\cap (F_0+Fx^2)=0$.~$\rhd$

\subsection{Rings with Flat Ideals}\label{section4.3}

The results of this subsection based on \cite{Tug21}.

\textbf{4.3.1. Rings of weak global dimension $\le 1$.}\\
For a ring $R$, we write $\text{w.gl.dim. }R\le 1$ if $R$ is a \textsf{ring of weak global dimension at most one}\label{wgd1} i.e., $R$ satisfies the following equivalent\footnote{The equivalence of the conditions is well known; e.g., see \cite[Theorem 6.12]{Tug98}.} conditions.
\begin{itemize}
\item
For every finitely generated right ideal $X$ of $R$ and each finitely generated left ideal $Y$ of $R$, the natural the group homomor\-phism $X\otimes_RY\to XY$ is an isomor\-phism.
\item
Every finitely generated right (resp., left) ideal of $R$ is a flat\footnote{A right $R$-module $X$ is said to be \textsf{flat}\label{flat} if for any left $R$-module $Y$, the natural the group homomor\-phism $X\otimes Y\to XY$ is an isomor\-phism.} right (resp., left) $R$-module.
\item
Every right (resp., left) ideal of $R$ is a flat right (resp., left) $R$-module.
\item
Every submodule of any flat (right or left) $R$-module is flat.
\end{itemize}
Since every projective module is a flat module, any right or left\\ (semi)hereditary\footnote{A module $M$ is said to be \textsf{hereditary} (resp., \textsf{semihereditary})\label{semiher} if all submodules (resp., finitely generated submodules) of $M$ are projective.} ring is of weak global dimension at most one. 
We also recall that a ring $R$ is of \textsf{weak global dimension zero}\label{wgd0} if and only if $R$ is a \textsf{Von Neumann regular}\label{vnrring} ring, i.e., $r\in rRr$ for every element $r$ of $R$. Von Neumann regular rings are widely used in mathematics; see \cite{Goo79}, \cite{Kap68}. 

\textbf{4.3.2. Theorem; \cite[Theorem]{Jen64}.}\label{theorem 1} \\
A commutative ring $R$ is a ring of weak global dimension at most one if and only if $R$ is an arithmetical semiprime ring.

It is clear that a commutative ring is right (resp., left) distributive if and only if the ring is arithmetical.

\textbf{4.3.3. Example.}\\
There exists a right hereditary ring $R$, which is neither right distributive nor semiprime; in particular, the right hereditary ring $R$ is of weak global dimension at most one.\\ 
Let $F$ be a field and let $R$ be the $5$-dimensional $F$-algebra consisting of all $3\times 3$ matrices of the form
$\left(\begin{array}{ccc}
f_{11}& f_{12}& f_{13} \\
0& f_{22}& 0 \\
0& 0& f_{33}
\end{array}\right)$, where $f_{ij}\in F$.
The ring $R$ is not semiprime, since the set $\left\{\left(\begin{array}{ccc}
0& f_{12}& f_{13} \\
0& 0& 0 \\
0& 0& 0
\end{array}\right)\right\}$ is a non-zero nilpotent ideal of $R$.
Let $e_{11}$, $e_{22}$ and $e_{33}$ be ordinary matrix units.
The ring $R$ is not right or left distributive, since every idempotent of a right or left distributive ring is central, see \cite{Ste74}, but the matrix unit $e_{11}$ of $R$ is not central. To prove that the ring $R$ is right hereditary, it is sufficient to prove that $R_R$ is a direct sum of hereditary right ideals. We have that $R_R=e_{11}R\oplus e_{22}R\oplus e_{33}R$, where $e_{22}R$ and $e_{33}R$ are projective simple $R$-modules; in particular, $e_{22}R$ and $e_{33}R$ are hereditary $R$-modules. Any direct sum of hereditary modules is hereditary; see \cite[39.7, p.332]{Wis91}. Therefore, it remains to show that the $R$-module $e_{11}R=e_{11}F+e_{12}F+e_{13}F$ is hereditary, which is directly verified.

The following lemma is well known; e.g., see \cite[Assertion 6.13]{Tug98}. 

\textbf{4.3.4. Lemma.} Let $R$ be a ring whose principal right ideals are flat. If $r$ and $s$ be two elements of $R$ with $rs=0$, then there exist two elements $a,b\in R$ such that $a+b=1$, $ra=0$, and $bs=0$. 

\textbf{4.3.5. Lemma.} There exists a right and left uniserial prime ring $R$, which has a non-flat principal right ideal.

$\lhd$
There exists a right and left uniserial prime ring $R$ with two non-zero elements $r,s\in R$ such that $rs=0$; see \cite[p.~234, Corollary]{Dub94}. The uniserial ring $R$ is local; therefore, the non-invertible elements of $R$ form the Jacobson radical $J(R)$ of $R$. The ring $R$ is not a ring whose principal right ideals are flat. Indeed, let's assume the contrary. By Lemma 4.3.4, there exist two elements $a,b\in R$ such that $a+b=1$, $ra=0$, and $bs=0$. We have that either $aR\subseteq bR$ or $bR\subseteq aR$; in addition, $aR+bR=R=Ra+Rb$. Therefore, at least one of the elements $a,b$ of the local ring $R$ is invertible; in particular, this invertible element is not a right or left zero-divisor. This contradicts to the relations $ra=0$ and $bs=0$.~$\rhd$

\textbf{4.3.6. Lemma.} Let $R$ be a centrally essential ring whose principal right ideals are flat. Then the ring $R$ does not have non-zero nilpotent elements. 

$\lhd$ Indeed, let's assume that there exists a non-zero element $r\in R$ with $r^2=0$. Since the ring $R$ is centrally essential, there exist two non-zero central elements $x,y\in R$ with $rx = y$. Since $r^2=0$, we have that $y^2=(rx)^2=r^2x^2=0$. Since $y^2=0$, it follows from Lemma 4.3.4 that there exist two elements $a,b\in R$ such that $a+b=1$, $ry=0$, and $by=yb=0$. Then $y=y(a+b)=ya+yb=0$. This is a contradiction.~$\rhd$

\textbf{4.3.7. Theorem.}\label{theorem 2}\\
For a centrally essential ring $R$, the following conditions are equivalent.
\begin{enumerate}
\item[\textbf{a}]
$R$ is a ring of weak global dimension at most one.
\item[\textbf{b}]
$R$ is a right (resp., left) distributive semiprime ring.
\item[\textbf{c}]
$R$ is an arithmetical semiprime ring
\end{enumerate}

$\lhd$
\textbf{a)}\,$\Rightarrow$\,\textbf{b)}. Since $R$ is a centrally essential ring of weak global dimension at most one, it follows from Lemma 4.3.6 that the ring $R$ does not have non-zero nilpotent elements. By Theorem 1.2.2, the centrally essential semiprime ring $R$ is commutative. By Theorem 4.3.2, $R$ is an arithmetical semiprime ring. Any commutative arithmetical ring is right and left distributive.

The implication \textbf{b)}\,$\Rightarrow$\,\textbf{c)} follows from the property that every right or left distributive ring is arithmetical.

\textbf{c)}\,$\Rightarrow$\,\textbf{a)}. Since $R$ is a centrally essential semiprime ring, it follows from Theorem 1.2.2 that the ring $R$ is commutative; in particular, $R$ is centrally essential. In addition, $R$ is arithmetical. By Theorem 4.3.2, the ring $R$ is of weak global dimension at most one.~$\rhd$

\textbf{4.3.8. Remark.}\\ 
It follows from Lemma 4.3.5 that the implication \textbf{b)}\,$\Rightarrow$\,\textbf{a)} of Theorem 4.3.7 is not true for arbitrary rings.

\textbf{4.3.9. Corollary.}\\
A ring $R$ is a right (resp., left) hereditary, right (resp., left) Noetherian, centrally essential ring if and only if $R$ is a finite direct product of commutative Dedekind domains.\\
Consequently, a ring $R$ is a right (resp., left) hereditary, right (resp., left) Noetherian, indecomposable, centrally essential ring if and only if $R$ is a commutative Dedekind domain.

$\lhd$ 
Since the right or left Noetherian ring $R$ is a finite direct product of right or left Noetherian rings, we can assume that $R$ is a right or left Noetherian indecomposable ring. In this case, it well known that $R$ is a commutative hereditary ring if and only if $R$ is a commutative Dedekind domain. Now we can use Theorem 4.3.7 and the well known property that every flat module over a Noetherian ring is projective.~$\rhd$

\subsection{Distributive Noetherian Rings}\label{section4.4}

The results of this subsection are based on \cite{Tug23}.

\textbf{4.4.1. Distributive Noetherian semiprime rings.}\\
A ring $R$ is a right (resp., left) distributive, right (resp., left) Noetherian, semiprime, centrally essential ring if and only if $R$ is a finite direct product of commutative Dedekind domains. \\
Consequently, a ring $R$ is a right (resp., left) distributive, right (resp., left) Noetherian, indecomposable, centrally essential ring if and only if $R$ is a commutative Dedekind domain.

$\lhd$ 
The assertion follows from Corollary 4.3.9 and the following well known property: a commutative ring is a Dedekind domain if and only if $R$ is a commutative distributive Noetherian domain.~$\rhd$

\textbf{4.4.2. The notation $\ld$ in 4.4.2, 4.4.3, and 4.4.5.}\\
Let $A$ be a ring, $X$ be a right $A$-module, and let $X_1$, $X_2$ be two subsets in $X$.\\
We denote by $(X_1\ld X_2)$\label{ld} the subset $\{a\in A\,\mid \, X_1a\subseteq X_2\}$ of the ring $A$. If $X_2$ is a submodule in $X$, then $(X_1\ld X_2)$ is a right ideal of the ring $A$. If $X_1$ and $X_2$ are two submodules in $X$, then $(X_1\ld X_2)$ is an ideal in $A$.\\

We use some familiar properties of distributive modules and rings. For convenience of readers these properties are gathered in Lemmas 4.4.3, 4.4.4, 4.4.5. 

\textbf{4.4.3. Lemma; \cite{Ste74}.}\\
Let $A$ be a ring and let $X$ be a distributive right $A$-module. 

\textbf{a.} For any two elements $x,y\in X$, there exist elements $a,b\in A$ such that $a+b=1$ and $xaA+ybA\subseteq xA\cap yA$. Consequently, $A=(x\ld yA)+(y\ld xA)$ for any elements $x,y\in X$. \\ In particular, if $xA\cap yA=0$, then there exist elements $a,b\in A$ such that $a+b=1$ and $xaA=ybA=0$.

\textbf{b.} $\text{Hom}\,(Y,Z)=0$ for any submodules $Y,Z$ of the module $X$ such that $Y\cap Z=0$.

\textbf{c.} All idempotents of the ring $\text{End}\,M$ are central.\\
In particular, all idempotents of any right distributive ring are central. Therefore, the right distributive ring $A$ is indecomposable into the ring direct product of if and only if $A$ does not have non-trivial idempotents.

\textbf{d.} If the ring $A$ is local, then $M$ is a uniserial module.\\
In particular, right distributive local rings coincide with right uniserial rings.

\textbf{e.} If $M$ is a Noetherian module, then $M$ is an invariant\footnote{A module $M$ is said to be \textbf{invariant}\label{invMod} if every submodule of $M$ is fully invariant in $M$.} module.\\
In particular, any right distributive right Noetherian ring is right invariant.

$\lhd$ \textbf{a.} We set $T=xA\cap yA$. Since
$(x+y)A=(x+y)A\cap xA+(x+y)A\cap yA$, there exist elements $b,d\in A$ such that
$$
(x+y)b\in xA,\; (x+y)d\in yA,\; x+y=(x+y)b+(x+y)d.
$$
Therefore, $yb=(x+y)b-xb\in T$ and $xd=(x+y)d-yd\in T$. We set $a=1-b$ and $z=a-d=1-b-d$. Then
$$
1=a+b,\; (x+y)z=(x+y)-(x+y)b-(x+y)d=0,
$$
$$
xa=xd+xz=xd+(x+y)z-yz=xd-yz,
$$
$$
yz=-xz\in T, \; xa\in T.
$$
\textbf{b.} Let $f\in \text{Hom}\,(Y,Z)$, $y\in Y$ and $z = f(y)\in Z$.
By \textbf{a}, there exists an element $a\in A$ such that
$$
yaA+z(1-a)A\subseteq yA\cap zA\subseteq Y\cap Z = 0,
$$
$$
ya = z(1-a) = 0, \; z = za = f(y)a = f(ya) = f(0) = 0.
$$
Therefore, $f\equiv 0$ and $\text{Hom}(X,Y)=0$.

\textbf{c.} With the use of \textbf{b}, the assertion is directly verified.

\textbf{d.} Let $x,y\in X$. It is sufficient to prove that the submodules $xA$ and $yA$ are comparable with respect to inclusion. By \textbf{a}, there exist elements $a,b,c,d\in A$ such that $1=a+b$ and $xaA+ybA\subseteq xA\cap yA$. Since the ring $A$ is local and $1=a+b$, at least one of the right ideals $aA$, $bA$ coincides with $A$. Therefore, at least one of inclusions $xA\subseteq yA$, $yA\subseteq xA$ holds.

\textbf{e.} The assertion is proved in \cite{Ste74}.~$\rhd$

The following Lemma 4.4.4 is a direct corollary of Lemma 4.4.3(d) and \cite[Proposition 2]{Cam75}

\textbf{4.4.4. Lemma; see \cite[Proposition 2]{Cam75}.}\\
$A$ is a right distributive, right Noetherian, semiprime ring if and only if $A$ is a finite direct product of right distributive, right Noetherian, right invariant domains.

\textbf{4.4.5. Lemma; see \cite[Lemma 20]{Tug95}.}\\
Let $A$ be a right invariant ring and let $M$ be a distributive right $A$-module.

\textbf{a.} $A = (Y\ld X)+(X\ld Y)$ for any finite generated submodules $X$, $Y$ of the module $M$.

\textbf{b.} For any submodule $Z'$ of an arbitrary finitely generated submodule $Z$ of the module $M$, there exists
an ideal $A'$ of the ring $A$ such that $ZA'=Z'$.

\textbf{c.} If $M$ is a finite generated module, then $M$ is an invariant module.

$\lhd$ \textbf{a.} Since $X + Y$ is a finitely generated module, there exist a positive integer $n$ and elements $x_i\in X$, $y_i\in Y$, $1\le i\le n$ such that $X + Y = (x_1+ y_1)A+\dots + (x_n+ y_n)A$. Since the module $M$ is distributive, $(X + Y)\bigcap Z = (X\bigcap Z) + (Y\bigcap Z)$ for any submodule $Z$ in $M$. 

Let $y\in Y$. For any $1\le i\le n$, we have 
$$
(x_i + y)A = (x_i + y)A\bigcap (X + Y) =
$$
$$
=[(x_i + y)A\bigcap X] + ((x_i + y)A\bigcap Y]. 
$$
Therefore, there exist elements $a\in A$ and $z\in Y$ such that 
$$
(x_i + y)a\in X,\quad x_i+ y = (x_i+ y)a + z. 
$$
Therefore, $x_i(1-a)\in Y$ and $ya\in X$. Consequently,
$$
A = (yA\ld X) +(x_iA\ld Y), \quad 1\le i\le n.
$$ 
Therefore, 
$$
A = (yA\ld X) +[(x_1A\ld Y)\bigcap \dots \bigcap (x_nA\ld Y)] = 
$$
$$
=(yA\ld X) + (X\ld Y).
$$ 
In particular, 
$$
A = (y_iA\ld X) + (X\ld Y) \quad (1\le i\le n). 
$$
Therefore, 
$$
A = [(y_1A\ld X)\bigcap \dots \bigcap (y_nA\ld X)] + (X\ld Y) =
$$
$$
=(Y\ld X) + (X\ld Y).
$$
\textbf{b.} Let $Z$ be an $n$-generated module, $n\in \mathbb{N}$. We use the induction on $n$. For $n=1$, we can identify the cyclic $A$-module $Z$ over the right invariant ring $A$ with the right invariant factor ring $A/\text{r}_A(Z)$ of the ring $A$. In this case, the assertion is directly verified.

Now we assume that the assertion is true for all $k$-generated submodules of the module $M$ for $k<n$. We can assume that $Z=X+Y$, where $X$ is a cyclic module and $Y$ is an $(n-1)$-generated module. By the induction hypothesis, there exist ideals $B$ and $C$ of the ring $A$ such that $X\cap Y=XB=YC$. Therefore, $X\bigcap Y=X(X\ld Y)=Y(Y\ld X)$
By \textbf{a}, $A=(Y\ld X) +(X\ld Y)$. Therefore,
$$
X=X((Y\ld X)+(X\ld Y)) =
$$
$$
= X(Y\ld X)+X(X\ld Y)= X(Y\ld X)+Y(Y\ld X)=ZB,
$$
where $B=(Y\ld X)$. Similarly $Y=ZC$, where $C=(X\ld Y)$.

Let $Z'$ be a submodule in $Z=X+Y$. We have to prove that there exists an ideal $H$ of the ring $A$ such that $Z'=(X+Y)H$. 
By assumption, $Z'= X\cap Z'+Y\bigcap Z'$. By the induction hypothesis, there exist ideals $D$ and $E$ of the ring $A$ such that $Z'\bigcap X=XD$ and $Z'\bigcap Y=YE$. In addition, $X=ZB$ and $Y=ZC$. Therefore, 
$$
Z'=X\bigcap Z'+Y\bigcap Z' = XD+YE = 
$$
$$
=ZBD + ZCE = Z(BD+CE)$$ and $BD+CE$ is the required ideal $A'$ of the ring $A$.~$\rhd$

\textbf{4.4.6. Lemma; see \cite[Proposition 2]{Tug95}.} For a ring $A$, the following conditions are equivalent.

\textbf{a)} $A$ is a right distributive, right Noetherian, left finite-dimensional semiprime ring.

\textbf{b)} 
$A$ is a left distributive, left Noetherian, right finite-dimensional semiprime ring.

\textbf{c)} $A$ is a finite direct product of invariant hereditary Noetherian domains.

\textbf{4.4.7. Lemma.}\\
Let $A$ be a right distributive, right Noetherian semiprime ring such that every non-zero left ideal of $A$ contains a non-zero central element. Then $A$ is a finite direct product of invariant hereditary Noetherian domains.

$\lhd$ 
By Lemma 4.4.6, it is sufficient to prove that the ring $A$ is left finite-dimensional. We assume the contrary. Then the ring $A$ contains a left ideal $B$, which is a countable direct sum of non-zero left ideals $B_k$, $k=1,\dots,+\infty$. By assumption, every left ideal $B_k$ contains non-zero central element $c_k$, where the sum of all ideals $Ac_k=c_kA$ is a direct sum. This contradicts to the property that the ring $A$ is right finite-dimensional.~$\rhd$ 

\textbf{4.4.8. Proposition.}\\
Let $A$ be a right distributive right Noetherian indecomposable ring with the prime radical $M$. Then the ring $A$ is right invariant, $A/M$ is a right distributive, right invariant, Noetherian domain, $M$ is a completely prime Noetherian nilpotent right ideal. In addition, the following assertions are true.

\textbf{a.} $M=xM$ for any element $x\in A\setminus M$.

\textbf{b.} For any submodule $N$ of the module $M_A$, there exists an ideal $D$ of the ring $A$ such that $N=MD=xN$ for any element $x\in A\setminus M$.

\textbf{c.} If $M$ contains a non-zero central element $m$, then $A$ is a right uniserial right Artinian ring with radical $M$.

\textbf{d.} If every non-zero left ideal of the ring $A$ contains a non-zero central element, then either $A$ is an invariant hereditary Noetherian domain or $A$ is a right uniserial, right Artinian ring.

$\lhd$ By Lemma 4.4.3(c), the ring $A$ is right invariant. Since $M$ is the prime radical of the right Noetherian ring $A$, the nil-ideal $M$ is nilpotent. Since $M$ is a nil-ideal, the idempotents of the factor ring $A/B$ are lifted to idempotents of the ring $A$. By Lemma 4.4.3(a), the indecomposable ring $A$ does not have non-trivial idempotents. Therefore, the factor ring $A/M$ does not have non-trivial idempotents. By Lemma 4.4.4, $A/B$ is a right distributive, right invariant, right Noetherian domain. Therefore, the Noetherian nilpotent right ideal $M$ is completely prime.

\textbf{a.} Let $x\in A\setminus M$ and let $y$ be an arbitrary element of the ideal $M$. By Lemma 4.4.2, there exist two elements $a,b\in A$ such that $a+b=1$, $xa\in yA$ and $yb\in xA$. The ideal $M$ is completely prime, $x\in A\setminus M$ and $xa\in M$. Therefore, $a\in M$ and element $a$ is nilpotent. Therefore, the element $b=1-a$ is invertible; in addition, $yb\in xA$. Then $y=xz$ for some $z\in A$. Since the element $xz$ is contained in the completely prime ideal $M$ and $x\in A\setminus M$, we have that $z\in M$ and $y=xz\in xM$. Since the element $y\in M$ is arbitrary, we have $M=xM$.

\textbf{b.} Let $N$ be a submodule of the module $M_A$ and let $x\in A\setminus M$. By Lemma 4.4.5, there exists an ideal $D$ of the ring $A$ with $MD=N$. In addition, $M=xM$. Therefore, $N=xMD=xN$.

\textbf{c.} 
Let $M$ contain a non-zero central element $m$. Since the ring $A$ is right invariant and $m$ is a non-zero central element, there exists a maximal ideal $X$ of the ring $A$ such that $A/X$ is a division ring and the ideal $mA$ properly contains the ideal $(mA)X=X(mA)$. If $X=M$, then the ring $A$ is local. By Lemma 4.4.3(b), $A$ is a right uniserial ring. In addition, $A$ is a right Noetherian ring with nilpotent Jacobson radical $M$. Therefore, $A$ is a right uniserial right Artinian ring.

\textbf{d.} If the completely prime ideal $M$ is equal to the zero, then $A$ is a domain and $A$ is an invariant hereditary Noetherian domain, by Lemma 4.4.7. We assume that $M\ne 0$. By assumption, the ideal $M$ contains a non-zero central element. It follows from \textbf{c} that $A$ is a right uniserial right Artinian ring.~$\rhd$

\textbf{4.4.9. Theorem.}\\
A centrally essential ring $A$ is a right distributive, right Noetherian ring if and only if $A=A_1\times\cdots\times A_n$, where every ring $A_k$ is either a commutative Dedekind domain or a (not necessarily commutative) Artinian uniserial ring.

$\lhd$ Let $A$ be a centrally essential ring.
If $A=A_1\times\cdots\times A_n$, where every ring $A_k$ is either a commutative Dedekind domain or a (not necessarily commutative) Artinian uniserial ring, then the assertion follows from Lemma 4.4.6.

Now let $A$ be a right distributive, right Noetherian ring. Without loss of generality, we can assume that $A$ is an indecomposable ring. Since the ring $A$ is centrally essential, every non-zero left or right ideal of the ring $A$ contains a non-zero central element. It follows from Proposition 4.4.8(d) that either $A$ is a right uniserial, right Artinian ring or $A$ is an invariant hereditary Noetherian domain. If $A$ is a right uniserial, right Artinian ring, then $A$ is a uniserial Artinian ring, by Lemma 4.2.4. If $A$ is a domain, then $A$ is a commutative Dedekind domain by 4.4.1.~$\rhd$

\section[Centrally Essential Semirings]{Centrally Essential Semirings}\label{section5}

\textsf{In Section 5, we consider only unital semirings and rings}, i.e., semirings and rings with $1$. 

In this section, we consider centrally essential semirings. 
Some semiring notions are defined below.
Another required information n semirings is contained in the book \cite{HebW98}. 

\subsection{General Information}\label{subsection5.1}

\textbf{5.1.1. Unital semirings and their centers.}\\
By a \textsf{semiring},\label{semiRi} we mean a structure that differs from an associative ring, possibly, by the irreversibility of the additive operation. In a semiring $S$, the zero is multiplicative by definition: we have $0s = s0 = 0$ for every $s\in S$. 

For a semiring $S$, the \textsf{center}\label{cenSemir} of $S$ is the set $Z(S) = \{s\in S\colon ss' = s's$ for all $s'\in S\}$. This set is not empty, since it contains $0$ and $1$; we also have that $Z(S)$ is a subsemiring in $S$. 

\textbf{5.1.2. Centrally essential semirings.}\\
A semiring $S$ is said to be \textsf{centrally essential}\label{cenEsSemi} if either $S$ is commutative or for every non-zero element $s\in S$, there are non-zero central elements $x,y$ with $sx = y$.

It is clear that any centrally essential associative ring is a centrally essential semiring. 

\textbf{5.1.3. Reduced semirings and non-zero-divisors.}\\
A semiring $S$ is said to be \textsf{reduced}\label{redSemi} if $x = y$ for all $x, y\in S$ with $x^2 + y^2 = xy + yx$. If $S$ is a ring, this is equivalent the property that $S$ has no non-zero nilpotent elements.

An element $a$ of a semiring $S$ is said to be \textsf{left} (resp., right) \textsf{zero-divisor} if $ab = 0$ (resp., $ba = 0$) for some $0\neq b\in S$. Similar to 1.1.2(a), it can be proved that one-sided zero-divisors are two-sided zero-divisors in a centrally essential semiring.

\textbf{5.1.4. Semiprime and semisubtractive semirings.}\\
A semiring $S$ without nilpotent ideals is said to be \textsf{semiprime}.\label{semisemipr} A semiring $S$ is said to be \textsf{semisubtractive}\label{semisubsSemi} if for all $a, b\in S$ with $a\neq b$, there exists an element $x\in S$ such that $a + x = b$ or $b + x = a$.

\textbf{5.1.5. Additively cancellative semirings and\\ rings of differences.}\\
A semiring $S$ is said to be \textsf{additively cancellative}\label{addCan} if the relation $x + z = y + z$ is equivalent to the relation $x = y$ for all $x, y, z\in S$. 

A ring $D(S)$ is called the \textsf{ring of differences}\label{ringDif} of the semiring if $S$ is a subsemiring in $D(S)$ and every element $a\in D(S)$ is the difference $x - y$ of some elements $x,y\in S$. 

It is well known that a semiring $S$ can be embedded in the ring of differences $D(S)$ if and only if $S$ is additively cancellative.

The class of additively cancellative semirings contains all rings. The ring of differences is unique up to isomor\-phism over $S$; see \cite[Chapter II]{HebW98} for details. 

\textbf{5.1.6. Idempotents of centrally essential semirings.}\\
By 1.1.4, idempotents of centrally essential rings are central. For semirings, a similar result is not true; see Example 5.2.2 below. 

For a semiring $S$, an idempotent $e$ of $S$ is said to be \textsf{complemented}\label{compIde} if there exists an idempotent $f\in S$ with $e + f = 1$.

\textbf{5.1.7. Proposition.}\\
In an additively cancellative centrally essential semiring $S$, any complemented idempotent is central. 

$\lhd$ Let $e^2 = e$ and $e + f = 1$ for some $f\in S$. Since $S$ is an additively cancellative semiring, it follows from $e = e + fe$ that $fe = 0$. Similarly, we have $ef = 0$. Let $x\in S$ and $xe\neq 0$. Then $x = ex + fx$ and $xe = exe + fxe$. 

First, we assume that $fxe = 0$, i.e., $xe = exe$. Since $x = xe + xf$, we have $ex = exe + exf$. If $exf\neq 0$, then there are $c, d\in Z(S)$ 
with $(exf)c =d\neq 0$. Then 
$$
0\neq d = ed = de = (exfc)e = (exc)fe = 0;
$$
this is a contradiction. Therefore, $exf = 0$ and $ex = xe = exe$.

Now let $fxe\neq 0$. Then $0\neq (fxe)c = d$ for some non-zero elements $c, d\in Z(S)$. In this case, 
$$
0\neq d = de = ed = ef(xec) = 0;
$$
this is a contradiction.~$\rhd$

\textbf{5.1.8. Remark.} If $S$ is an additively cancellative semiring, then the semiring $M_n(S)$ of all matrices and the semiring $T_n(S)$ of all upper triangular matrices over $S$ is not centrally essential for $n\ge 2$.

$\lhd$ For the identity matrices of the above semirings, we have $E = E_{11} + \ldots + E_{nn}$, where $E_{11}, \ldots, E_{nn}$ are matrix units. It follows from \cite[Example 4.19]{Gol99} that $M_n(S)$ is an additively cancellative semiring. The idempotents $E_{11}, \ldots, E_{nn}$ are non-central complemented idempotents. Consequently, the semirings $M_n(S)$ and $T_n(S)$ are not centrally essential.~$\rhd$

\subsection{Examples, Constructions and Remarks}\label{subsection5.2}

\textbf{5.2.1. Proposition.}\\ Let $S$ be an additively cancellative semisubtractive centrally essential semiring with center $C=Z(S)$. The following conditions are equivalent.
\begin{itemize}
\item
$S$ is a semiprime semiring.
\item
$C$ is a semiprime semiring.
\item
$S$ does not have non-zero nilpotent elements.
\item
$S$ is a commutative semiring without non-zero nilpotent elements.
\end{itemize}

$\lhd$ By 5.1.5, the semiring $S$ can be embedded in the ring of differences $D(S)$. In addition, the relation $D(S) = -S\cup S$ holds if and only if $S$ is a semisubtractive semiring; see \cite[Chapter II, Remark 5.12]{HebW98}. Then the assertion follows from Theorem 1.2.2.~$\rhd$

\textbf{5.2.2. Example.}\\
We consider a semigroup $(M, \cdot)$ with multiplication table
\begin{center}
\begin{tabular}{|c|c|c|c|c|}
 \hline
$\cdot$ & $1$ & $a$ & $b$ & $c$ \\
\hline
$1$ & $1$ & $a$ & $b$ & $c$ \\
\hline
$a$ & $a$ & $a$ & $a$ & $c$ \\
\hline
$b$ & $b$ & $b$ & $b$ & $c$ \\
\hline
$c$ & $c$ & $c$ & $c$ & $c$ \\
 \hline
\end{tabular}\;.
\end{center}
For a quick test of associativity, it is convenient to use the Light's associativity test; see \cite[p.7]{CliP61}. 

Let $S = 2^M$ be the set of all subsets of the semigroup $M$. For any $A, B\in S$, operations $A + B = A\cup B$ and $AB = \{ab\,|\, a\in A, b\in B\}$ are defined; then $S$ is a semiring with zero $\emptyset$ and the identity element $1 = 1_M$; see \cite[Example 1.10]{Gol99}. 
We have $|S| = 2^4 = 16$. We note that $S$ does not contain zero sums, i.e., the relation $A + B = \emptyset$ implies the relation $A = B = \emptyset$. In addition, $S$ is additively idempotent and multiplicatively idempotent. The center $Z(S)$ is of the form
$$
Z(S) = \{\emptyset, \{1\}, \{c\}, \{1, c\}\}.
$$
If $A\in S\backslash Z(S)$, then $\emptyset\neq A\cdot \{c\}\in Z(S)$. Consequently, $S$ is a non-commutative centrally essential semiring.~$\rhd$

\textbf{5.2.3. Remark.}\\
It follows from Example 5.2.2 that the assertion of Proposition 5.2.1 is not true without the assumptions of additive cancellativity and semisubtractivity.

\textbf{5.2.4. Example.}\\
We consider the semiring $S$ generated by the matrices
$$
\left(\begin{matrix}
\alpha & a & b\\
0 & \alpha & c\\
0 & 0 & \alpha\\
\end{matrix}\right), \left(\begin{matrix}
0 & 0 & b\\
0 & 0 & 0\\
0 & 0 & 0\\
\end{matrix}\right), \left(\begin{matrix}
0 & 0 & 0\\
0 & 0 & 0\\
0 & 0 & 0\\
\end{matrix}\right), \left(\begin{matrix}
\alpha & 0 & 0\\
0 & \alpha & 0\\
0 & 0 & \alpha\\
\end{matrix}\right), 
$$
where $\alpha, a, b, c\in \mathbb{Z}^{+}$.
Let $A = (a_{ij})$ and $B = (b_{ij})$, where $a_{12} = b_{23} = a$, $b_{12} = a_{23} = c$, $a\neq c$, and the remaining components are equal to each other. 
Then $AB\neq BA$, i.e., $S$ is a non-commutative semiring. It is directly verified that the center $Z(S)$ consists of matrices of the form 
$$
\left(\begin{matrix}
\alpha & 0 & b\\
0 & \alpha & 0\\
0 & 0 & \alpha\\
\end{matrix}\right),
$$
where $\alpha, b\in \mathbb{Z}^{+}\cup \{0\}$. Since $0\neq AD\in Z(S)$, where $0\neq A\in S\backslash Z(S)$, $0\neq D\in Z(S)$ with $\alpha = 0$, we have that $S$ is a non-commutative centrally essential semiring. However, the ring of differences $D(S) = M_3(\mathbb{Z})$ is not a centrally essential ring, since the ring has non-central idempotents. In addition, by Remark 3.6.4, any centrally essential subalgebra of a local triangular $3\times 3$ matrix algebra is commutative.~$\rhd$

We give an example of a centrally essential ring $R$, which is the ring of differences for two proper subsemirings $S_1$ and $S_2$ of $R$ such that $S_1$ is not a centrally essential semiring and $S_2$ is a centrally essential semiring.

\textbf{5.2.5. Example.}\\
Let $R$ be a ring consisting of matrices of the form 
$$
\left(\begin{matrix}
\alpha & a & b & c & d & e & f\\
0 & \alpha & 0 & b & 0 & 0 & d\\
0 & 0 & \alpha & 0 & 0 & 0 & e\\
0 & 0 & 0 & \alpha & 0 & 0 & 0\\
0 & 0 & 0 & 0 & \alpha & 0 & a\\
0 & 0 & 0 & 0 & 0 & \alpha & b\\
0 & 0 & 0 & 0 & 0 & 0 & \alpha\\
\end{matrix}\right) \eqno(*)
$$
over the ring $\mathbb{Z}$ of integers. In 3.6.8, it is proved that $R$ is a non-commutative centrally essential ring. 
Let $S_1$ be a semiring generated by matrices of the form $(*)$ over $\mathbb{Z}^{+}$ and scalar matrices with 
$\alpha\in \mathbb{Z}^{+}\cup \{0\}$ and zeros on the remaining positions. 
Since $Z(S_1)$ consists of scalar matrices, $S_1$ is not a centrally essential semiring. We note that $S_1$ is a semiring without zero-divisors. At the same time, the semiring $S_2$ of matrices of the form $(*)$ over the semiring $\mathbb{Z}^{+} \cup \{0\}$ is a centrally essential semiring.

\textbf{5.2.6. Proposition.}\\
Let $S$ be a centrally essential semiring without zero-divisors. If the ring $D(S)$ does not contain zero-divisors, the semiring $S$ is commutative.

$\lhd$ Let $0\neq a = x - y\in D(S)$. By assumption, $0\neq xc = d$ and $0\neq yf = g$ for some $c, d, f, g\in Z(S)$. Then
$$
a(cf) = (x - y)cf = (xc)f - (yf)c = df - gc.
$$
We need the following familiar property, \cite[Chapter II, Theorem 5.13]{HebW98}: In any semiring $S$ with ring of differences $D(S)$, any central element of $S$ is contained in the center of $D(S)$.

Therefore,$c, d, f, g\in Z(D(S))$ and $ac'\in Z(D(S))$, where $c' = cf$. In addition, $ac'\neq 0$, since $D(S)$ does not contain zero-divisors. Then $D(S)$ is a commutative ring by 1.2.2.~$\rhd$

\textbf{5.2.7. The upper central series and the nilpotence class.} \\
We recall that the \textsf{upper central series} of a group $G$ is the chain of subgroups 
$$
\{e\} = C_0(G)\subseteq C_1(G)\subseteq \ldots,
$$
where $C_i(G)/C_{i-1}(G)$ is the center of the group $G/C_{i-1}(G)$, $i\ge 1$. For the group $G$, the \textsf{nilpotence class} of $G$ is the least positive integer $n$ with $C_n(G) = G$ provided such an integer $n$ exists.

\textbf{5.2.8. Proposition; cf. Proposition 3.2.4.}\\ 
Let $G$ be a finite the group of nilpotence class $n\le 2$ and let $S$ be a commutative semiring without zero-divisors or zero sums. Then $SG$ is a centrally essential group semiring.

$\lhd$ If $n = 1$, then the group $G$ is Abelian and $SG$ is a centrally essential group semiring; see 1.4.1(a).

Let $n = 2$. Similar to the case of the group rings (e.g., see, \cite[Part 2]{Pas77}), the center $Z(SG)$ is a free $S$-semimodule with basis
$$
\left\{\sum_{K} \, | \, K \mbox{are the conjugacy classes in the group $G$}\right\}.
$$
It is sufficient to verify that $SG\sum_{Z(G)}\subseteq Z(SG)$, where $Z(G)$ is the center of the group $G$. Indeed, if $g, h\in G$, then 
$$
(gh)^{-1}hg\sum_{Z(G)} = \sum_{Z(G)},
$$
since $h^{-1}g^{-1}hg\in G'\subseteq Z(G)$.~$\rhd$

In Example 5.2.9 below, a non-commutative centrally essential semiring without zero-divisors is constructed; this semiring is additively cancellative but is not semisubtractive.

\textbf{5.2.9. Example.}\\
Let $Q_8$ be the quaternion the group, i.e., the group with two generators $a$, $b$ and defining relations 
$a^4 = 1$, $a^2 = b^2$ and $aba^{-1} = b^{-1}$; e.g., see, \cite[Section 4.4]{Hal59}. We have 
$$
Q_8 = \{e, a, a^2, b, ab, a^3, a^2b, a^3b\},
$$
the conjugacy classes of $Q_8$ are
$$
K_e = \{e\}, K_{a^2} = \{a^2\}, K_a = \{a, a^3\}, K_b = \{b, a^2b\}, K_{ab} = \{ab, a^3b\},
$$
and the center $Z(Q_8)$ is $\{e, a^2\}$. We consider the group semiring $SQ_8$, where $S = \mathbb{Q}^{+}\cup \{0\}$. 
Since $Q_8$ is a group of nilpotence class $2$, it follows from Proposition 5.2.8 that $SQ_8$ is a centrally essential group semiring.
To illustrate the above, we have 
$$
a\sum_{Z(Q_8)} = \sum_{K_a},\quad b\sum_{Z(Q_8)} = \sum_{K_b},
$$
$$
ab\sum_{Z(Q_8)} = \sum_{K_{ab}},\quad a^3\sum_{Z(Q_8)} = \sum_{K_a},
$$ 
$$
a^2b\sum_{Z(Q_8)} = \sum_{K_b},\quad a^3b\sum_{Z(Q_8)} = \sum_{K_{ab}}.
$$
The the group ring of differences $\mathbb{Q}Q_8$ is a reduced ring; see \cite[Theorem 3.5]{Seh75}. Then $SQ_8$ is a reduced semiring. Indeed, if $x^2 + y^2 = xy + yx$ and $x\neq y$, then $x^2 + y^2 - xy - yx = (x - y)^2 = 0$ in the ring $\mathbb{Q}Q_8$; this is not true. Thus, $SQ_8$ is a non-commutative reduced centrally essential semiring without zero-divisors. 
We note that the ring $\mathbb{Q}Q_8$ is not centrally essential, since centrally essential reduced rings are commutative.~$\rhd$

\textbf{5.2.10. Theorem.}\\ 
There exists a non-commutative additively cancellative reduced centrally essential semiring without zero-divisors.
An additively cancellative reduced semiring $S$ is commutative if and only if the ring of differences of $S$ is a centrally essential ring.
 
$\lhd$ It follows from Example 5.2.9 that there exists a non-commutative additively cancellative reduced centrally essential semiring without zero-divisors.

If a semiring $S$ is commutative, then $D(S)$ is a commutative ring, i.e., $D(S)$ is centrally essential. Conversely, let $D(S)$ be a centrally essential ring. Since $S$ is a reduced semiring, $D(S)$ is a reduced ring. Indeed, let $0\neq a = x - y\in D(S)$.
If $a^2 = 0$, then $x^2 + y^2 = xy + yx$. Therefore, $x = y$, $a = 0$, and we have a contradiction. Then the ring $D(S)$ is commutative, since $D(S)$ is a reduced centrally essential ring. Consequently, $S$ is a commutative semiring.~$\rhd$

\textbf{5.2.11. Multiplicatively cancellative semirings.}\\
An element $x$ of a semiring is said to be \textsf{left} (resp., \textsf{right}) \textsf{multiplicatively cancellative}\label{mulCanEle} if $y = z$ for all $y, z\in S$ with $xy = xz$ (resp., $yx = zx$).
A semiring $S$ is said to be \textsf{left} (resp., \textsf{right}) \textsf{multiplicatively cancellative}\label{mulCanSem} if every element $x\in S\setminus \{0\}$ is left (resp., right) multiplicatively cancellative. A left and right multiplicatively cancellative semiring is said to be \textsf{multiplicatively cancellative};\label{mulCanSem} e.g., see, \cite[Chapter I]{HebW98}. 

\textbf{5.2.12. Remark.}\\
A left (resp., right) multiplicatively cancellative centrally essential semiring $S$ is commutative.

$\lhd$ Let $a$ and $b$ be two non-zero elements of the semiring $S$. Since $S$ is a centrally essential semiring, there exists an element $c\in Z(S)$ with $0\neq ac\in Z(S)$. A left multiplicatively cancellative semiring does not contain left zero-divisors; see \cite[Chapter I, Theorem 4.4]{HebW98}. Therefore, $acb\neq 0$. Then
$$
(ac)b = c(ab) = (ca)b = b(ca)= c(ba),
$$
whence we have $ab = ba$. A similar argument is true for right multiplicatively cancellative semirings.~$\rhd$

\textbf{5.2.13. Division semirings and semifields.}\\ 
A semiring with division, which is not a ring, is called a \textsf{division semiring}.\label{divSemi}\\
A commutative division semiring is said to be \textsf{semifield}.\label{semiF}

Any centrally essential division semiring is a semifield. Indeed, it follows from \cite[Chapter I, Theorem 5.5]{HebW98} that a division semiring with at lest two elements is multiplicatively cancellative. Therefore, our assertion follows from Remark 5.2.12.

\section[Non-Associative Rings]{Non-Associative Rings}\label{section6}

\textsf{In this section, the considered rings are not necessarily associative.}

We use notation and terminology from \cite{Zhe82}; also see \cite{Sch66}.

\subsection[Types of Centers and Central Essentiality]{Types of Centers and\\ Central Essentiality}\label{subsection6.1}

\textsf{In this subsection, we consider rings which are not necessarily unital or associative.}

Let $R$ be a ring. We denote by $R^1$ the union of $R$ with adjoint external unit.\label{R1}

The \textsf{associator}\label{ass-tor} of three elements $a,b,c$ of the ring $R$ is the element $(a,b,c)=(ab)c-a(bc)$ and the \textsf{commutator} of two elements $a,b\in R$ is the element $[a,b]=ab-ba$. 

\textbf{6.1.1. Associative and commutative centers.}\\
For a ring $R$, the \textsf{associative center},\label{assCen} \textsf{commutative center}\label{comCen} and \textsf{center}\label{CenNonass} of $R$ (in the sense of \cite[$\S$7.1]{Zhe82}) are the sets
$$
\begin{array}{l}
N(R)=\{x\in R:\forall a,b\in R,\;(x,a,b)=(a,x,b)=(a,b,x)=0\},\\
K(R)=\{x\in R:\forall a\in R,\;[x,a]=0\},\\
Z(R)=N(R)\cap K(R),
\end{array}
$$
respectively. It is clear that $N(R)$ and $Z(R)$ are subrings in $R$ and the ring $R$ is a unitary (left and right) $N(R)$-module and $Z(R)$-module.

\textbf{6.1.2. The centroid.} \\ 
For a ring $R$, we denote by $\widehat Z(R)$\label{centro} the \textsf{centroid} of $R$,\label{centro} i.e., the set of endomor\-phisms of the additive the group $(R,+)$, which commute with the left and right multiplications by elements of $R$. 

It is clear that $R$ can be considered as a left or right module over the associative commutative ring $Z(R)$; $R$ can be also considered as a unitary module over the unital associative commutative ring $Z(R)^1$ and as a unitary module over over the centroid $\widehat Z(R)$.

\textbf{6.1.3. Remark.}\label{centroid}\\
The associative center $N(R)$, the commutative center $K(R)$ and the center $Z(R)$ of the ring $R$ are $\widehat Z(R)$-submodules in $R$.

$\lhd$ Let $n\in N(R)$ and $c\in \widehat Z(R)$. For any $a,b\in R$, we have
$$
(c(n),a,b)=c(n)a\cdot b-c(n)\cdot
ab=c(na)b-c(n\cdot ab)=
$$
$$
=c((n,a,b))=0,
$$
$$
(a,c(n),b)=ac(n)\cdot b-a\cdot
c(n)b=c(an)b-ac(nb)=
$$
$$
=c((a,n,b))=0,
$$
$$
(a,b,c(n))=ab\cdot c(n)-a\cdot bc(n)=c(ab\cdot n)-c(a\cdot
bn)=
$$
$$
=c((a,b,n))=0.
$$
Consequently, $c(n)\in N(R)$.

Similarly, if $k\in K(R)$, then for any $a\in R$, we have
$$[c(k),a]=c(k)a-ac(k)=c(ka)-c(ak)=c([k,a])=0.
$$
Consequently, $c(k)\in K(R)$.

Finally, the assertion about the center $Z(R)$ directly follows from two previous assertions, since $Z(R)=N(R)\cap K(R)$ by definition.~$\rhd$

\textbf{6.1.4. Centrally essential, strongly centrally essential rings, weakly centrally essential, $N$-essential, and\\ $K$-essential rings.} \\ 
A ring $R$ with center $C=Z(R)$ is said to be \textsf{centrally essential}\label{Z-ess} if $Cr\cap C\neq 0$ for any non-zero element $r\in R$ (equivalently, $K\cap C\neq 0$ for any non-zero submodule $K$ of the module $R_C$, i.e., $C$ is an essential submodule of the module ${}_CR$).

A ring $R$ with the center $C=Z(R)$ is said to be \textsf{strongly centrally essential}\label{sCEring} (resp., \textsf{weakly centrally essential})\label{wCEring} if $Cr\cap C\neq 0$ (resp., $\widehat Z(R)r\cap C\neq 0$) for any non-zero element $r\in R$.

In the definition of a strongly centrally essential ring, we can formally replace $Z(R)$ by $N(R)$; in thise case, the ring $R$ is called a \textsf{left $N$-essential}\label{Ness} ring\label{Kess} ring).

A ring $R$ is said to be \textsf{left $K$-essential}\label{Kess} if $K(R)r\cap K(R)\neq 0$ for any non-zero element $r\in R$, i.e., $K=K(R)$ is an essential submodule of the module ${}_KR$.

The following proposition is well known in the associative case.

\textbf{6.1.5. Proposition.}\label{relations}\\
Let $R$ be a ring with the center $C=Z(R)$.

\textbf{a.} Any strongly centrally essential ring $R$ is centrally essential.

\textbf{b.} Any centrally essential ring $R$ is weakly centrally essential.

\textbf{c.} Any unital ring $R$ is strongly centrally essential if and only if $R$ is centrally essential, if and only if $R$ is weakly centrally essential.

$\lhd$ \textbf{a}. Since $C$ is a subring in $C^1$, we obtain the assertion.

\textbf{b.} It is sufficient to note that multiplications by central elements and multiplications by integers belong to the centroid of $R$.

\textbf{c.} It follows from \textbf{a} and \textbf{b} that it is sufficient to verify that if $R$ is a weakly essential ring with the identity element 1, then $R$ is strongly centrally essential. Let $R$ be weakly centrally essential and $r\in R\setminus \{0\}$. There exists an element $\widehat c\in \widehat Z(R)$ with $\widehat c(r)\in Z(R)\setminus \{0\}$. Then 
$$
0\neq \widehat c(r)=\widehat c(1\cdot r)=\widehat c(1)r\in Z(R)r,
$$
since $\widehat c(1)\in Z(R)$ by Remark 6.1.3.\label{centroid}
Thus, $Z(R)r\cap Z(R)\neq 0$. Therefore, $R$ is strongly centrally essential.~$\rhd$

We give Examples 6.1.6 and 6.1.8 which show that the classes of strongly centrally essential, centrally essential and weakly centrally essential rings are distinct in the general case.

\textbf{6.1.6. Example.}\\ 
Any non-zero ring $R$ with zero multiplication is a centrally
essential ring which is not strongly centrally essential. Indeed, $R=Z(R)$ and for any non-zero element $r\in R$, we have $r\in R^1r\cap R$ but $Z(R)r=0$.

For the next example, we need Remark 6.1.7.

\textbf{6.1.7. Remark.}\label{r-cube}\\
Let $R$ be a ring such that $R\cdot R^2=R^2\cdot R=0$ and let $\varphi\colon R\rightarrow R$ be an endomor\-phism of the group $(R,+)$ such that 
$$
\varphi(R)\subseteq R^2\subseteq \text{Ker}\,\varphi.
$$
Then $\varphi\in \widehat Z(R)$. Indeed, $\varphi(ab)=0$ for any two elements $a,b\in R$, since $ab\in R^2$; we also have
$a\varphi(b)=\varphi(a)b=0$, since $aR^2=R^2b=0$.

\textbf{6.1.8. Example.}\\ 
Let $F=\mathbb{Z}/3\mathbb{Z}$ be the field of order $3$, $\Lambda(F^2)$ be the Grassmann algebra of the two-dimensional linear space over $F$. Let $e_1,e_2$ be a basis of the space $F^2$ and let $R$ be the subalgebra of the algebra $\Lambda(F^2)$ with basis $e_1,e_2,e_1\wedge e_2$. Let $r=\alpha_1e_1+\alpha_2e_2+\alpha_3e_1\wedge e_2$ be an arbitrary element of the ring $R$. It is easy to see that $r\in Z(R)$ if and only if $\alpha_1=\alpha_2=0$, i.e., $Z(R)=R^2$ and $Z(R)^1r=Fr$ for any
$r\in R$. In particular, $Z(R)^1e_1=Fe_1$ and $Fe_1\cap R^2=0$; therefore, the ring $R$ 
is not centrally essential. Now let $r\neq 0$. If $r\in Z(R)$, then
$r\in \widehat Z r\cap Z(R)$, since $\widehat Z(R)$ contains the identity automor\-phism
of the group $(R,+)$. Let $\pi\colon R\rightarrow R/R^2$ be the canonical homomor\-phism. If $r\not\in Z(R)$, then $\pi(r)$ is a non-zero element of the two-dimensional space $R/R^2$ and there exists a linear mapping $\psi\colon R/R^2\rightarrow R^2$, for which $\psi(\pi(r))\neq 0$. If $\varphi = \psi\pi$, then $\varphi\in \widehat Z(R)$ by Remark 6.1.7 and $0\neq \varphi(r)\in \widehat Z(R)r\cap R^2=Z(R)$.
Consequently, $R$ is a weakly centrally
essential ring.~$\rhd$

\subsection{Reduced and Semiprime Rings}\label{subsection6.2}

A ring is said to be \textsf{reduced}\label{reNonass} if it does not contain non-zero elements with zero square. We note that associative reduced rings are exactly the rings without non-zero nilpotent elements.

A ring $R$ is said to be \textsf{semiprime}\label{semipNonass} if $R$ does not contain a non-zero ideal with zero multiplication; see \cite[\S 8.2]{Zhe82}.

\textbf{6.2.1. Theorem.} \label{redcenter1} \\
Let $R$ be a weakly centrally essential ring such that its center $C=Z(R)$ is a reduced ring.

\textbf{a.}\label{redcenter1} $R$ is a strongly centrally essential ring.

\textbf{b.}\label{redcenter2} $R$ is an associative ring.

\textbf{c.} $R$ is a commutative ring.

$\lhd$ \textbf{a.} Let $r\in R\setminus \{0\}$, $\varphi\in\widehat{C}$ and $\varphi(r)=d\in C\setminus\{0\}$.
Then $0\neq d^2 =d\varphi(r)=\varphi(dr)=\varphi(d)r$. By Remark 6.1.3,\label{centroid} $\varphi(d)\in C$. It is also clear that $d^2\in C$. Consequently, $0\neq \varphi(d)r\in Cr\cap C$, i.e., $R$ is a strongly centrally essential ring.

\textbf{b.} By \textbf{a}, $R$ is a strongly centrally essential ring. We assume that the ring $R$ is not associative and some elements $x,y,z$ of the ring $R$ have the non-zero associator $(x,y,z)=(xy)z-x(yz)$. Then there exist two elements $c,d\in C$ such that 
$$
d=(x,y,z)c\in C\setminus \{0\}.
$$
We note that $xd\neq 0$; otherwise, $d^2=(x,y,z)c\cdot d=$
$$
=(x,y,z)\cdot cd=(x,y,z)\cdot dc=((xy\cdot z)d - (x\cdot yz)d)c=
$$
$$
=(d(xy\cdot z)-d(x\cdot yz))c=((dx\cdot y)z-dx\cdot
yz)c=0,
$$ 
which is impossible. Therefore, there exists an element $b\in C$ such that $xd\cdot b=x\cdot db\in C\setminus \{0\}$.
We consider the set $I=\{c\in C: cx\in C\}$. It is clear that $db\in I$. Now we assume that $dI = 0$. Then 
$$
d(db)=0,\, (db)^2=db\cdot db=(db\cdot d)b=(d\cdot db)b=0,\, db=0,
$$
a contradiction. Therefore, $di\neq 0$ for some $i\in I$. However, 
$$
di=(xy\cdot z-x\cdot yz)c\cdot i=c((xi\cdot y)z - xi\cdot yz) = 0;
$$
this is a contradiction. Thus, $R$ is an associative ring.

\textbf{c.} We assume that the ring $R$ is not commutative and we have two elements $x,y\in R$ with $xy-yx\neq 0$. Then there exist two elements $c,d\in C$ such that $d=(xy-yx)c\in C\setminus \{0\}$. We note that $xd\neq 0$; otherwise, 
$$
d^2=(xy-yx)cd=c((xd)y-y(xd))=0;
$$
this is impossible. Therefore, there exists an element $z\in C$ such that $xdz\in C\setminus \{0\}$. We consider the set $I=\{c\in C\,|\, cx\in C\}$. It is clear that $dz\in I$. Now we assume that $dI = 0$. Then 
$$
d(dz)=0,\; (dz)^2=0,\; dz=0,
$$
a contradiction. Therefore, $di\neq 0$ for some $i\in I$. However, 
$$
di=(xy-yx)ci=c((xi)y - y(xi)) = 0;
$$
this is a contradiction. Thus, $R$ is a commutative ring.~$\rhd$

\textbf{6.2.2. Remark.}\\ 
It is clear that the center of a semiprime ring is a reduced ring; the converse is not always true since the ring of upper triangular matrices over a field is a non-semiprime ring with reduced center.

\textbf{6.2.3. Remark.}\\ 
Let $R$ be a centrally essential associative ring. In 1.1.4 and 1.2.5, it is proved that all idempotents of the ring $R$ are central and the ring $R$ is commutative if $R$ is semiprime. In the introduction, examples of finite non-commutative centrally essential associative unital rings are constructed.

\textbf{6.2.4. Alternative rings.}\\ 
A ring $R$ is said to be \textsf{right alternative}\label{altRing} (resp., \textsf{left alternative}) if $(ab)b=a(bb)$ (resp., $(aa)b=a(ab)$) for any elements $a,b\in R$.

Right and left alternative rings are called \textsf{alternative} rings. A ring $R$ is alternative if and only if $(a,a,b)=(a,b,b)$ for any elements $a,b\in R$, where $(a,b,c)$ denotes associator $(a,b,c)=(ab)c-a(bc)$ of elements $a,b,c$ of the ring $R$.

By the Artin theorem \cite[Theorem 2.3.2]{Zhe82}, the ring $R$ is alternative if and only if any two elements of $R$ generate the associative subring.

In connection to Remark 6.2.3, we prove Theorem 6.2.5.

\textbf{6.2.5. Theorem.}\\
Let $R$ be a centrally essential ring.

\textbf{a.} If the center $Z(R)$ of the ring $R$ is semiprime, then the ring $R$ is commutative and associative.

\textbf{b.} If the ring $R$ is alternative and $e$ is an idempotent of the ring $R$, then $e\in Z(R)$.

$\lhd$ \textbf{a.} It follows from 6.2.1 and 6.2.2 that any weakly centrally essential semiprime ring is associative and commutative.

\textbf{b.} If $R$ is a weakly centrally essential alternative ring and $e$ is an idempotent of the ring $R$, then we have to prove that $e\in Z(R)$.\\
Let $r$ be an arbitrary element of $R$. Further, we use several times the associativity of the subring generated by two elements $e$ and $r$ in the alternative ring $R$.
If $c\in \widehat Z(R)$ is an element of the centroid of the ring $R$ such that $c(ere-re)=d\in Z(R)$, then 
$$
de=c(ere-re)e=c((ere-re)e)=c(ere-re)=d.
$$
On the other hand, 
$$
ed=ec(ere-re)=c(e(ere-re))=c(0)=0.
$$
Therefore, $d=0$ and $\widehat Z(R)(ere-re)\cap Z(R)=0$. Since the ring $R$ is weakly centrally essential, we have $ere-re=0$. It can be similarly verified that $ere-er=0$, whence $re=er$.~$\rhd$

\textbf{6.2.6. Open questions.}\\
\textbf{a.} Is it true that any $N$-essential\footnote{See 6.1.4.} ring is associative? Our assumption: It is not true.

\textbf{b.} Is it true that any semiprime $N$-essential ring is associative?

The following example shows that analogous questions for $K$-essential\footnote{See 6.1.4.} rings have negative answers.

\textbf{6.2.7. Example.}\\ Let $F$ be an arbitrary field and let $R$ be the algebra over $F$ with basis $\{e,f,x_1,y_1,x_2,y_2,\ldots\}$ such that and its multiplication is defined on basis elements by the relations
$$
\begin{array}{l}
e^2=e,\;f^2=f,\;ef=x_1,\;fe=y_1,\\
ex_i=x_ie=x_i,\;y_jf=fy_j=y_j,\\
x_ix_j=x_{i+j},\;y_iy_j=y_{i+j},\\
x_if=fx_i=y_je=ey_j=x_iy_j=y_jx_i=0\text{\
for all\ }i,j\in\mathbb{N}.
\end{array}
$$
We set $x=x_1$, $y=y_1$. It is easy to verify that $K(R)=F[x]x+F[y]y$. Indeed, if $r=ae+bf+s$, where $a,b\in F$ and $s\in F[x]x+F[y]y$, then $[r,e]=b(y-x)$ and $[r,f]=a(x-y)$. Therefore, $K(R)\subseteq F[x]x+F[y]y$; the converse inclusion follows from the definition of the multiplication in $R$.

Now we note that $K(R)$ is an ideal in $R$ and the ring $K(R)\cong F[x]x\oplus F[y]y$ is reduced. Therefore, if $r\in K(R)\setminus\{0\}$, then $r^2\in K(R)r\setminus\{0\}$. If $r\notin K(R)$, then 
$$
r=ae+bf+\sum_{i=1}^\infty a_ix^i+\sum_{i=1}^\infty b_iy^i,
$$
where $a,b,a_i,b_i\in F$ for all $i\in\mathbb{N}$ and $a\neq 0$ or $b\neq 0$. If $a\neq 0$, then $xr=ax+\sum_{i=1}^\infty a_ix^{i+1}\in (K(R)r\cap K(R))\setminus\{0\}$; similarly, if $b\neq 0$, then 
$$
yr=by+\sum_{i=1}^\infty b_iy^{i+1}\in (K(R)r\cap K(R))\setminus\{0\}.
$$
Thus, $R$ is $K$-essential.
We have that $R/K(R)\cong F\oplus F$ and $K(R)$ are associative reduced rings. Therefore, $r=0$ for any element $r\in R$ with $r^2=0$. Especially, $R$ does not have a non-zero ideal with zero multiplication, i.e., $R$ is semiprime. In the same time, $e$ and $f$ are non-central idempotents of $R$; this is impossible in any associative semiprime weakly centrally essential ring.

\textbf{6.2.8. Remark.}\\ If $R$ is an alternative ring without elements of order 3 in the additive the group, then $R$ is $K$-essential if and only if $R$ is centrally essential, since $3K(R)\subseteq N(R)$ \cite[Corollary 7.1.1]{Zhe82}.

\subsection[Cayley-Dickson Process and Associative Centers]{Cayley-Dickson Process and\\ Associative Centers}\label{subsection6.3}


Let $R$ be a ring, $M$ a left $R$-module, and $S$ a subset of $M$. Then $\text{Ann}_R(S)$ denotes the annihilator of $S$ in $R$, i.e., $\text{Ann}_R(S)=\{r\in R\,|\, rS=0\}$.
We denote by $[A,A]$ the ideal\label{[AA]} of the ring $A$ generated by commutators of all elements in $A$.

The following definition slightly generalizes the definition of the Cayley-Dickson process given in \cite[\S 2.2]{Zhe82}, see \cite{Alb42}.

\textbf{6.3.1. The Cayley-Dickson process\label{CDpro} and the rings $(A,\alpha)$.}\label{pCD}\\
Let $A$ be a ring with involution $*$\footnote{We recall that the ring anti-endomor\-phism is called an \textsf{involution}\label{inv} if its double application is the identity mapping.} and $\alpha$ an invertible symmetrical element of the center of the ring $A$.
We define a multiplication operation on the Abelian group $A\oplus A$ as follows:
$$\label{mCD}
(a_1,a_2)(a_3,a_4) = (a_1a_3 +\alpha
a_4a_2^*,\; a_1^*a_4 + a_3a_2) \eqno (*)
$$
for any $a_1,\ldots,a_4\in A$. We denote the obtained ring by $(A,\alpha)$.\label{Aalpha}

The subset in $(A,\alpha)$ consisting of elements in $(A,\alpha)$ of the form $(a,0)$ ($a\in A$) is a subring in $(A,\alpha)$ which is isomorphic to the ring $A$; we identify the elements with the corresponding elements of $A$. We set
$\nu=(0,1)\in (A,\alpha)$. Then $a*\nu=(0,a)=\nu a$ for any
$a\in A$ and $\nu^2=\alpha$. Thus, $(A,\alpha)=A+A\nu$.

Note that many works were devoted to the study of the structure and the properties of rings and algebras obtained by this process, for instance, \cite{Bro67}, \cite{Fla06}, \cite{FlaSh13}, \cite{FlaSt09}, \cite{PumA06}, \cite{Sch66}, \cite{Wat87}.

The following properties are directly verified with the use of the above relation $(*)$.\label{elprop}

\textbf{a.} $\nu^2=\alpha$ and $\nu a=a^*\nu$ for any element $a\in A$.

\textbf{b.} $(1,0)$ is the identity element of the ring $(A,\alpha)$.

\textbf{c.} The set $\{(a,0)\,|\, a\in A\}$ is a subring of the ring $(A,\alpha)$, which is isomorphic to the ring $A$.

\textbf{d.} The mapping $(a,b)\mapsto (a^*,-b)$, $a,b\in A$, is an involution of the ring $(A,\alpha)$.

Up to the end of Section 6.3, we fix a ring $A$ and an element $\alpha$, which satisfy the conditions of the Cayley-Dickson process from 6.3.1 \label{pCD}; we also set $R=(A,\alpha)$.

\textbf{6.3.2. Lemma.}\label{ident}\\
An element $(x,y)\in R$ belongs to the ring $N(R)$ if and only if for any elements $u,v\in A$, the following relation systems hold
$$\label{id_x}
\begin{array}{l}
(xu)v=x(uv),(ux)v=u(xv),(uv)x=u(vx),\\
v(ux)=x(vu),(xu)v=u(vx),(vu)x=(xv)u,\\
v(xu)=(vu)x,v(ux)=(vx)u,x(uv)=u(xv),\\
(ux)v=(uv)x,v(xu)=(xv)u,x(vu)=(vx)u;
\end{array} \eqno (*)
$$
$$\label{id_y}
\begin{array}{l}
(uy)v=y(vu),(uy)v=(yv)u,y(vu)=u(yv),\\
v(yu)=y(uv),(yu)v=(vy)u,y(uv)=(vy)u,\\
v(uy)=(uv)y,v(uy)=u(vy),(vu)y=u(vy),\\
(yu)v=(vu)y,v(yu)=u(yv),(uv)y=(yv)u.
\end{array} 
\eqno (**)
$$

$\lhd$ Let $(x,y)\in R$. Since associators are linear, \mbox{$(x,y)\in N(R)$} if and only if for any two elements $u,v\in A$, we have
$$\label{associators}
((x,y)(u,0)(v,0))=((u,0),(x,y),(v,0))=
$$
$$
=((u,0),(v,0),(x,y))=0,
$$
$$
((x,y)(u,0)(0,v))=((u,0),(x,y),(0,v))=
$$
$$
=((u,0),(0,v),(x,y))=0,\eqno (***)
$$
$$
((x,y)(0,u)(v,0))=((0,u),(x,y),(v,0))=
$$
$$
=((0,u),(v,0),(x,y))=0,
$$
$$
((x,y)(0,u)(0,v))=((0,u),(x,y),(0,v))=
$$
$$
=((0,u),(0,v),(x,y))=0.
$$
By calculating associators from $(***)$, we obtain the following system consisting of 12 relations
$$
\begin{array}{l}
((xu)v,v(uy))=(x(uv),(uv)y),\\
((ux)v,v(u^*y))=(u(xv),u^*(vy)),\\
((uv)x,(v^*u^*)y)=(u(vx),u^*(v^*y)),\\
({\alpha}v(y^*u^*),(u^*x^*)v)=({\alpha}(u^*v)y^*,x^*(u^*v)),\\
({\alpha}v(y^*u),(x^*u^*)v)=({\alpha}u(vy^*),u^*(x^*v)),\\
({\alpha}y(v^*u),x(u^*v))=({\alpha}u(yv^*),u^*(xv)),\\
({\alpha}(uy^*)v,v(x^*u))=({\alpha}(vu)y^*,x^*(vu)),\\
({\alpha}(yu^*)v,v(xu))=({\alpha}(vy)u^*,(xv)u),\\
({\alpha}y(u^*v^*),x(vu))=({\alpha}(v^*y)u^*,(vx)u),\\
({\alpha}v(u^*x),{\alpha}(yu^*)v)=({\alpha}x(vu^*),{\alpha}(vu^*)y),\\
({\alpha}v(u^*x^*),{\alpha}(uy^*)v)=({\alpha}(x^*v)u^*,{\alpha}(vy^*)u),\\
({\alpha}(vu^*)x,{\alpha}(uv^*)y)=({\alpha}(xv)u^*,{\alpha}(yv^*)u).
\end{array}
$$
By equating components of equal elements of the ring $R$ and considering that the element $\alpha$ is invertible, we obtain the following equivalent system
$$
(xu)v=x(uv),v(uy)=(uv)y),
(ux)v=u(xv),v(u^*y)=u^*(vy),
$$
$$
(uv)x=u(vx),(v^*u^*)y=u^*(v^*y),
 v(y^*u^*)=(u^*v)y^*,
 $$
$$
(u^*x^*)v=x^*(u^*v), v(y^*u)=u(vy^*),(x^*u^*)v=u^*(x^*v),
$$
$$
y(v^*u)=u(yv^*),u^*(xv)=x(u^*v), (uy^*)v=(vu)y^*,
$$
$$
v(x^*u)=x^*(vu), (yu^*)v=(vy)u^*,
v(xu)=(xv)u, 
$$
$$
y(u^*v^*)=(v^*y)u^*,x(vu)=(vx)u, v(u^*x)=x(vu^*),
$$
$$
(yu^*)v=(vu^*)y, v(u^*x^*)=(x^*v)u^*,(uy^*)v=(vy^*)u,
$$
$$
(vu^*)x=(xv)u^*,(uv^*)y=(yv^*)u.
$$
We replace the equations, both parts of which contain $x^*$ or $y^*$, by relations of
conjugate elements. We note that either $u$ or $u^*$ stands in every equation. Therefore, we can put $u$ instead of $u^*$, since $A^*=A$. Similarly, we replace $v^*$ by $v$. By choosing equations containing $x$, we obtain $(*)$,
the remaining equations form the system $(**)$.~$\rhd$

\textbf{6.3.3. Lemma.}\label{x_in_Z}\\
Let $x\in A$. The relations $(6.3.2(*))$ hold for all $u,v\in A$ if and only if $x\in Z(A)$.

$\lhd$ Let $x\in A$ and let relations~$(6.3.2(*))$ hold for all $u,v\in A$. The first three relations mean that $x\in N(A)$. It follows from the fourth relation for $u=1$ that $x\in K(A)$. Consequently, $x\in Z(A)$.

Conversely, if $x\in Z(A)$, then each of the relations in $(6.3.2(*))$ is transformed into one of the true relations $x(uv)=x(uv)$ or $x(vu)=x(vu)$, i.e., relations
$(6.3.2(*))$ hold for all $u,v\in A$.~$\rhd$

\textbf{6.3.4. Lemma.}\label{y_in_Ann_I}\\
Let $y\in A$. The relations $(6.3.2(**))$
hold for all $u,v\in A$ if and only if $y\in \text{Ann}_{Z(A)}([A,A])$.

$\lhd$ Let $y\in A$ and let relations $(6.3.2(**))$ hold for all $u,v\in A$. First of all, we note that for $v=1$, the first equation of $(6.3.2(**))$
turns into equation $uy=yu$; this is equivalent to the inclusion $y\in K(A)$, since the element $u$ of $A$ is arbitrary.

We verify that $y\in N(A)$. For any two elements $u,v\in A$, we have
$$
\begin{array}{l}
(yu)v\stackrel{\overline 1}{=}(vy)u\stackrel{\overline 2}{=}y(uv),\\
(uy)v\stackrel{1}{=}y(vu)\stackrel{3}{=}u(yv)\\
(uv)y=y(uv)\stackrel{4}{=}v(yu)=v(uy)\stackrel{8}{=}u(vy).
\end{array}
$$
In these transformations, the number over the relation sign is the number of used equation in $(6.3.2(**))$ (equations are numbered in rows from the left to right beginning with the first row). The number underscore denotes that, instead of the given equation, we use the equivalent equation obtained by the permutation of the variables $u,v$.

Consequently, $y\in N(A)\cap K(A)=Z(A)$.

Finally, we take into account the proven to see that already the first equation of $(6.3.2(**))$ implies that $y[u,v]=0$ for any $u,v\in A$, i.e., $y\in\text{Ann}_C([A,A])$.

Conversely, if $y\in \text{Ann}_{Z(A)}([A,A])$, then each of the relations $(6.3.2(**))$ is transformed into the true relation $y(uv)=y(vu)$, i.e., relations
$(6.3.2(**))$ hold for all $u,v\in A$.~$\rhd$

\textbf{6.3.5. Theorem.}\\
Let $A$ be a ring with center $C=Z(A)$, $I=\text{Ann}_C([A,A])$, $R=(A,\alpha)$. Then
 $N(R)=\{(x,y)\colon x\in C,\;y\in I\}$.

$\lhd$ The assertion follows from Lemma 6.3.2, Lemma 6.3.3 and Lemma 6.3.4.~$\rhd$

\textbf{6.3.6. Remark.}\\
Theorem 6.3.5 implies the following classical result (cf. \cite[Exercise 2.2.2(a)]{Zhe82}): the ring $R=(A,\alpha)$ is associative if and only if the ring $A$ is associative and commutative.

\subsection[Cayley-Dickson Process and Central Essentiality]{Cayley-Dickson Process and\\ Central Essentiality}\label{subsection6.4}

\textbf{6.4.1. Lemma.}\label{ess_in_R}\\ Let $B$ be a subring of the center of the ring $A$ and $I$ an essential ideal of $B$. If $B$ is an essential $B$-submodule of the module ${}_BA$, then $I$ is an essential $B$-submodule of the module ${}_BR$.

$\lhd$ If $r$ is a non-zero element of the ring $R$, then there exists an element $b\in B$ with $0\neq br\in B$. Therefore, there exists an element $d\in B$ such that $0\neq dcr\in I$ and $Br\cap I\neq 0$.~$\rhd$

\textbf{6.4.2. Theorem.}\\
Let $A$ be a ring with center $C=Z(A)$, $I=\text{Ann}_C([A,A])$, $R=(A,\alpha)$.
The ring $R$ is left (resp., right) $N$-essential if and only if
$A$ is centrally essential and $I$ is an essential ideal in $C$.

$\lhd$ Let the ring $A$ and the element $\alpha$ satisfy the conditions of 6.3.1 (the Cayley-Dickson process).\label{pCD} It is obvious that $C^*=C$, $I^*=I$, $\alpha C=C$ and $\alpha I=I$.

Let the ring $R=(A,\alpha)$ be $N$-essential. Then for any non-zero element $a\in A$, there exists an element $(x,y)\in N(R)$ such that $(x,y)(a,0)=(xa,ay)\in N(R)\setminus\{0\}$. By Theorem 6.3.5\label{main}, $x\in C$ and $y\in I$. If $xa\neq 0$, then $xa\in C\setminus\{0\}$; otherwise, $ya\in C\setminus\{0\}$. In the both cases, $Ca\cap C\neq 0$. Thus, $A$ is centrally essential.

We prove that $I$ is an essential ideal of $C$. Let $c\in C\setminus \{0\}$. If $Ic\neq 0$, then $Ic\subseteq I$ and $Cc\cap I\neq 0$. Let $Ic=0$. We consider the element $(0,c)$. There exists an element $(x,y)\in N(R)$ such that $(x,y)(0,c)=(\alpha cy,x^*c)\in N(R)\setminus\{0\}$. Since $\alpha y\in I$, $\alpha c y=0$, we have that $x^*c\neq 0$ and $x^*c\in I$ by Theorem 6.3.5. Consequently, $Cc\cap I\neq 0$, which is required.

Conversely, let's assume that $A$ is a centrally essential ring and $I$ is an essential ideal in $C$.

Let $(x,y)\in R\setminus\{0\}$. There exists an element $c\in C$ such that $cx\in C\setminus\{0\}$. Since $(c,0)\in N(R)$, we have 
$$
0\neq (c,0)(x,y)=(cx,c^*y)\in N(R)(x,y).
$$
If $c^*y=0$, then 
$$
0\neq (cx,0)\in N(R)(x,y)\cap N(R).
$$ If $c^*y\neq 0$, then by Lemma 6.4.1\label{ess_in_R} (for $B=C$), there exists an element $d\in C$ such that $dc^*y\in I\setminus \{0\}$. Then
$$
(d^*,0)(c,0)(x,y)=(d^*,0)(cx,c^*y)=
$$
$$
=(d^*cx,dc^*y)\in N(R)(x,y)\cap N(R)\setminus\{0\}.
$$
Thus, the ring $R$ is $N$-essential.~$\rhd$

\textbf{6.4.3. Remark.} \\
Up to the end of this subsection, we fix a ring $A$ with center, $C=Z(A)$, and an element $\alpha$, which satisfy 6.3.1 (the Cayley-Dickson process).\label{pCD} 

We set $R=(A,\alpha)$,
$$\label{centerdefs}
\begin{array}{l}I=\text{Ann}_C([A,A]),\quad
B=\{a\in C:\;a=a^*\},\\
J=\text{Ann}_B(\{a-a^*\;|\;a\in A\}).
\end{array}
$$
We note that the sets $B$ and $J$ are invariant with respect to the involution and are closed by the multiplication by $\alpha$.

\textbf{6.4.4. Theorem.}\label{main2.1}\\
$Z(R)=\{(x,y)\,|\, x\in B,\,y\in I\cap J\}$.

$\lhd$ Let $(x,y)\in Z(R)$. Since $Z(R)\subseteq N(R)$, it follows from
Theorem 6.3.5 \label{main} that $x\in C$ and $y\in I$. The relations $(0,1)(x,y)=(x,y)(0,1)$ imply the relations $\alpha y=\alpha y^*$ and $x=x^*$. Consequently, $x\in B$ and $y\in B\cap I$. Next, the relation $(a,0)(x,y)=(x,y)(a,0)$ ($a\in A$) implies the relations $ax=xa$ and $ay=a^*y$. The first relation holds for any $x\in C$ and the second relation means that $y(a-a^*)=0$, i.e., $y\in J$. Consequently, $y\in I\cap J$.

Conversely, if $x\in B$ and $y\in I\cap J$, then $(x,y)\in N(R)$ and for any $a,b\in A$
we have
\begin{equation*}
\begin{array}{l}
(x,y)(a,b)=(xa+\alpha by^*,
x^*b+ay)=(xa+\alpha yb, xb+ay),\\
(a,b)(x,y)=(ax+\alpha
yb^*,a^*y+xb)=(ax+\alpha yb,xb+ay)
\end{array}
\end{equation*}
Thus, $(x,y)\in K(R)$, whence $(x,y)\in Z(R)$.~$\rhd$

\textbf{6.4.5. Theorem.}\label{main2.2}\\
A ring $R=(A,\alpha)$ is centrally essential if and only if $B$ is an essential $B$-submodule of the ring $R$ and $J'=J\cap I$ is an essential ideal of the ring $B$.

$\lhd$ Let the ring $R=(A,\alpha)$ be centrally essential. Then for any $a\in A\setminus\{0\}$, there exists an element $(x,y)\in Z(R)$ such that $(x,y)(a,0)=(xa,ay)\in Z(R)\setminus\{0\}$. By Theorem 6.4.4\label{2.1}, $x,xa\in B$ and $y,ay=ya\in J'$. If $xa\neq 0$, then $xa\in B\setminus\{0\}$; otherwise, $ya\in B\setminus\{0\}$. In the both cases, we have $Ba\cap B\neq 0$. Thus, $B$ is an essential submodule of the module ${}_BA$.

We prove that $J'=$ is an essential ideal of the ring $B$. Let $b\in B\setminus \{0\}$. If $J'b\neq 0$, then $J'b\subseteq J'$ and $Bb\cap J'\supseteq J'b\cap J'\neq
0$. Let $J'b=0$. We consider element $(0,b)$. There exists an element $(x,y)\in Z(R)$ such that $(x,y)(0,b)=(\alpha by,x^*b)\in Z(R)\setminus\{0\}$. Since $\alpha
y\in J'$ and $\alpha b y=0$, we have that $x^*b\neq 0$, $x\in B$ and $x^*b=xb\in J'$, by Theorem 6.4.4\label{main2.1}. Consequently, $Bb\cap J'\neq 0$, which is required.

Conversely, let's assume that $B$ is an essential $B$-submodule of the ring $R$ and $J'$ is an essential ideal of the ring $B$.

Let $(x,y)\in R\setminus\{0\}$. First, we assume that $x\neq 0$. There exists an element $b\in B$ such that $bx\in B\setminus\{0\}$. Since $(b,0)\in Z(R)$, we have $0\neq
(b,0)(x,y)=(bx,b^*y)\in Z(R)(x,y)$. If $b^*y=0$, then $0\neq (bx,0)\in Z(R)(x,y)\cap Z(R)$. If $b^*y\neq 0$, then by Lemma 6.4.1\label{ess_in_R} (for $I=J'$) there exists an element $d\in B$ such that $db^*y\in J'\setminus \{0\}$. Then
$$
(d^*,0)(b,0)(x,y)=(d^*,0)(bx,b^*y)=
$$
$$
=(d^*bx,db^*y)\in Z(R)(x,y)\cap Z(R)\setminus\{0\}.
$$
Now let $x=0$. Then $y\neq 0$, and there exists an element $d\in B$ such that $dy\in J'\setminus \{0\}$. We obtain 
$$
(d^*,0)(0,b)=(0,db)\in Z(R)\setminus \{0\},\quad (d^*,0)\in Z(R).
$$ Thus, the ring $R$ is centrally essential.~$\rhd$

\subsection{Quaternion and Octonion Algebras}\label{quat-oct}\label{subsection6.5}

\textbf{6.5.1. Remarks and notation.}\\
Let $K$ be a commutative associative ring with the identity involution and $a$ an invertible element of the ring $R$. We consider the ring $A_1=(K,a)$. Then $A_1$ is a commutative associative ring, since $B=C=I=J=K$, under the notation of Section 6.4. It is natural to write elements of the ring $A_1$ in the form $x+yi$, where $x,y$ are elements of the ring $K$, $i=(0,1)$. On the ring $A_1$, an involution is defined by the relation $(x+yi)^*=x-yi$ for any $x,y\in K$. We choose an invertible element $b\in K$. Then $b$ is an invertible symmetrical element of the center of the ring $A_1$ and we can construct the ring $A_2=(A_1,b)$. We consider the $K$-basis of the algebra $A_2$, which is formed by the elements $1=(1,0)$, $i=(i,0)$, $j=(0,1)$ and $k=(0,-i)$. The relations 
$$
i^2=a,\, j^2=b,\, ij=-ji=k,\, ik=-ki=aj,\, kj=-jk=bi
$$
are directly verified. Consequently, the obtained ring is the generalized quaternion algebra $\left(\dfrac{a,b}{K}\right)$. It is well known (and also follows from Theorem 6.3.5\label{main}) that the ring $A_2$ is associative (e.g., see \cite[Example 7.2.III]{Zhe82}). The center of the ring $A_2$ is of the form $K+Ni+Nj+Nk$, where $N=\text{Ann}_K(2)$ (see \cite[Lemma 2(b)]{Tug93}). Let $B,I,J$ be defined by equations from $6.4.3$ for $A=A_2$. It is easy to verify that 
$$
B=C=Z(A_2),\; I=J=N+Ni+Nj+Nk.
$$

\textbf{6.5.2. Lemma.}\label{essK}
Under the above notation, the ideal $I$ is an essential ideal in $B$ if and only if $N$ is an essential ideal in $K$.

$\lhd$ Let $I$ be an essential ideal in $B$. If $x\in K\setminus\{0\}$, then there exists an element $y\in B$ such that $xy\in I\setminus\{0\}$. We set
$$
y=y_1+y_2i+y_3j+y_4k,\;\text{where } y_1\in K,\, y_2,y_3,y_4\in N.
$$ 
If $xy_1\neq 0$, then $xK\cap N\neq 0$. Otherwise, at least one of the elements $xy_2,xy_3,xy_4$ is not equal to $0$ and each of them belong to the ideal $N$, whence $xK\cap N\neq 0$ in this case too.

Conversely, if $N$ is an essential ideal in $K$ and 
$$
x=x_1+x_2i+x_3j+x_4k\in I\setminus\{0\},
$$
then $x_2,x_3,x_4\in N$. If $x_1\neq 0$, then there exists an
element $y\in K$ with $yx_1\in N\setminus\{0\}$. Then $yx\in Bx\cap
I\setminus\{0\}$. If $x_1=0$, then $x=1\cdot x\in Bx\cap I$. Thus,
$I$ is an essential ideal of the ring $B$.~$\rhd$

From the above argument, we obtain the following

\textbf{6.5.3. Proposition.}\label{csQuat}\\
The quaternion algebra $((K,a),b)$ is a non-commutative centrally essential ring if and only if $\text{Ann}_K(2)$ is a proper essential ideal of the ring $K$.

Now we consider an arbitrary invertible element $c\in K$ and the ring $A_3=(A_2,c)$. We set 
$$
f_1=i,\, f_2=j,\, f_3=k,\, f_4=l=(0,1),
$$
$$
f_5=(0,-i),\, f_6=(0,-j),\, f_8=(0,-k).
$$
It can be directly verified that the basis $\{1,f_1,f_2,\ldots,f_7\}$ of the $K$-module $A_3$ satisfies the relations from \cite{FlaSt09} for basis elements of the generalized octonion algebra ${\mathbb O}(\alpha, \beta, \gamma)$ (for $\alpha=-a, \beta=-b, \gamma=-c$).

Similar to Proposition 6.5.3\label{csQuat}, we obtain Proposition 6.5.4.

\textbf{6.5.4. Proposition.}\label{csOct}\\
The octonion algebra $(((K,a),b),c)$ is a non-associative centrally essential ring if and only if $\text{Ann}_K(2)$ is a proper essential ideal of the ring $K$.

\textbf{6.5.5. Example.}
Let $K=\mathbb{Z}_4$ and $R=(((K,1),1),1)$. We prove that $R$ is a non-associative non-commutative centrally essential ring.\\
Indeed, $\text{Ann}_K(2)=2K$ is an essential proper ideal in $K$. Therefore, the non-commutativity of the ring $((K,1),1)$ (and the non-commutativity of the ring $R$ containing $((K,1),1)$) follows from Proposition 6.5.3\label{csQuat} and the non-associativity of the ring $R$ follows from Proposition 6.5.4.\label{csOct}

We note that $R=(((K,1),1),1)$ is an alternative ring and the ring $(R,1)$ is not even right alternative, i.e., $(R,1)$ does not satisfy the identity $(x,y,y)=0$ \cite[Exercise 7.2.2]{Zhe82}. Thus, there are alternative non-associative finite centrally essential rings and non-alternative finite centrally essential rings.~$\rhd$

\textbf{6.5.6. Open Questions.}

\textbf{a.} Is it true that there exists a left $N$-essential ring which is not right $N$-essential?

\textbf{b.} Is it true that there exists a commutative $N$-essential (equivalently, centrally essential) non-associative ring?

\textbf{c.} Is it true that there exists a right alternative centrally essential or $N$-essential non-alternative ring?

\textbf{d.} How can we generalize the obtained results to the case of non-unital rings and the case, where the element $\alpha$ in definition 6.3.1 is not supposed to be invertible?

\textbf{e.} Since the Cayley-Dickson process gives non-associative division algebras (see, e.g., \cite{Bro67,Fla06}), it seems natural to state the following question:\\
What can be said about $N$-essentiality of these division rings? 

Note that centrally essential semiprime rings are commutative, but it is not known whether $N$-essential semiprime rings are associative.

\subsection*{Symbols}
\addcontentsline{toc}{subsection}{\textbf{Symbols}}



\vspace{0.8mm}

$\ell_R(S)$ and $\text{r}_R(S)$ \hfill left and right annihilators of the subset $S$ \hfill \pageref{ann}

\vspace{0.8mm}

$\text{Sing}\,M$ \hfill singular submodule of the module $M$ \hfill \pageref{sinsub}

\vspace{0.8mm}

$\Sigma_S$ \hfill element $\sum_{x\in S}x$ of the group ring $AG$ \hfill \pageref{sigmaS}

\vspace{0.8mm}

$\text{supp}\,(r)$ \hfill support $\{g\in G|a_g\neq 0\}$ of the element $r$ \hfill \pageref{supp}

\vspace{0.8mm}

$g^G$ \hfill class of conjugate elements containing the element $g$ \hfill \pageref{gG}

\vspace{0.8mm}

$\mathbb{Q}\text{End}\,A$ \hfill quasi-endomor\-phism ring of the group $A$ \hfill \pageref{quasEndRing}

\vspace{0.8mm}

$Z(S)$ or $C(S)$ \hfill center of a semiring or the group $S$ \hfill \pageref{center}

\vspace{0.8mm}

$(X\ld Y)$ \hfill subset $\{a\in A\,\vert \,Xa\subseteq Y\}$ of the right
$A$-module \hfill \pageref{ld}

\vspace{0.8mm}

$R^1$ \hfill union of the ring $R$ with adjoint external unit\hfill \pageref{R1}

\vspace{0.8mm}

$(a,b,c)=(ab)c-a(bc)$ \hfill associator of elements $a,b,c$ of the ring $R$\hfill \pageref{ass-tor}

\vspace{0.8mm}

$N(R)$ \hfill associative center of the not necessarily associative ring $R$\hfill \pageref{assCen} 

\vspace{0.8mm}

$K(R)$ \hfill commutative center of the not necessarily associative ring $R$\hfill \pageref{comCen}

\vspace{0.8mm}

$Z(R)=N(R)\cap K(R)$ \hfill center of the non-associative ring $R$\hfill \pageref{CenNonass}

\vspace{0.8mm}

$\widehat Z(R)$ \hfill centroid of the not necessarily associative ring $R$\hfill \pageref{centro}

\vspace{0.8mm}

$[A,A]$ \hfill ideal generated by commutators of all elements of $A$ \hfill \pageref{[AA]}

\vspace{0.8mm}

$(A,\alpha)$ \hfill ring from the Cayley-Dickson process \hfill \pageref{Aalpha}

\begin{theindex}
\addcontentsline{toc}{section}{\textbf{Index}}




\item additively cancellative semiring \hfill \pageref{addCan}

\item alternative ring \hfill \pageref{altRing}

\item annihilator \hfill \pageref{ann}

\item arithmetical ring \hfill \pageref{ariring}

\item Artinian module \hfill \pageref{artinmod}

\item associative center \hfill \pageref{assCen} 

\item associator \hfill \pageref{ass-tor}









\indexspace




\item Cayley-Dickson process \hfill \pageref{CDpro}


\item CE ring \hfill \pageref{cering}

\item center of a non-associative ring \hfill \pageref{CenNonass}

\item center of a semiring \hfill \pageref{cenSemir}

\item centrally essential non-associative ring \hfill \pageref{Z-ess}

\item centrally essential ring \hfill \pageref{ceness}

\item centrally essential semiring \hfill \pageref{cenEsSemi}

\item centroid \hfill \pageref{centro}


\item classical identity \hfill \pageref{cliden}



\item classical ring of fractions \hfill \pageref{clfra}





\item closed submodule \hfill \pageref{closed}


\item commutative center of a\\ non-associative ring \hfill \pageref{comCen}

\item commutator \hfill \pageref{comm}

\item complemented idempotent\hfill \pageref{compIde}








\indexspace

\item derivation \hfill \pageref{deriv}



\item distributive module \hfill \pageref{dismod}

\item divisible Abelian group \hfill \pageref{divGr}

\item division semiring \hfill \pageref{divSemi}

\item domain \hfill \pageref{dom}


\indexspace

\item element algebraic over the center \hfill \pageref{algcenele}






\item essential extension \hfill \pageref{essext}

\item essential submodule \hfill \pageref{esssub}




\indexspace



\item $FC$-group \hfill \pageref{FCgr}


\item finite-dimensional module \hfill \pageref{findimmod}





\item flat module \hfill \pageref{flat}



\indexspace

\item generalized anticommutative ring \hfill \pageref{genanti}





\item graded ring \hfill \pageref{grad}




\indexspace





\item hereditary module \hfill \pageref{hermod}



\item homogeneously faithful ring \hfill \pageref{homfai}


\indexspace



\item integrality over the center \hfill \pageref{celcenRing} 

\item invariant module \hfill \pageref{invMod}

\item invariant ring \hfill \pageref{invRing}



\item involution \hfill \pageref{inv}




\indexspace

\item $K$-essential ring \hfill \pageref{Kess}

\indexspace

\item large center \hfill \pageref{larcen}





\item local ring \hfill \pageref{locring}




\indexspace






\item multiplicatively cancellative\\ semiring \hfill \pageref{mulCanSem}

\indexspace

\item $N$-essential ring \hfill \pageref{Ness}

\item nil-index \hfill \pageref{nilInd}

\item nilpotence class \hfill \pageref{nilpCla}

\item nilpotence index \hfill \pageref{nilpInd}

\item Noetherian module \hfill \pageref{netermod}

\item non-reduced Abelian group \hfill \pageref{nonRed}

\item non-singular module \hfill \pageref{nonsin}


\item non-zero-divisor \hfill \pageref{nzd}





\indexspace


\item perfect ring \hfill \pageref{perf}

\item PI ring \hfill \pageref{pirin}

\item polynomial identity \hfill \pageref{poliden}

\item prime ring \hfill \pageref{priring}


\item projective module \hfill \pageref{promod}


\item pseudo-socle \hfill \pageref{pseSoc}

\item pure subgroup \hfill \pageref{purSub}


\indexspace

\item quasi-contained the group \hfill \pageref{quaCon}


\item quasi-decomposition of\\ an Abelian group \hfill \pageref{quaDec}

\item quasi-endomor\-phism \hfill \pageref{quasEnd}

\item quasi-endomor\-phism ring \hfill \pageref{quasEndRing}

\item quasi-equal groups \hfill \pageref{quaEq}




\indexspace


\item reduced Abelian group \hfill \pageref{redGr}

\item reduced semiring \hfill \pageref{redSemi}


\item reduced non-associative ring \hfill \pageref{reNonass}

\item regular element \hfill \pageref{regele}



\item ring algebraic over the center \hfill \pageref{algcenRing}


\item ring of differences \hfill \pageref{ringDif}

\indexspace

\item $S$-graded ring \hfill \pageref{sgrad}



\item semi-Artinian module \hfill \pageref{semiart}

\item semifield \hfill \pageref{semiF}


\item semihereditary module \hfill \pageref{semiher}

\item semilocal ring \hfill \pageref{semiloc}


\item semiperfect ring \hfill \pageref{semiperfring}


\item semiprime ring \hfill \pageref{semipr}

\item semiprime non-associative ring \hfill \pageref{semipNonass}

\item semiprime semiring \hfill \pageref{semisemipr}

\item semiring \hfill \pageref{semiRi}


\item semisubtractive semiring \hfill \pageref{semisubsSemi}




\item singular ideal \hfill \pageref{sinide}


\item singular submodule \hfill \pageref{sinsub}


\item socle \hfill \pageref{soc}


\item strongly centrally essential ring \hfill \pageref{sCEring}

\item strongly indecomposable group \hfill \pageref{strInd}



\item support \hfill \pageref{supp}

\indexspace

\item $T$-nilpotent ideal \hfill \pageref{tnilp}




\indexspace


\item uniform module \hfill \pageref{unifo}

\item uniserial module \hfill \pageref{unisermod}


\item upper central series \hfill \pageref{upcen}

\indexspace

\item von Neumann regular ring \hfill \pageref{vnrring}



\indexspace

\item weak global dimension \hfill \pageref{wgd1}

\item weakly centrally essential ring \hfill \pageref{wCEring}



\end{theindex}

\end{document}